\documentclass[11pt,a4paper]{article}
\usepackage{mathrsfs}
\usepackage[T1]{fontenc}
\usepackage{lmodern}
\usepackage{amssymb,latexsym,amsmath, amsthm}
\usepackage{amsmath,amssymb,amsthm,mathtools}
\usepackage[margin=27mm]{geometry}
\usepackage{booktabs,microtype,xcolor}
\usepackage[colorlinks=true,linkcolor=blue!50!black,citecolor=blue!50!black,urlcolor=blue!50!black]{hyperref}
\hypersetup{pdftitle={Forelli--Rudin-Type Operators in Infinite Dimensions, I},pdfauthor={Yu Hou and Jiayang Yu},pdfsubject={Weighted boundedness, Schatten classes, and trace formulas}}
\numberwithin{equation}{section}
\newtheorem{theorem}{Theorem}[section]
\newtheorem{proposition}[theorem]{Proposition}
\newtheorem{lemma}[theorem]{Lemma}
\newtheorem{corollary}[theorem]{Corollary}
\theoremstyle{remark}
\newtheorem{remark}[theorem]{Remark}
\newcommand{\C}{\mathbb C}
\newcommand{\R}{\mathbb R}
\newcommand{\N}{\mathbb N}
\newcommand{\dd}{\,\mathrm d}
\newcommand{\ip}[2]{\langle #1,#2\rangle}
\newcommand{\norm}[1]{\lVert #1\rVert}
\newcommand{\abs}[1]{\lvert #1\rvert}
\newcommand{\ind}{\mathbf 1}
\newcommand{\HH}{\ell^2}
\newcommand{\TT}{T_{a,b,c}}
\newcommand{\SSS}{S_{a,b,c}}
\newcommand{\Lt}[2]{L^{#1}_{#2}}

\DeclareMathOperator{\Ree}{Re}
\DeclareMathOperator{\Log}{Log}
\title{\textbf{Forelli--Rudin-Type Operators\\ in Infinite Dimensions, I}\\[5pt]
\large  }
\author{Yu Hou\thanks{School of Mathematics, Sichuan University,
Chengdu 610064, China. Email: \texttt{houyu020426@stu.scu.edu.cn}.}\and
Jiayang Yu\thanks{Corresponding author. School of Mathematics, Sichuan University,
Chengdu 610064, China. Email: \texttt{jiayangyu@scu.edu.cn}.}}
\date{}
\begin{document}
\maketitle

\begin{abstract}
We introduce Forelli--Rudin-type integral operators on the open unit ball of \(\ell^2\). We completely characterize the boundedness of these operators between weighted \(L^p\) and \(L^q\) spaces for all \(p,q\in[1,\infty]\), as well as their membership in every Schatten class of operators between the corresponding weighted \(L^2\) spaces. We also establish trace formulas for these operators that cannot be obtained from their finite-dimensional counterparts by a formal passage to the limit as the dimension tends to infinity. Our results are not merely infinite-dimensional extensions of those obtained by Cheng, Fang, Wang, and the second author (2017), Zhao and Zhou (2022), and Ding (2022). Rather, they reveal new phenomena that are unique to infinite dimensions and arise from the fundamental geometric and measure-theoretic differences between finite- and infinite-dimensional spaces. A key ingredient in our proofs is the establishment of new Forelli--Rudin-type asymptotic estimates in infinite dimensions. These estimates also provide new insights into the development of function theory on the infinite-dimensional unit ball.
\end{abstract}

\noindent\textbf{Keywords.} Gaussian measure; Forelli--Rudin-type
operator;  Forelli--Rudin-type asymptotic
estimates; $L^p-L^q$ boundedness; Schatten classes; Trace formulas.\\
 2020 MSC: 26E15, \, 47G10, \, 47B38, \,47B10, \, 46G99, \,47B34
\section{Introduction}

The Bergman projection and its associated integral operators play a central role in the theory of holomorphic function spaces and operator theory; see, for instance, \cite{AC,PR,Zhu05,Zhu}. A natural class of operators arising in this context consists of Forelli--Rudin-type operators, whose integral kernels are constructed from the reproducing kernels of Bergman spaces, together with suitable weight factors. The boundedness properties of these operators on Lebesgue and weighted Lebesgue spaces have been extensively investigated in finite-dimensional complex domains.

An important starting point for the present work is the study of the
hyper-singular cousins of the Bergman projection by Cheng, Fang, Wang, and
the second author \cite{CFWY}. In particular, on the unit disk \(\mathbb D\subset\mathbb C\), they
considered the operators
\[
    K_c f(z)
       \triangleq \int_{\mathbb D}
          \frac{f(w)}{(1-z\overline w)^c}\,\dd A(w),
       \qquad c\in\mathbb R,\,z\in\mathbb D,
\]
$\mathrm d A(w)$ denotes the normalized area measure on $\mathbb D$. They obtained a complete characterization of the pairs \((p,q)\in[1,\infty]^2\) for which
$
K_c:L^p(\mathbb D)\longrightarrow L^q(\mathbb D)
$
is bounded; see Table~\ref{tab:standard1}. This work demonstrates that the mapping properties of Bergman-type integral operators are strongly influenced by the geometry of the underlying domain and the structure of the associated reproducing kernel.
\begin{table}[htbp]
\centering
\renewcommand{\arraystretch}{1.6}
\begin{tabular}{@{}p{32mm}p{111mm}@{}}
\toprule
Range of \(c\)& Necessary and sufficient conditions\\
\midrule
\( c\leq 0\)
 &\(1\le p, q\leq\infty\)\\
\(0<c\leq 2\)
 &\(p=1 , 1\leq q<\frac{2}{c} ; 1<p<\frac{2}{2-c}, \text{and }\frac{1}{q}\geq \frac{1}{p}+\frac{c}{2}-1;p=\frac{2}{2-c},1\leq q<\infty;\frac{2}{2-c}<p\leq \infty, 1\leq q\leq \infty\)(Here $\frac{2}{0}\triangleq \infty$)\\
 \(2<c<3\)
 & \(\frac{1}{3-c}<p<\infty,\frac{1}{q}>\frac{1}{p}+c-2; p=\infty,1\leq q<\frac{1}{c-2}\)\\
\bottomrule
\end{tabular}
\caption{The common boundedness region for \(K_{c}\).
}
\label{tab:standard1}
\end{table}
\newpage

The finite-dimensional theory was subsequently developed on the unit
ball of \(\mathbb C^n\). For $a,b,c\in\mathbb R$, the Forelli--Rudin-type integral operators are defined by
\begin{align}
T^n_{a,b,c} f(\mathbf z)
 &\triangleq(1-||\mathbf z||^2)^a\int_{\mathbb B_n}
       \frac{(1-||\mathbf y||^2)^b f(\mathbf y)}{(1-\langle\mathbf z, \mathbf y\rangle)^c}\dd v_n(\mathbf y),\label{eq:standardTn}\\
S^n_{a,b,c} f(\mathbf z)
 &\triangleq(1-||\mathbf z||^2)^a\int_{\mathbb B_n}
       \frac{(1-||\mathbf y||^2)^b f(\mathbf y)}{|1-\langle\mathbf z, \mathbf y\rangle|^c}\dd v_n(\mathbf y),\label{eq:standardSn}
\end{align}
where $\mathbb B_n$ denotes the open unit ball in $\mathbb C^n$, and $\mathrm d v_n(\mathbf y)$ is the normalized volume measure on $\mathbb B_n$ such that $v_n(\mathbb B_n)=1$. For $1\leq p<\infty$ and $-1<\alpha<\infty$, we denote by $L^p(\mathbb B_n, v_{n,\alpha})$ the weighted Lebesgue space which contains measurable functions $f$ on $\mathbb B_n$ such that
$$
\int_{\mathbb B_n}|f |^p\,\mathrm d v_{n,\alpha} <\infty,
$$
where $d v_{n,\alpha}(\mathbf z)\triangleq \frac{\Gamma(n+\alpha+1)}{n!\Gamma(\alpha+1)}(1-||\mathbf z||^2)^{\alpha}\mathrm d v_n(\mathbf z)$ and the normalization constant is chosen so that $v_{n,\alpha}(\mathbb B_n)=1$.
We denote by $L^{\infty}(\mathbb B_n, v_{n,\alpha})$ the corresponding essentially bounded measurable functions.

Forelli and Rudin \cite{FR} initiated the study of weighted integral operators on the unit ball, and Kures and Zhu \cite{KZ} subsequently investigated a general class of Forelli--Rudin-type operators. More recently, Zhao and Zhou \cite{ZZ} obtained a complete characterization of the \(L^p(\mathbb B_n,v_{n,\alpha})\)-to-\(L^q(\mathbb B_n,v_{n,\beta})\) boundedness of two such classes of operators on the unit ball for all \(1\leq p,q\leq\infty\), including the endpoint cases (see Table~\ref{tab:standard2}). In a different direction, Schatten class properties of Bergman-type operators on finite-dimensional bounded symmetric domains have also been investigated; see, for instance, Ding \cite{Ding}. These studies rely, to varying degrees, on Forelli--Rudin-type asymptotic estimates. Obtaining sharp Forelli--Rudin-type asymptotic estimates in the generality needed for the study of integral operators on bounded symmetric domains can be quite subtle even in finite dimensions; see \cite{Ding,FK}.
\begin{table}[htbp]
\centering
\renewcommand{\arraystretch}{1.6}
\begin{tabular}{@{}p{32mm}p{111mm}@{}}
\toprule
Exponent & Necessary and sufficient conditions\\
\midrule
\(1<p \le q<\infty\)
 & \(p(b+1)>\alpha+1,\quad aq+\beta>-1,c\leq  n+1+a+b+\frac{n+1+\beta}{q}-\frac{n+1+\alpha}{p}\)\\
\(p=1\), \(1\le q<\infty\)
 & \(b>\alpha,\quad aq+\beta>-1,c\leq  a+b-\alpha+\frac{n+1+\beta}{q};b=\alpha,\quad aq+\beta>-1,c< a+ \frac{n+1+\beta}{q} \)\\
\(1\leq q<p<\infty\)
 & \(p(b+1)>\alpha+1,\quad aq+\beta>-1,c<n+1+a+b+\frac{1+\beta}{q}-\frac{1+\alpha}{p}\)\\
\(p=\infty\), \(1\le q<\infty\)
 & \(b>-1,\quad aq+\beta>-1, c<n+1+a+b+\frac{1+\beta}{q}\)\\
\(p=1\), \(q=\infty\)
 & \(a\ge0,\quad b\ge\alpha,\quad c\le a+b-\alpha\)\\
\(1<p<\infty\), \(q=\infty\)
 & \(a=0,\quad p(b+1)>\alpha+1,\,c<n+1+b-\frac{n+1+\alpha}{p};a>0,p(b+1)>\alpha+1,\,c\leq n+1+a+b-\frac{n+1+\alpha}{p}\)\\
\(p=\infty\), \(q=\infty\)
 & \(a=0,\quad b>-1,c<n+1+b;a>0,b>-1,c\leq n+1+a+b\)\\
\bottomrule
\end{tabular}
\caption{The common boundedness region for \(T_{a,b,c}^n\) and \(S^n_{a,b,c}\).
}
\label{tab:standard2}
\end{table}

The first objective of this paper is to develop an infinite-dimensional counterpart of the aforementioned theory. A major obstacle to achieving this goal is the lack of Forelli--Rudin-type asymptotic estimates in infinite dimensions. To overcome this difficulty, we take the open unit ball of \(\ell^2\), equipped with a Gaussian measure, as our underlying space. This choice is also motivated by our previous works \cite{WYZ1,YZ}, in which we established \(L^2\)-estimates for the \(\overline{\partial}\)-operator on this space and developed a fairly comprehensive and convenient analytical framework. Building on this framework, we successfully establish Forelli--Rudin-type asymptotic estimates in infinite dimensions. These estimates play a crucial role in the proofs of the main results of this paper.

The passage from a finite-dimensional unit ball to its infinite-dimensional counterpart is far from being a mere matter of replacing the dimension \(n\) by infinity. To illustrate the fundamental difficulties arising in infinite dimensions, we recall two relevant observations.

\begin{itemize}
\item Tao observed in \cite[p.~86]{Tao} that ``the advantage of the Baire category approach is that it also works well in infinite dimensional complete normed vector spaces, i.e., Banach spaces, whereas the measure-theoretic approach runs into significant difficulties in infinite dimensions'';

\item Lempert remarked in \cite[p.~775]{Lem99} that ``Why is the study of the \(\overline{\partial}\)-equation in infinite dimensions so much harder than in \(\mathbb C^n\)? We can discern two reasons. One has to do with the fact that in \(\mathbb C^n\) solutions of the \(\overline{\partial}\)-equation are constructed using integrals, and most of the time integrals with respect to volume measure; here translation invariance of the Lebesgue measure is essential. In infinite dimensions translation invariant measures are all pathological, hence tools different from integrals must be found.''
\end{itemize}

In a nutshell, there is no nontrivial translation-invariant Borel measure in infinite-dimensional Banach spaces that plays the role of Lebesgue measure in finite dimensions. This fundamental obstruction is closely related to the striking geometric differences between finite- and infinite-dimensional normed spaces. Indeed, every nonempty open ball in an infinite-dimensional normed space contains infinitely many pairwise disjoint nonempty open balls of the same radius. Consequently, the geometry and measure theory of domains in infinite dimensions must be approached in a fundamentally different way. For a more detailed introduction, we refer the reader to \cite{YZ-b}, lecture notes based on our research and teaching experience in infinite-dimensional analysis at Sichuan University since 2016.

Another reason for choosing the open unit ball of \(\ell^2\), equipped with a Gaussian measure, as our underlying space is that this framework provides a more transparent and convenient connection with finite-dimensional analysis than many existing approaches to infinite-dimensional analysis, including abstract Wiener spaces, white noise analysis, and Malliavin calculus. It is particularly well suited to explicit computations and the investigation of various analytical problems in infinite-dimensional spaces.

The second objective of this paper is to develop function theory on the infinite-dimensional unit ball. In comparison with the classical theory of holomorphic functions on the unit ball of \(\mathbb C^n\), as presented in \cite{Rud}, several remarkable results have already been established in infinite dimensions. For instance, the work \cite{IS} provides a comprehensive infinite-dimensional counterpart of Chapter 2 of \cite{Rud}.

A fundamental tool in Rudin's approach is the Cauchy integral formula and its higher-dimensional analogue, namely Koppelman's Cauchy formula (see \cite[p.~347]{Rud}; for a more general version, see also \cite[p.~27]{Dem}), which we recall below.
\begin{theorem}[Koppelman's Cauchy formula]
Let $\Omega$ be a bounded domain in $\mathbb C^n$ with $C^2$-smooth boundary. If $z\in\Omega$ and $f\in C^1(\overline{\Omega})$, then
\begin{equation}\label{20261114for1}
nc_nf(z)=\int_{\partial \Omega}f(\zeta)K_s(z,\zeta)-\int_{\Omega}(\overline{\partial}f)\wedge K_s(z,\zeta).
\end{equation}
Here, \(c_n\) is a dimensional constant depending only on \(n\), and \(K_s(z,\zeta)\) is the differential-form kernel appearing in Koppelman's formula.
\end{theorem}
The availability of such integral representations is one of the fundamental ingredients in the development of function theory on the finite-dimensional unit ball. This naturally raises the question of whether an analogous integral representation can be established in infinite dimensions.
However, a direct extension of this formula to the infinite-dimensional setting is highly nontrivial, primarily because the corresponding theories of the \(\overline{\partial}\)-equation, differential forms, Sobolev spaces, surface measures, Stokes-type theorems, integral operators, holomorphy, etc., have yet to be fully developed in infinite dimensions. Accordingly, the relevant developments can be summarized as follows:
\begin{itemize}
\item
For many years, research on the \(\overline{\partial}\)-equation in infinite dimensions focused primarily on strong solutions; see, for instance, \cite{Hen,Lem99}. Establishing weak solvability via H\"{o}rmander's \(L^2\) method \cite{Hor65} remained a longstanding open problem. In a series of works \cite{WYZ1,YZ}, the second author and X. Zhang introduced new notions of differential forms and weak partial derivatives in infinite dimensions. Within this framework, they established the weak solvability of the \(\overline{\partial}\)-equation using the \(L^2\) method. Addressing certain problems left unresolved in \cite{WYZ1,YZ} required the development of new tools from infinite-dimensional Sobolev theory. This led us to develop a new and fairly comprehensive theory of infinite-dimensional Sobolev spaces in \cite{WYZZZ}, revealing a variety of new phenomena specific to infinite dimensions.
\item
Section 1.4 of \cite[pp. 12--19]{Rud} is devoted to integration over spheres. More generally, the theory of integral operators often relies on the measures defined on both the interior and the boundary of a smooth domain. In infinite-dimensional settings, however, a fundamental difficulty is to construct an appropriate surface measure on the boundary that allows for a satisfactory development of the corresponding integral calculus. In a recent preprint \cite{WYZ}, Z. Wang, the second author and X. Zhang developed a convenient framework for constructing surface measures and establishing Gauss--Green and Stokes-type theorems in infinite dimensions. Based on this framework, we derive the key Forelli--Rudin-type estimates needed for the study of the integral operators considered in this paper.
\item
Furthermore, motivated by the following Lemma~\ref{20261113lem1}, we introduce an infinite-dimensional analogue of the classical Bergman space (see \cite[Section 3.1.1, p.~36]{Rud}) in \eqref{20261114def1}. An infinite-dimensional Hardy space was also introduced earlier in \cite[p.~545]{YZ}. These two spaces are closely related to the theory of monomial expansions, which dates back to the work of von Koch and Hilbert; see, for instance, \cite{DMP,Din99,Mujica}.

Historically, in the study of infinite-dimensional holomorphy, expansions in homogeneous polynomials were often preferred to monomial expansions. One of the reasons was the lack of suitable analytical tools for investigating holomorphy through the latter. However, Lempert's work \cite{Lem99} shed new light on the role of monomial expansions in infinite-dimensional holomorphy. Subsequently, in \cite{YZ19}, the second author and X. Zhang developed further tools for studying holomorphy in terms of monomial expansions. We expect these developments to open up new avenues for further research in infinite-dimensional holomorphy.
\end{itemize}

We hope that the present work, together with our previous contributions, will provide new insights into the appropriate formulation of integral formulas in infinite dimensions and, more broadly, contribute to the development of function theory on the infinite-dimensional unit ball.

The remainder of this paper is organized as follows. Subsection~\ref{20261105sbusection1} introduces the basic notation and presents the main results of the paper. Subsection~\ref{20261115subsec1} highlights several new phenomena revealed by our results and investigates the underlying mechanisms responsible for their occurrence in infinite dimensions. Section~\ref{20261105sbusectionnew2} is devoted to preliminary results and the key Forelli--Rudin-type estimates needed in the subsequent sections. Finally, Sections~\ref{20261105sbusection3} and \ref{20261105sbusection4} contain the proofs of the main results, which rely on Forelli--Rudin-type estimates, spectral estimates, and Fourier analysis.
\subsection{Basic notations and main results}\label{20261105sbusection1}
Suppose that $G$ is a
nonempty set and $(G,\cal J)$ is a topological space on $G$. The smallest $\sigma$-algebra
generated by all sets in $\cal J$ is called the
Borel $\sigma$-algebra of
$G$, denoted by $\mathscr{B}(G)$.
Denote by $\mathbb{R}$ and $\mathbb{C}$ respectively the sets of all real and complex numbers, and by $\mathbb{N}$ the set of all positive integers. Write $\mathbb{N}_0\triangleq \mathbb{N}\cup \{0\}$ and $\mathbb{C}^{\infty}\triangleq\prod\limits_{i=1}^{\infty}\mathbb{C}$, which is the countable Cartesian product of $\mathbb{C}$, i.e., the space of sequences of complex numbers (Nevertheless, in the sequel we shall write an element $\textbf{z}\in \mathbb{C}^{\infty }$ explicitly as $\textbf{z}=( z_i)_{i\in \mathbb N }$ rather than $\textbf{z}=\{ z_{i}\}^{\infty }_{i=1}$, where $z_i\in\mathbb{C}$ for each $i\in \mathbb{N}$). Clearly, $\mathbb{C}^{\infty}$ can be identified with $(\mathbb{R}^2)^{\infty}=\prod\limits_{i=1}^{\infty}\mathbb{R}^2$. Let
\[
 \HH\triangleq\left\{\mathbf z=(z_j)_{j\in\mathbb N}\in\C^{\infty}:
                 \sum_{j=1}^{\infty}|z_j|^2<\infty\right\}.
\]
Then, $\ell^2$ is a separable Hilbert space with the canonical inner product $ \langle \cdot,\cdot\rangle$ (and hence  the canonical norm $\left| \left|\;\cdot\;\right|  \right|$). For any $i,j\in\mathbb{N}$, set $\delta_{ij} \triangleq 1$, if $i=j$, and $\delta_{ij} \triangleq 0$, if $i\neq j.$ For each $k\in\mathbb{N}$, define $\textbf{e}_k\triangleq(\delta_{k,i})_{i\in\mathbb{N}}$. Denote by $B\left( \textbf{x} ,r \right)$ the open ball in $\ell^{2}$ centered at $\textbf{x}(\in \ell^{2})$ with radius $r>0$, and we simply write $B\left( \textbf{0} ,r \right)$ as $B_r$. Fix a sequence of positive numbers \(\{a_j\}_{j=1}^{\infty}\) satisfying \(\sum_{j=1}^{\infty}a_j^2<\infty\). Set
$$
 \vartheta_j\triangleq 2a_j^2,\qquad j\in\mathbb N,\qquad \text{and}\qquad  \vartheta \triangleq (\vartheta_j)_{j\in\mathbb N}.
$$
For $a\in(0,\infty)$, let
\[
 \dd\mathcal N_{a } ( z)
\triangleq\frac{1}{2\pi a^2}
   \exp\left(-\frac{\| z\|^2}{2a^2}\right)\dd  z,
\]
where \(\dd  z\) is the Lebesgue measure on $\mathbb C$. Let $\mathcal N=\prod_{i=1}^{\infty}\mathcal N_{a_i }$. By the same argument as in Section 2.2 in \cite{YZ} and Section 2.3 in \cite{WYZ1}, we obtain a probability measure $P$ on $(\ell^2, \mathscr{B} (\ell^2))$  by
$$
P(E)\triangleq \mathcal N(E),\qquad E\in\mathscr B(\ell^2).
$$
Set
\[
 B\triangleq\{\mathbf z\in \ell^2:\norm {\mathbf z}<1\},\qquad
 \delta(\mathbf z)\triangleq 1-\norm {\mathbf z}^2,\quad \mathbf z\in \ell^2,\qquad
 \mu\triangleq P|_B,\qquad m\triangleq \mu(B).
\]
For \(\gamma\in\R\), let
\begin{equation}\label{eq:weightedmeasure}
 \dd\nu_\gamma(\mathbf z)\triangleq \delta(\mathbf z)^\gamma\dd\mu(\mathbf z).
\end{equation}
We prove below that \(\nu_\gamma(B)<\infty\) precisely when
\(\gamma>-1\). For all real \(\gamma\) the measure is sigma-finite
and has the same null sets as \(\mu\).
Write
\[
 \Lt p\gamma\triangleq L^p(B,\nu_\gamma),\qquad
 \norm f_{p,\gamma}
 \triangleq\left(\int_B|f(\mathbf z)|^p\delta(\mathbf z)^\gamma\dd\mu(\mathbf z)\right)^{1/p},
 \quad 0<p<\infty,\qquad \mathcal H_\gamma \triangleq L^2(B,\nu_\gamma).
\]
For \(p=\infty\) the norm is the essential supremum; it is independent
of \(\gamma\).

The map
\begin{equation}\label{20261031for1}
 U_\gamma:\mathcal H_\gamma\longrightarrow L^2(B,\mu),
 \qquad U_\gamma f\triangleq\delta^{\gamma/2}f,
\end{equation}
is unitary, including when \(\nu_\gamma(B)=\infty\).

Similar to \eqref{eq:standardTn} and \eqref{eq:standardSn}, we have the following definition for $a,b,c\in\mathbb R$:
\begin{align}
 \TT f(\mathbf z)
 &\triangleq\delta(\mathbf z)^a\int_B
       \frac{\delta(\mathbf y)^b f(\mathbf y)}{(1-\langle\mathbf z, \mathbf y\rangle)^c}\dd\mu(\mathbf y),\label{eq:standardT}\\
 \SSS f(\mathbf z)
 &\triangleq\delta(\mathbf z)^a\int_B
       \frac{\delta(\mathbf y)^b f(\mathbf y)}{|1-\langle\mathbf z, \mathbf y\rangle|^c}\dd\mu(\mathbf y),\label{eq:standardS}
\end{align}
where $\mathbf z\in B$ and
for nonintegral \(c\), we use the principal logarithm:
\[
 (1-\langle\mathbf z, \mathbf y\rangle)^{-c}
   =\exp\bigl(-c\Log(1-\langle\mathbf z, \mathbf y\rangle)\bigr).
\]
Since
\[
 \Ree(1-\langle\mathbf z, \mathbf y\rangle)\ge1-\norm{ \mathbf z}\cdot\norm{\mathbf y}>0,\qquad \mathbf z,\mathbf y\in B,
\]
this definition is single-valued and
\(\abs{(1-\langle\mathbf z, \mathbf y\rangle)^{-c}}=|1-\langle\mathbf z, \mathbf y\rangle|^{-c}\).

 The index \(r\) below is a Schatten exponent, not an exponent of a
function space. For Hilbert spaces \( H_1, H_2\), write
\(\mathcal S_r( H_1, H_2)\), \(0<r<\infty\), for
the compact operators whose singular values $\{s_j(A):j\in J\subset \mathbb N\}$ satisfy
\[
 \|A\|_{\mathcal S_r}^r\triangleq\sum_{j\in J}s_j(A)^r<\infty.
\]
This definition applies to operators between different Hilbert spaces;
an ordinary operator trace, in contrast, requires an endomorphism.
We use the standard ideal, approximation-number, and tensor-product
properties of Schatten classes; see Simon \cite{Simon}.




\begin{theorem}\label{thm:standard}
Let \(a,b,c,\alpha,\beta\in\R\) and \(1\le p,q\le\infty\).
The following are equivalent:
\begin{enumerate}
 \item \(\TT:\Lt p\alpha\to\Lt q\beta\) is bounded.
 \item \(\SSS:\Lt p\alpha\to\Lt q\beta\) is bounded.
 \item The parameters satisfy the conditions in the corresponding row
       of Table~\ref{tab:standard}.
\end{enumerate}
\end{theorem}

\begin{table}[htbp]
\centering
\renewcommand{\arraystretch}{1.6}
\begin{tabular}{@{}p{32mm}p{111mm}@{}}
\toprule
Exponent & Necessary and sufficient conditions\\
\midrule
\(p=1\), \(1\le q<\infty\)
 & \(b\ge\alpha,\quad aq+\beta>-1\)\\
\(1<p<\infty\), \(1\le q<\infty\)
 & \(p(b+1)>\alpha+1,\quad aq+\beta>-1\)\\
\(p=\infty\), \(1\le q<\infty\)
 & \(b>-1,\quad aq+\beta>-1\)\\
\(p=1\), \(q=\infty\)
 & \(a\ge0,\quad b\ge\alpha,\quad c\le a+b-\alpha\)\\
\(1<p<\infty\), \(q=\infty\)
 & \(a\ge0,\quad p(b+1)>\alpha+1\)\\
\(p=\infty\), \(q=\infty\)
 & \(a\ge0,\quad b>-1\)\\
\bottomrule
\end{tabular}
\caption{The common boundedness region for \(T_{a,b,c}\) and \(S_{a,b,c}\).
}
\label{tab:standard}
\end{table}
\newpage
\begin{theorem} \label{thm:schatten2}
Let all five operator and weight parameters be real. Put
\begin{equation}\label{eq:ball-effective}
 u\triangleq a+\frac\beta2,\quad v\triangleq b-\frac\alpha2,\quad
 t\triangleq u+v=a+b+\frac{\beta-\alpha}{2}.
\end{equation}
Either of
\(T_{a,b,c},S_{a,b,c}:\mathcal H_\alpha\to\mathcal H_\beta\)
is bounded if and only if
\begin{equation}\label{eq:hilbertadmissible}
 u>-\tfrac12,\qquad v>-\tfrac12,
 \quad\text{equivalently}\quad
 2a+\beta>-1,\quad 2b-\alpha>-1.
\end{equation}
 Under \eqref{eq:hilbertadmissible}, both operators are compact,
and their complete classification is as follows.
\begin{enumerate}
\item For every \(1<r<\infty\), both operators belong to
\(\mathcal S_r\), for every real \(c\).
\item Both operators belong to \(\mathcal S_1\) if and only if
\begin{equation}\label{eq:tracecriterion}
 c<t+1.
\end{equation}
In this case,
\begin{equation}\label{20261105for1}
 \operatorname{Tr}(U_\beta T_{a,b,c}U_\alpha^{-1})=\operatorname{Tr}(U_\beta S_{a,b,c}U_\alpha^{-1})
 =\int_B\delta^{a+b+\frac{\beta-\alpha}{2}-c}\dd\mu.
\end{equation}

\item For \(0<r<1\), the exact criteria are
\begin{align}
 T_{a,b,c}\in\mathcal S_r
 &\quad\Longleftrightarrow\quad
 \left\{\begin{array}{l}
 c=0,\quad\text{or}\\
 c\in\{-1,-2,\ldots\}\ \text{and}\
       \displaystyle\sum_j\vartheta_j^r<\infty;
 \end{array}\right.\label{eq:subunitT}\\
 S_{a,b,c}\in\mathcal S_r
 &\quad\Longleftrightarrow\quad
 \left\{\begin{array}{l}
 c=0,\quad\text{or}\\
 c\in\{-2,-4,\ldots\}\ \text{and}\
       \displaystyle\sum_j\vartheta_j^r<\infty.
 \end{array}\right.\label{eq:subunitS}
\end{align}
\end{enumerate}
For \(c=0\), both operators have rank one, with the single nonzero
singular value \(\sqrt{I_{2u}I_{2v}}\). Every nonconstant polynomial
kernel appearing in \eqref{eq:subunitT} or \eqref{eq:subunitS} has
infinite rank.
\end{theorem}


\subsection{New phenomena unique to infinite dimensions}\label{20261115subsec1}
In this subsection, we compare our results with their finite-dimensional counterparts and investigate the underlying mechanisms responsible for the emergence of these phenomena in infinite dimensions.


Zhao and Zhou \cite{ZZ} characterized the boundedness of the corresponding operators on the unit ball of \(\mathbb C^n\) with weighted volume measures; see Table~\ref{tab:standard2}. Comparing Tables~\ref{tab:standard2} and \ref{tab:standard}, we observe two striking differences:
\begin{itemize}
\item The endpoint condition remains unchanged;
\item All other restrictions on the parameter \(c\) disappear in infinite dimensions.
\end{itemize}

This phenomenon is essentially a consequence of the differences in measure and geometry between finite- and infinite-dimensional spaces, rather than a formal substitution of \(n=\infty\) into finite-dimensional results. To understand why the admissible range of \(c\) differs between the two settings, we examine the key Forelli--Rudin-type estimates underlying the boundedness results.

In finite dimensions, for each fixed \(n\in\mathbb N\), \cite[Corollary 5.24, p.~183]{Zhu05} yields the following estimate:
\[
\sup_{\mathbf z\in\mathbb B_n}
v_n\bigl(\{\mathbf y\in\mathbb B_n:
|1-\langle\mathbf z,\mathbf y\rangle|<\varepsilon\}\bigr)
\leq C\varepsilon^{\kappa_n}.
\]
Here, the exponent \(\kappa_n\) depends on the dimension \(n\). Consequently, the restrictions on \(c\) in the \(L^p\)-\(L^q\) boundedness results for \((p,q)\neq(1,\infty)\) are dimension-dependent.

In contrast, as established in Lemma~\ref{lem:caps}, the corresponding estimate in infinite dimensions takes the following form:
\begin{equation}\label{20261109for1}
\sup_{\mathbf z\in B}
\mu\bigl(\{\mathbf y\in B:
|1-\langle\mathbf z,\mathbf y\rangle|<\varepsilon\}\bigr)
\leq C_N\varepsilon^N,
\qquad \text{for every }N\in\mathbb N.
\end{equation}
Thus, the measure decays faster than every prescribed polynomial power of \(\varepsilon\). This super-polynomial decay is the key reason why no additional dimension-dependent restriction on \(c\) is required for \((p,q)\neq(1,\infty)\).

We next explain why the endpoint \(L^1\to L^\infty\) behaves differently. At this endpoint, the operator norm is governed by the pointwise supremum of the integral kernel, rather than by averaged estimates of the kernel over small sets. Indeed, Theorem~\ref{thm:endpoint} establishes that
\[
\begin{aligned}
\|T_{a,b,c}\|_{L^1_\alpha\to L^\infty_\beta}
&=\|S_{a,b,c}\|_{L^1_\alpha\to L^\infty_\beta}\\
&=\sup_{\mathbf z,\mathbf y\in B}
\frac{\delta(\mathbf z)^a\delta(\mathbf y)^{b-\alpha}}
{|1-\langle\mathbf z,\mathbf y\rangle|^c}.
\end{aligned}
\]
The expression on the right is independent of the dimension. Consequently, the endpoint condition agrees with its finite-dimensional counterpart. In other words, unlike the non-endpoint cases, the endpoint is determined entirely by a pointwise kernel-size condition and is therefore unaffected by the dimension-dependent Forelli--Rudin-type estimates.

In a nutshell, the dimension-dependent restrictions in finite dimensions are replaced by super-polynomial decay in infinite dimensions, whereas the endpoint remains governed by the pointwise size of the kernel.

Setting \(a=b=\alpha=\beta=0\) gives
\[
T_{0,0,c}f(\mathbf z)=\int_B(1-\ip {\mathbf z}{\mathbf y})^{-c}f(\mathbf y)\dd\mu(\mathbf y),\qquad S_{0,0,c}f(\mathbf z)=\int_B|1-\ip {\mathbf z}{\mathbf y}|^{-c}f(\mathbf y)\dd\mu(\mathbf y).
\]
Theorem~\ref{thm:standard} shows that every pair \((p,q)\in[1,\infty]^2\) is admissible when \(c\le0\), whereas every such pair except \((1,\infty)\) is admissible when \(c>0\). The corresponding results for the unit disk, obtained in \cite{CFWY}, are summarized in Table~\ref{tab:standard1}. Remarkably, the boundedness criteria for this unweighted problem are considerably more intricate in one complex dimension than in infinite dimensions. Although we do not yet have a complete explanation for this contrast, the Forelli--Rudin-type asymptotic estimate in \eqref{20261109for1} accounts for at least part of the difference.

In both finite and infinite dimensions, the two Forelli--Rudin-type
operators satisfy the same \(L^p\)-\(L^q\) boundedness criteria.
Their Schatten class behavior, however, exhibits a new distinction
in infinite dimensions when \(c\) is a negative odd integer and
\(0<r<1\), as described in Conclusion~3 of
Theorem~\ref{thm:schatten2}. In this case, the complex-kernel
operator \(T_{a,b,c}\) has a polynomial kernel and is governed by
the covariance summability condition in \eqref{eq:subunitT},
whereas the modulus-kernel operator \(S_{a,b,c}\) does not have a
polynomial kernel and fails to belong to \(\mathcal S_r\).

In finite
dimensions, a polynomial kernel defines a finite-rank integral
operator. Thus, for negative integer \(c\), the kernel of
\(T^n_{a,b,c}\) is polynomial; for negative even integer \(c\), the
kernel of \(S^n_{a,b,c}\) is also polynomial. These cases therefore
give finite-rank operators in finite dimensions. When \(c\) is a
negative odd integer, however, the kernel
\(|1-\langle\mathbf z,\mathbf y\rangle|^{-c}\) is not a polynomial,
so no finite-rank conclusion for \(S^n_{a,b,c}\) follows from this
argument. In the infinite-dimensional setting, every nonconstant
polynomial kernel appearing in \eqref{eq:subunitT} or
\eqref{eq:subunitS} defines an infinite-rank operator. This
difference in rank structure explains the covariance-summability
conditions in the polynomial cases of Theorem~\ref{thm:schatten2}.

Ding \cite[Theorems 3.1, 3.2, 5.2, and 5.3, pp.~8 and 26]{Ding} characterized the Schatten-class membership of Bergman-type and Szeg\"{o}-type operators on bounded symmetric domains for certain special parameter values. In comparison with Ding's finite-dimensional results, Conclusion~3 of Theorem~\ref{thm:schatten2} reveals another phenomenon that is unique to infinite dimensions.



For each \(n\in\mathbb N\), using polar coordinates, we obtain
\begin{equation*}\label{20261108for1}
\int_{\mathbb B_n}(1-\|\mathbf z\|^2)^{t-c}\,\mathrm d v_n(\mathbf z)
=n\int_0^1(1-s)^{t-c}s^{n-1}\,\mathrm d s.
\end{equation*}
In contrast, the corresponding trace formula in infinite dimensions, given by \eqref{20261105for1}, takes the form
\begin{equation*}\label{eq:gaussian-trace}
\int_B(1-\|\mathbf z\|^2)^{t-c}\,\mathrm d\mu(\mathbf z)
=\int_0^1(1-s)^{t-c}\widetilde h(s)\,\mathrm d s,
\end{equation*}
where \(\widetilde h\) is a continuous function on \(\mathbb R\).

Thus, the infinite-dimensional trace formula is not obtained as a formal limit of the finite-dimensional formula as \(n\to\infty\). Indeed, the density \(\widetilde h\) arising from the Gaussian measure has a fundamentally different structure from the finite-dimensional density \(ns^{n-1}\).








\section{Preliminaries and Forelli--Rudin-type estimates}\label{20261105sbusectionnew2}
Let
\(\mathcal D\) be the bounded Borel measurable functions that vanish outside
some ball \(B_R\), where \(R<1\).
Every displayed integral in \eqref{eq:standardT} and \eqref{eq:standardS} converges absolutely on \(\mathcal D\)
at each \(\mathbf z\in B\). For \(1\leq p<\infty\), by \cite[Theorem 2.1, p. 8]{WYZ1}, this class is dense in
\(\Lt p\gamma\).

\begin{proposition}\label{prop:support2}
It holds that \(m \in(0,1)\), and $P(O)>0$ for any non-empty open subset $O$ of $\ell^2$.
For every \(\gamma\in\mathbb R\), \(\nu_\gamma\) is sigma-finite measure, is equivalent
to \(\mu\).
\end{proposition}

\begin{proof}
The first conclusion follows from \cite[Lemma 2.1, p. 9]{WYZ1}. For every \(R<1\), on \(B_R\), the function \(\delta^\gamma\) is bounded
above and below by positive constants.
Consequently \(\nu_\gamma(B_R)<\infty\); exhaustion by such balls
proves sigma-finiteness. Since \(0<\delta^\gamma<\infty\) throughout
\(B\), equivalence follows.
\end{proof}

Set
$$
X(\mathbf z)\triangleq  \sum_{j=1}^{\infty}|z_j|^2,\qquad X_j(\mathbf z)\triangleq  |z_j|^2,\qquad \mathbf z=(z_j)_{j\in\mathbb N}\in\ell^2,\qquad\lambda_j\triangleq\frac{1}{2a_j^2},\quad j\in\mathbb N.
$$
\begin{lemma} \label{lem:radialdensity}
The distribution of $X$ under $P$ is absolutely continuous with respect
to Lebesgue measure on $(0,\infty)$, with a nonnegative density
$h$. Thus, for every Borel set $A\subset(0,\infty)$,
\[
P(\{\mathbf z\in\ell^2:X(\mathbf z)\in A\})
=\int_A h(s)\,\mathrm ds.
\]
Moreover, there exist constants $0<\ell\le L<\infty$ such that
\begin{equation}\label{eq:radialdensity}
 h(s)\le L\quad(s>0),\qquad
 h(s)\ge\ell\quad(1/2\le s\le1).
\end{equation}
\end{lemma}

\begin{proof}
Set
\[
 R(\mathbf z)\triangleq  \sum_{j=2}^{\infty}|z_j|^2,\qquad \mathbf z=(z_j)_{j\in\mathbb N}\in\ell^2.
\]

For every $j\in\mathbb N$, nonnegative Borel function $\varphi$, polar coordinates and the substitution $u=r^2$ give
\begin{align*}
 \int_{\ell^2}\varphi(X_j)\,\mathrm{d}P
 &=\frac{1}{2\pi a_j^2}\int_{\mathbb C}
   \varphi(| w|^2)e^{-| w|^2/(2a_j^2)}\,\mathrm d  w
 =\int_0^\infty\varphi(u)\lambda_j e^{-\lambda_j u}\,\mathrm du.
\end{align*}
Thus the distribution of $X_j$ under $P$ has density $\lambda_j e^{-\lambda_j u}\ind_{\{u>0\}}$ with respect to Lebesgue measure on $(0,\infty)$.

Let $R_3\triangleq\sum_{j=3}^{\infty}X_j$. Since
$\int_{\ell^2} R_3\,\mathrm d P=\sum_{j=3}^{\infty}\vartheta_j<\infty$,
the random variable $R_3$ is finite almost surely. It is independent
of $X_2$, whose distribution has an exponential density. Hence
$R=X_2+R_3$ has a density with respect to Lebesgue measure, namely
\[
\rho(r)\triangleq
\int_{\ell^2} \left[
\lambda_2e^{-\lambda_2(r-R_3)}
\ind_{\{R_3<r\}}
\right]\,\mathrm d P,\qquad r>0.
\]
We may therefore define
\[
 h(s)\triangleq\int_{[0,\infty)}
 \lambda_1 e^{-\lambda_1(s-r)}\ind_{\{r<s\}}\,\rho( r)\,\mathrm d r,
 \qquad s>0.
\]
The integrand is a nonnegative Borel function of $(s,r)$, so $h$ is Borel measurable.

For a Borel set $A\subset(0,\infty)$, Tonelli's theorem yields
\begin{align*}
P(\{\mathbf z\in\ell^2:X(\mathbf z)\in A\})
 &=\int_{[0,\infty)}\int_{[0,\infty)}
   \ind_A(u+r)\lambda_1 e^{-\lambda_1 u}\,\mathrm d u\,\rho( r)\,\mathrm d r\\
 &=\int_{[0,\infty)}\int_{[r,\infty)}
   \ind_A(s)\lambda_1 e^{-\lambda_1(s-r)}\,\mathrm d s\,\rho( r)\,\mathrm d r\\
 &=\int_A h(s)\,\mathrm d s.
\end{align*}
Since $X>0$ almost surely, the preceding identity shows that $h$ is a probability density of the distribution of $X$ with respect to Lebesgue measure.
For every $s>0$ and $r<s$ we have $e^{-\lambda_1(s-r)}\le1$. Hence
$
 0\le h(s)\le\lambda_1,
$
so we may take $L\triangleq\lambda_1$.

By Proposition \ref{prop:support2},
\[
 p_0\triangleq P(\{\mathbf z:R(\mathbf z)<1/4\})
 \ge P(\{\mathbf z:\|\mathbf z\|<1/2\})>0.
\]
If $\tfrac12\le s\le1$ and $0\le r<\tfrac14$, then $r<s$ and $0<s-r\le1$. Consequently,
\[
 h(s)\ge\int_{[0,1/4)}\lambda_1 e^{-\lambda_1(s-r)}\,\rho(r)\mathrm dr
 \ge\lambda_1 e^{-\lambda_1}p_0.
\]
The required bounds therefore hold pointwise for
$\ell\triangleq\lambda_1 e^{-\lambda_1}p_0>0.
$
\end{proof}

\begin{corollary} \label{cor:radial}
For \(t\in\R\), define
$
 I_t\triangleq\int_B\delta(\mathbf z)^t\dd\mu(\mathbf z).
$
Then
\begin{equation}\label{eq:radialthreshold}
 I_t<\infty\quad\Longleftrightarrow\quad t>-1.
\end{equation}
For \(t>-1\), \(I_t\le L/(t+1)\). Also,
\[
 \delta^u\in L^\infty(B,\mu)
 \quad\Longleftrightarrow\quad u\ge0.
\]
The \(L^\infty\)-norm is $1$ when \(u\ge0\).
\end{corollary}

\begin{proof}
By Lemma~\ref{lem:radialdensity},
\[
 I_t=\int_0^1(1-s)^t h(s)\dd s.
\]
The upper bound and the divergence for \(t\le-1\) follow from
\eqref{eq:radialdensity}.
For the essential-supremum assertion, \(0<\delta\le1\), and neighborhoods
of \(0\) have positive measure, giving norm one for \(u\ge0\).
For \(u<0\), \eqref{eq:radialdensity} makes
\(\delta^u\) essentially unbounded.
\end{proof}

Write
$d(\mathbf z,\mathbf y)\triangleq|1-\langle\mathbf z,\mathbf y\rangle|,
 \text{ for every } \mathbf z\in\overline B,\,\mathbf y\in B.
$

\begin{lemma} \label{lem:caps}
For every \(N\in\mathbb N\),
\begin{equation}\label{eq:caps2}
 \sup_{\mathbf z\in\overline B}
 \mu(\{\mathbf y\in B:d(\mathbf z,\mathbf y)<u\})\le C_Nu^N,\quad u>0 ,
 \qquad
 C_N\triangleq\frac{1}{N!\prod_{j=1}^Na_j^2}.
\end{equation}
\end{lemma}

\begin{proof}
For $N\in\mathbb N$, set
\[
\pi_N:\ell^2\longrightarrow\C^N,\qquad
\pi_N(\mathbf z)\triangleq (z_1,\ldots,z_N),\qquad \mathbf z=(z_i)_{i\in\mathbb N}\in\ell^2.
\]
Let
\[
f_N(\mathbf w)
\triangleq\prod_{j=1}^N\frac{1}{2\pi a_j^2}
\exp\!\left(-\frac{|w_j|^2}{2a_j^2}\right)
\le\frac{1}{(2\pi)^N\prod_{j=1}^Na_j^2},\qquad \mathbf w=(w_1,\ldots,w_N)\in \mathbb C^N.
\]
For $\mathbf x=(x_i)_{i\in\mathbb N}\in \ell^2$ and $r>0$,
\[
\sum_{j=1}^N|y_j-x_j|^2
\le\sum_{j=1}^{\infty}|y_j-x_j|^2
=\norm{\mathbf y-\mathbf x}^2,\qquad \mathbf y=(y_i)_{i\in\mathbb N}\in\ell^2,
\]
and hence
\[
\{\mathbf y \in\ell^2:\norm{\mathbf y-\mathbf x}<r\}
\subseteq
\{\mathbf y \in\ell^2:\norm{\pi_N(\mathbf y)-\pi_N(\mathbf x)} <r\}.
\]
Consequently,
\begin{align}
P(B (\mathbf x,r))
&\le\int_{\{\mathbf w\in\mathbb C^N:\norm{\mathbf w-\pi_N(\mathbf x)}<r\}}f_N(\mathbf w)\,\mathrm d\mathbf w \notag
\le\frac{\int_{\{\mathbf w\in\mathbb C^N:\norm{\mathbf w-\pi_N(\mathbf x)}<r\}} \,\mathrm d\mathbf w}
{(2\pi)^N\prod_{j=1}^Na_j^2}\notag\\
&=\frac{\pi^Nr^{2N}}
{\Gamma(N+1)(2\pi)^N\prod_{j=1}^Na_j^2}
=\frac{r^{2N}}{2^N(N!)\prod_{j=1}^Na_j^2}.
\label{eq:ball}
\end{align}
The same estimate holds for closed balls.

For $\mathbf z\in\overline B$ and $u>0$, put
$
C(\mathbf z,u)\triangleq\{\mathbf y\in B:|1-\ip{\mathbf z}{\mathbf y}|<u\}.
$
For every $\mathbf y\in C(\mathbf z,u)$,
\begin{align*}
\norm{\mathbf y-\mathbf z}^2
&=\norm{\mathbf y}^2+\norm{\mathbf z}^2-2\operatorname{Re}\ip{\mathbf z}{\mathbf y}\
\le 2-2\operatorname{Re}\ip{\mathbf z}{\mathbf y}
=2\operatorname{Re}(1-\ip{\mathbf z}{\mathbf y})
\le 2|1-\ip{\mathbf z}{\mathbf y}|
<2u.
\end{align*}
Thus
$
C(\mathbf z,u)\subseteq B (\mathbf z,\sqrt{2u}).
$
Using $\mu=P|_B$ and \eqref{eq:ball},
\begin{align*}
\mu(C(\mathbf z,u))
=P(C(\mathbf z,u))
\le P(B (\mathbf z,\sqrt{2u}))
\le\frac{(\sqrt{2u})^{2N}}
{2^N(N!)\prod_{j=1}^Na_j^2}
=\frac{u^N}{N!\prod_{j=1}^Na_j^2}.
\end{align*}
Therefore, \eqref{eq:caps2} holds.
\end{proof}

\begin{lemma} \label{lem:unweighted}
For every \(s\in\mathbb R\),
\[
 M_s\triangleq \sup_{\mathbf z\in\overline B}\int_Bd(\mathbf z,\mathbf y)^{-s}\dd\mu(\mathbf y)<\infty.
\]
For \(s>0\) and any integer \(N>s\),
$
 M_s\le m+\frac{sC_N}{N-s}.
$
For \(s\le0\), \(M_s\le2^{-s}m\).
\end{lemma}

\begin{proof}
The case \(s\le0\) follows from \(d(\mathbf z,\mathbf y)\le2\).
For \(s>0\), use
\[
 v^{-s}\le1+s\int_0^1u^{-s-1}\ind_{\{v<u\}}\dd u,
 \qquad v>0.
\]
After integration against \(\mu\), Tonelli's theorem and
\eqref{eq:caps2} give
\[
 \int_Bd(\mathbf z,\mathbf y)^{-s}\dd\mu(\mathbf y)
 \le m+sC_N\int_0^1u^{N-s-1}\dd u
 =m+\frac{sC_N}{N-s}.
\]
\end{proof}

\begin{theorem} \label{thm:FR}
For every \(t,s\in\R\), put
$
 W_{t,s}\triangleq\sup\limits_{\mathbf z\in\overline B}
          \int_B\frac{\delta(\mathbf y)^t}{d(\mathbf z,\mathbf y)^s}\dd\mu(\mathbf y).
$
Then
\begin{equation}\label{eq:FR}
 W_{t,s}<\infty\quad\Longleftrightarrow\quad t>-1.
\end{equation}
For \(t>-1\), choose any \(r>1\) with \(tr>-1\), and put
\(r'\triangleq r/(r-1)\). One explicit bound is
\begin{equation}\label{eq:weightedmomentbound}
 W_{t,s}\le I_{tr}^{1/r}M_{sr'}^{1/r'}.
\end{equation}
All these estimates hold with the integration and fixed variables interchanged.
\end{theorem}

\begin{proof}
If \(t>-1\), an \(r>1\) satisfying \(tr>-1\) exists.
H\"older's inequality gives, uniformly in \(\mathbf z\in\overline B\),
\[
 \int_B\delta(\mathbf y)^t d(\mathbf z,\mathbf y)^{-s}\dd\mu(\mathbf y)
 \le
 \left(\int_B\delta(\mathbf y)^{tr}\dd\mu(\mathbf y)\right)^{1/r}
 \left(\int_Bd(\mathbf z,\mathbf y)^{-sr'}\dd\mu(\mathbf y)\right)^{1/r'}.
\]
The two factors are finite by Corollary~\ref{cor:radial} and
Lemma~\ref{lem:unweighted}.
Conversely, setting \(\mathbf z=\mathbf 0\) gives \(d(\mathbf 0,\mathbf y)=1\), so
\(W_{t,s}\ge I_t=\infty\) if \(t\le-1\).
The last assertion follows from \(d(\mathbf z,\mathbf y)=d(\mathbf y,\mathbf z)\).
\end{proof}


Let \(\mathbb{N}_0^{(\mathbb{N})}\) be the countable set of finitely supported
multi-indices in \(\mathbb N_0^{\mathbb N}\). For
\(\eta=(\eta_j)_{j\in\mathbb N}\in \mathbb{N}_0^{(\mathbb{N})},\) $s\in\mathbb R$, and $\mathbf z=(z_j)_{j\in\mathbb N}\in \ell^2$, put
\begin{equation}\label{2026111def1}
 |\eta|\triangleq \sum_j\eta_j,\quad \eta!\triangleq\prod_j\eta_j!,\quad
 \mathbf z^\eta\triangleq \prod_j z_j^{\eta_j},\quad
 \vartheta^\eta\triangleq\prod_j\vartheta_j^{\eta_j},\quad
 M_\eta(s)\triangleq\int_B\delta(\mathbf z)^s|\mathbf z^\eta|^2\dd\mu(\mathbf z).
\end{equation}
Also define, for \(s>-1\),
\begin{equation}\label{eq:Jdefinition}
 J_n(s)\triangleq\int_B\delta(\mathbf z)^s\|\mathbf z\|^{2n}\dd\mu(\mathbf z)
       =\int_0^1(1-x)^s x^n h(x)\dd x,\qquad n\in\mathbb N_0,
\end{equation}
where $h$ is defined in Lemma \ref{lem:radialdensity}.

\begin{theorem}\label{lem:moments}
The density \(h\) given in Lemma \ref{lem:radialdensity} has a continuous version \(\widetilde{h}\) on \(\mathbb R\), extended
by zero on \(( -\infty,0]\), with \(\widetilde{h}(1)>0\). For each \(s>-1\),
\begin{align}
 J_n(s)&\sim \widetilde{h}(1)\Gamma(s+1)n^{-s-1},\quad n\to\infty ,
       \label{eq:Jasymptotic}\\
 M_\eta(s)&\le \frac{L}{s+1}\,\eta!\vartheta^\eta,
       \label{eq:momentupper}\\
 \frac{M_{\mathbf e_j}(s)}{\vartheta_j}&\longrightarrow I_s,
       \quad j\to\infty .
       \label{eq:firstmoment}
\end{align}
Here \(L\) is the bound in Lemma~\ref{lem:radialdensity}, and
\(\mathbf e_j\) is the \(j\)-th unit multi-index.
\end{theorem}

\begin{proof}
\renewcommand{\qedsymbol}{$\square$}

\medskip
\noindent\textbf{Step 1. A bounded continuous version of $h$.}

By the proof of Lemma \ref{lem:radialdensity}, the density of $X_j$ is
\[
 f_j(t)=\frac1{\vartheta_j}e^{-t/\vartheta_j}
 \mathbf 1_{(0,\infty)}(t),\qquad j\in\mathbb N,\quad t\in\mathbb R.
\]
Tonelli's theorem gives
\[
 \int X\,\mathrm d P=\sum_{j=1}^{\infty}\int X_j\,\mathrm d P
       =\sum_{j=1}^{\infty}\vartheta_j<\infty.
\]
Thus $X<\infty$ almost everywhere with respect to $P$.

Let $g\triangleq f_1*f_2$, which is the density of $X_1+X_2$. Explicitly,
\[
 g(t)=
 \begin{cases}
 \displaystyle
 \frac1{\vartheta_1\vartheta_2}
 \int_0^t e^{-u/\vartheta_1}e^{-(t-u)/\vartheta_2}\dd u,
 &t>0,\\[7pt]
 0,&t\le0.
 \end{cases}
\]
Since $f_1$ and $f_2$ are probability densities,
$
 0\le g(t)\le
 \min\left\{\frac1{\vartheta_1},\frac1{\vartheta_2}\right\},\text{ for all } t\in\mathbb R.
$

Also, for $t>0$,
$
 0\le g(t)\le\frac{t}{\vartheta_1\vartheta_2}.
$
Consequently, $g$ is continuous at $0$.

We will use the stronger fact that $g$ is globally Lipschitz. For $t>0$,
differentiation of the integral formula yields
\[
 g'(t)=\frac1{\vartheta_2}\bigl(f_1(t)-g(t)\bigr).
\]
Both $f_1(t)$ and $g(t)$ belong to $[0,\vartheta_1^{-1}]$, so
$
 |g'(t)|\le\frac1{\vartheta_1\vartheta_2}.
$
The mean value theorem gives the same Lipschitz bound between two positive
arguments. If $u\le0<t$, then
\[
 |g(t)-g(u)|=g(t)
 \le\frac{t}{\vartheta_1\vartheta_2}
 \le\frac{t-u}{\vartheta_1\vartheta_2}.
\]
For two nonpositive arguments the difference is zero. Hence
\begin{equation}
 |g(t)-g(u)|\le C_0|t-u|,\quad t,u\in\R,
 \qquad C_0\triangleq\frac1{\vartheta_1\vartheta_2}.
 \label{eq:g-lipschitz}
\end{equation}

Now set
\[
 R_2\triangleq\sum_{j=3}^{\infty}X_j.
\]
This is finite almost everywhere with respect to $P$ and independent of $X_1+X_2$. Define
\begin{equation}
\widetilde{ h}(t)\triangleq\int g(t-R_2)\,\mathrm d P ,\qquad t\in\R.
 \label{eq:h-representative}
\end{equation}
Independence and Tonelli's theorem show that this is also a density function of $X$.
Indeed, for every nonnegative Borel function $\varphi$,
\begin{align*}
 \int\varphi(X)\,\mathrm d P
 =\int\left(\int_{\R}\varphi(u+R_2)g(u)\dd u \right)\,\mathrm d P =\int_{\R}\varphi(t)\left(\int g(t-R_2) \,\mathrm d P\right) \dd t
  =\int_{\R}\varphi(t)\widetilde{h}(t)\dd t.
\end{align*}
The bounds for $g$ give
$
 0\le \widetilde{h}(t)\le\frac1{\vartheta_1}=L.
$
Furthermore, \eqref{eq:g-lipschitz} implies
\begin{align*}
 |\widetilde{h}(t)-\widetilde{h}(u)|
 \le \int|g(t-R_2)-g(u-R_2)|\,\mathrm d P
 \le C_0|t-u|.
\end{align*}
Thus
\begin{equation}
 0\le \widetilde{h}\le L,\qquad
 |\widetilde{h}(t)-\widetilde{h}(u)|\le C_0|t-u|,\quad t,u\in\R .
 \label{eq:h-bounds}
\end{equation}
In particular, $\widetilde{h}$ is bounded and uniformly continuous. Since $R_2\ge0$
and $g$ vanishes on $(-\infty,0]$, \eqref{eq:h-representative} also gives
$
\widetilde{ h}(t)=0,$ for all $ t\le0 .
$

Lemma~\ref{lem:radialdensity} provides a constant $\ell>0$ such that the density $h$ constructed
there satisfies
\[
 h(t)\ge\ell,\qquad \frac12\le t\le1 .
\]
Since $\widetilde{h}=h$ almost everywhere, and therefore $\widetilde{h}$ satisfies this lower bound almost
everywhere on $[1/2,1]$. Continuity extends the bound to every point of
that interval, including the endpoints. In particular,
$
 \widetilde{h}(1)\ge\ell>0.
$

\medskip
\noindent\textbf{Step 2. Proof of
\eqref{eq:Jasymptotic}.}

Because $X$ has density $\widetilde{h}$,
\begin{equation*}
 J_n(s)=\int_0^1(1-x)^s x^n \widetilde{h}(x)\dd x.
 \label{eq:radial-integral}
\end{equation*}
These integrals are finite, since $0\le x^n\le1$, $0\le\widetilde{h}\le L$, and
$
 \int_0^1(1-x)^s\dd x=\frac1{s+1}<\infty.
$
For $n\ge1$, put $t\triangleq n(1-x)$. Then
\begin{equation}
 n^{s+1}J_n(s)
 =\int_0^n t^s\left(1-\frac tn\right)^n
 \widetilde{h}\left(1-\frac tn\right)\dd t.
 \label{eq:radial-rescaling}
\end{equation}
Extend the integrand by zero outside $(0,n)$:
\[
 G_n(t)\triangleq
 \begin{cases}
 \displaystyle t^s\left(1-\frac tn\right)^n
 \widetilde{h}\left(1-\frac tn\right),&0<t<n,\\[6pt]
 0,&t\ge n.
 \end{cases}
\]
For each fixed $t>0$, when $n\to\infty$,
\[
 \left(1-\frac tn\right)^n\longrightarrow e^{-t},
 \qquad
 \widetilde{h}\left(1-\frac tn\right)\longrightarrow \widetilde{h}(1).
\]
The second limit follows from continuity of $\widetilde{h}$ at $1$. Consequently,
\begin{equation}
 G_n(t)\longrightarrow \widetilde{h}(1)t^se^{-t}.
 \label{eq:radial-pointwise}
\end{equation}
For $0<t<n$, the inequality $\log(1-u)\le-u$ for $0<u<1$ gives
$
 \left(1-\frac tn\right)^n\le e^{-t}.
$
Hence
\[
 0\le G_n(t)\le Lt^se^{-t},\qquad t>0.
\]
The majorant is integrable because $s>-1$:
\[
 \int_0^\infty Lt^se^{-t}\dd t=L\Gamma(s+1)<\infty.
\]
Dominated convergence applied to
\eqref{eq:radial-rescaling}--\eqref{eq:radial-pointwise} yields
\[
 \lim_{n\to\infty}n^{s+1}J_n(s)
 =\widetilde{h}(1)\int_0^\infty t^se^{-t}\dd t
 =\widetilde{h}(1)\Gamma(s+1).
\]
Since $\widetilde{h}(1)>0$, this is precisely \eqref{eq:Jasymptotic}. The argument
also covers $-1<s<0$: although $t^s$ is unbounded near zero, it is
integrable there. 

\medskip
\noindent\textbf{Step 3. The change from exponential to gamma distributions.}

Fix $\eta=(\eta_j)_{j\in\mathbb N}\in \mathbb{N}_0^{(\mathbb{N})}$, and let
$
 K\triangleq \{j:\eta_j>0\}.
$
This is a finite set. For every integer $m\ge0$,
\begin{align*}
 \int X_j^m\,\mathrm d P
 =\int_0^\infty t^m\frac1{\vartheta_j}e^{-t/\vartheta_j}\dd t
 =\vartheta_j^m\int_0^\infty u^me^{-u}\dd u
  =m!\,\vartheta_j^m.
\end{align*}
Thus
\begin{equation}
 q_{j,m}(t)\triangleq
 \frac{t^m}{m!\,\vartheta_j^m}f_j(t)
 =\frac{t^me^{-t/\vartheta_j}}
 {m!\,\vartheta_j^{m+1}}\mathbf 1_{(0,\infty)}(t),\qquad t\in\mathbb R,
 \label{eq:gamma-density}
\end{equation}
is a probability density. It is the gamma density of shape $m+1$ and
scale $\vartheta_j$.

For completeness, $q_{j,m}$ is the density of a sum of $m+1$ independent
exponential variables of mean $\vartheta_j$. This follows by induction.
The case $m=0$ is $q_{j,0}=f_j$. For $t>0$,
\begin{align*}
 (q_{j,m}*f_j)(t)
 =\frac{e^{-t/\vartheta_j}}{m!\,\vartheta_j^{m+2}}
   \int_0^t u^m\dd u
 =\frac{t^{m+1}e^{-t/\vartheta_j}}
 {(m+1)!\,\vartheta_j^{m+2}}
 =q_{j,m+1}(t).
\end{align*}

On an enlarged probability space, introduce independent variables
\[
 E_{j,k},\qquad j\in K,\quad1\le k\le\eta_j,
\]
each $E_{j,k}$ exponential with mean $\vartheta_j$. Choose these variables
independently of the original sequence $\{X_j\}_{j=1}^{\infty}$, and set
\begin{equation}
 R_\eta\triangleq \sum_{j\in K}\sum_{k=1}^{\eta_j}E_{j,k}.
 \label{eq:extra-sum}
\end{equation}
This is a finite sum and is independent of $X$. When $\eta=\mathbf 0$, set
$R_0\triangleq 0$.

Equation~\eqref{eq:gamma-density} says that multiplying the density of
$X_j$ by $t^{\eta_j}$ produces the factor
$\eta_j!\vartheta_j^{\eta_j}$ and replaces that density by the density of
\[
 X_j+\sum_{k=1}^{\eta_j}E_{j,k}.
\]
Applying this identity in the finitely many coordinates $j\in K$ and
using independence gives
\begin{equation}
 \int G(X)\prod_jX_j^{\eta_j}
 =\eta!\,\vartheta^\eta\,\int G(X+R_\eta)
 \label{eq:change-density}
\end{equation}
for every nonnegative Borel function $G$ on $[0,\infty)$.


\medskip
\noindent\textbf{Step 4. Proof of \eqref{eq:momentupper}.}

Define the nonnegative Borel function
\begin{equation*}
 W_s(t)\triangleq
 \begin{cases}
 (1-t)^s,&0\le t<1,\\
 0,&t\ge1.
 \end{cases}
 \label{eq:weight}
\end{equation*}
By the definition of $M_\eta(s)$ at \eqref{2026111def1},
$
 M_\eta(s)=\int W_s(X)\prod_jX_j^{\eta_j}\,\mathrm d P.
$
Taking $G=W_s$ in \eqref{eq:change-density} gives
\begin{equation*}
 M_\eta(s)=\eta!\,\vartheta^\eta\,
 \int W_s(X+R_\eta).
 \label{eq:moment-shift}
\end{equation*}
For $r\ge0$, define
$
 F_s(r)\triangleq \int W_s(X+r) .
$
Because $X$ has density $\widetilde{h}$,
\begin{equation}
 F_s(r)=
 \begin{cases}
 \displaystyle
 \int_0^{1-r}(1-r-x)^s \widetilde{h}(x)\dd x,&0\le r<1,\\[8pt]
 0,&r\ge1.
 \end{cases}
 \label{eq:Fs-definition}
\end{equation}
Independence of $X$ and $R_\eta$ now yields the exact identity
\begin{equation}
 M_\eta(s)=\eta!\,\vartheta^\eta\,\int F_s(R_\eta).
 \label{eq:moment-representation}
\end{equation}
For $0\le r<1$, use $0\le \widetilde{h}\le L$ to obtain
\begin{align*}
 0\le F_s(r)
 \le L\int_0^{1-r}(1-r-x)^s\dd x
 =\frac{L}{s+1}(1-r)^{s+1}
 \le\frac{L}{s+1}.
 \label{eq:Fs-bound-sharp}
\end{align*}
The integral is finite because $s>-1$, and the last inequality holds
because $s+1>0$. For $r\ge1$, the same bound holds since $F_s(r)=0$.
Consequently,
\begin{equation*}
 0\le F_s(r)\le\frac{L}{s+1},\qquad r\ge0 .
 \label{eq:Fs-bound}
\end{equation*}
Substituting this into \eqref{eq:moment-representation} proves
$
 M_\eta(s)\le\frac{L}{s+1}\,\eta!\,\vartheta^\eta,
$
which is \eqref{eq:momentupper}.

\medskip
\noindent\textbf{Step 5. Continuity of $F_s$, including $-1<s<0$.}

The representation \eqref{eq:Fs-definition} has an upper endpoint $1-r$
that depends on $r$. If $s<0$, the factor $(1-r-x)^s$ is singular as
$x\uparrow1-r$. To treat this precisely, substitute
\[
 t\triangleq 1-r-x.
\]
Since $\widetilde{h}(x)=0$ for $x\le0$, the resulting formula can be written on a fixed
interval:
\begin{equation}
 F_s(r)=\int_0^1t^s \widetilde{h}(1-r-t)\dd t,\qquad r\ge0.
 \label{eq:Fs-fixed}
\end{equation}
If
\[
 k_s(t)\triangleq t^s\mathbf 1_{(0,1)}(t),\qquad t\in\mathbb R,
\]
then
\begin{equation}
 F_s(r)=(\widetilde{h}*k_s)(1-r),\qquad
 \|k_s\|_{L^1(\R)}=\frac1{s+1}.
 \label{eq:Fs-convolution}
\end{equation}

Continuity now follows directly. If $r_n\to r$ in $[0,\infty)$,
continuity of $\widetilde{h}$ gives
\[
 t^s\widetilde{h}(1-r_n-t)\longrightarrow t^s\widetilde{h}(1-r-t),
 \qquad 0<t<1 .
\]
Moreover,
$
 0\le t^s\widetilde{h}(1-r_n-t)\le Lt^s,
$
and $Lt^s$ is integrable on $(0,1)$ because $s>-1$. Dominated convergence
therefore yields $F_s(r_n)\to F_s(r)$. 

In fact, the Lipschitz estimate \eqref{eq:h-bounds} gives
\begin{align*}
 |F_s(r)-F_s(q)|
 \le\int_0^1t^s
 |\widetilde{h}(1-r-t)-\widetilde{h}(1-q-t)|\dd t
 \le C_0|r-q|\int_0^1t^s\dd t
  =\frac{C_0}{s+1}|r-q|.
\end{align*}
Thus
\begin{equation}
 |F_s(r)-F_s(q)|
 \le\frac{|r-q|}{(s+1)\vartheta_1\vartheta_2},
 \qquad r,q\ge0.
 \label{eq:Fs-lipschitz}
\end{equation}
Finally, setting $r=0$ in \eqref{eq:Fs-definition} gives
\begin{equation}
 F_s(0)=\int_0^1(1-x)^s\widetilde{h}(x)\dd x=I_s.
 \label{eq:Fs-zero}
\end{equation}

\medskip
\noindent\textbf{Step 6. Proof of \eqref{eq:firstmoment}.}

Take $\eta=\mathbf e_j$. Then $\eta!=1$ and $\vartheta^\eta=\vartheta_j$, while
$R_{\mathbf e_j}$ is a single exponential random variable of mean $\vartheta_j$.
Formula~\eqref{eq:moment-representation} becomes
\begin{equation*}
 \frac{M_{\mathbf e_j}(s)}{\vartheta_j}=\int F_s(R_{\mathbf e_j}) .
 \label{eq:coordinate-representation}
\end{equation*}
Since $\sum_j\vartheta_j<\infty$, we have $\vartheta_j\to0$. Using
\eqref{eq:Fs-lipschitz} and \eqref{eq:Fs-zero}, we find
\begin{align*}
 \left|\frac{M_{\mathbf e_j}(s)}{\vartheta_j}-I_s\right|
 &=\left|\int F_s(R_{\mathbf e_j})-F_s(0)\right|
 \le\int |F_s(R_{\mathbf e_j})-F_s(0)| \\
 &\le\frac{\int R_{\mathbf e_j}}{(s+1)\vartheta_1\vartheta_2}
 =\frac{\vartheta_j}{(s+1)\vartheta_1\vartheta_2}
 \longrightarrow0.
\end{align*}
This proves \eqref{eq:firstmoment}. 
All assertions of Theorem~\ref{lem:moments} follow.
\end{proof}
Set
\[
 Q(\mathbf z,\mathbf y)\triangleq\sum_{j\ge1}|z_j||y_j|,
 \qquad \mathbf z=(z_i)_{i\in\mathbb N}\in\overline B,\quad \mathbf y=(y_i)_{i\in\mathbb N}\in B.
\]
Then \(|\langle\mathbf z,\mathbf y\rangle|\le Q(\mathbf z,\mathbf y)<1\). The following estimate, which may be viewed as a refinement of Lemma~\ref{lem:caps}, will be used to establish the absolute convergence of the homogeneous expansions.

\begin{lemma} \label{lem:moduluscap}
For every \(D>0\) there is \(C_D<\infty\) such that
\[
 \sup_{\mathbf z\in\overline B}\mu(\{\mathbf y\in B:1-Q(\mathbf z,\mathbf y)<\varepsilon\} )
 \le C_D\varepsilon^D,\qquad 0<\varepsilon<1.
\]
Consequently, for \(s>-1\) and \(L_0\ge0\),
\begin{equation*}\label{eq:modulusmoment}
 \sup_{\mathbf z\in\overline B}\int_B
 \delta(\mathbf y)^s(1-Q(\mathbf z,\mathbf y))^{-L_0}\dd\mu(\mathbf y)<\infty.
\end{equation*}
\end{lemma}
\begin{proof}
Set $R_j \triangleq \sqrt{ X_j }$ for every $j\in\mathbb N$. For any Borel set $E\subset(0,\infty)$, polar
coordinates give
\begin{align*}
P(\{\mathbf z\in \ell^2:R_j(\mathbf z)\in E\})
&=\frac{1}{2\pi a_j^2}
\int_0^{2\pi}\int_E
\exp\!\left(-\frac{r^2}{2a_j^2}\right)r\dd r\dd\theta
=\int_E\frac{r}{a_j^2}
\exp\!\left(-\frac{r^2}{2a_j^2}\right)\dd r.
\end{align*}
Hence $R_j$ has the density
\[
g_j(r)=\frac{r}{a_j^2}
\exp\!\left(-\frac{r^2}{2a_j^2}\right)
\ind_{(0,\infty)}(r),\qquad r\in\mathbb R.
\]
The function $r\mapsto r\exp(-r^2/(2a_j^2))$ attains its maximum at
$r=a_j$. Therefore
\begin{equation}
\norm{g_j}_{L^\infty(\R)}=\frac{e^{-1/2}}{a_j}.
\label{eq:radius-density-bound}
\end{equation}

Fix an integer $N\ge1$. Since $X_1,\ldots,X_N$ are independent under
$P$, their square root $R_1,\ldots,R_N$ are also independent. Consequently,
the random vector
\[
R^{(N)}\triangleq(R_1,\ldots,R_N)\in\R^N
\]
has joint density
\[
F_N(r_1,\ldots,r_N)\triangleq \prod_{j=1}^N g_j(r_j),\qquad (r_1,\ldots,r_N)\in\mathbb R^N,
\]
and \eqref{eq:radius-density-bound} implies
\begin{equation*}
F_N(r_1,\ldots,r_N)
\le e^{-N/2}\prod_{j=1}^N a_j^{-1}.
\label{eq:joint-density-bound}
\end{equation*}
Let $\omega_N\triangleq  \frac{\pi^{N/2}}{\Gamma(N/2+1)}$ which is the Lebesgue measure of the open unit ball in $\mathbb R^N$.
For every center $\mathbf c\in\R^N$ and every radius $\rho>0$, we obtain
\begin{equation}
\begin{aligned}
P\left(\left\{\mathbf z\in\ell^2:\|R^{(N)}(\mathbf z)-\mathbf c\|<\rho\right\}\right)
&=\int_{\left\{\mathbf r\in\mathbb R^N:\|\mathbf r-\mathbf c\|<\rho\right\}}F_N(\mathbf r)\dd \mathbf r\\
&\le e^{-N/2}
\left(\prod_{j=1}^N a_j^{-1}\right)\omega_N\rho^N.
\end{aligned}
\label{eq:finite-ball-bound}
\end{equation}
The bound is independent of the center~$\mathbf c$.


For $\mathbf z=(z_i)_{i\in\mathbb N}\in\overline B$ and $0<\varepsilon<1$, put
\[
A_\varepsilon(\mathbf z)\triangleq\{\mathbf y\in B:1-Q(\mathbf z,\mathbf y)<\varepsilon\}.
\]
If $\mathbf y=(y_i)_{i\in\mathbb N}\in A_\varepsilon(\mathbf z)$, then
\begin{equation*}
\begin{aligned}
\sum_{j=1}^{\infty}\bigl(|y_j|-|z_j|\bigr)^2
&=\sum_{j=1}^{\infty}|y_j|^2
 +\sum_{j=1}^{\infty}|z_j|^2
 -2\sum_{j=1}^{\infty}|y_j||z_j|\\
&=\norm{\mathbf y}^2+\norm{\mathbf z}^2-2Q(\mathbf z,\mathbf y)\\
&\le 2-2Q(\mathbf z,\mathbf y)\\
&<2\varepsilon.
\end{aligned}
\label{eq:modulus-distance}
\end{equation*}
All these series converge because $(|y_j|)_{j\in\mathbb N},(|z_j|)_{j\in\mathbb N}\in\ell^2$.

In particular, for every integer $N\ge1$,
\[
\sum_{j=1}^N\bigl(|y_j|-|z_j|\bigr)^2<2\varepsilon.
\]
Writing $c_N(\mathbf z)\triangleq(|z_1|,\ldots,|z_N|)$, we therefore have the inclusion
\begin{equation*}
\begin{split}
A_\varepsilon(\mathbf z)\subset
\bigl\{\mathbf y=(y_j)_{j\in\mathbb N}\in\ell^2:
(|y_1|,\ldots,|y_N|)
\in B_{\R^N}\bigl(c_N(\mathbf z),\sqrt{2\varepsilon}\bigr)\bigr\},
\end{split}
\label{eq:cap-inclusion}
\end{equation*}
where $B_{\R^N}\bigl(c_N(\mathbf z),\sqrt{2\varepsilon}\bigl)\triangleq \left\{\mathbf r\in\mathbb R^N:\|\mathbf r-c_N(\mathbf z)\|<\sqrt{2\varepsilon}\right\}$.
Because $A_\varepsilon(\mathbf z)\subset B$ and $\mu=P|_B$, by \eqref{eq:finite-ball-bound}, it follows that
\begin{equation}
\begin{aligned}
\mu(A_\varepsilon(\mathbf z))
&=P(A_\varepsilon(\mathbf z))\\
&\le P\!\left(\left\{\mathbf y\in\ell^2:R^{(N)}(\mathbf y)
\in B_{\R^N}\bigl(c_N(\mathbf z),\sqrt{2\varepsilon}\bigr)\right\}\right)\\
&\le e^{-N/2}
\left(\prod_{j=1}^N a_j^{-1}\right)
\omega_N(2\varepsilon)^{N/2}\\
&=K_N\varepsilon^{N/2},
\end{aligned}
\label{eq:cap-N}
\end{equation}
where
\begin{equation*}
K_N\triangleq\omega_N\left(\frac2e\right)^{N/2}
\left(\prod_{j=1}^N a_j^{-1}\right)<\infty.
\label{eq:constant-N}
\end{equation*}


Now let $D>0$ be arbitrary. Choose an integer $N>2D$. Since
$0<\varepsilon<1$,
\[
\varepsilon^{N/2}\le\varepsilon^D.
\]
Thus \eqref{eq:cap-N} yields
\begin{equation}
\sup_{\mathbf z\in\overline B}\mu(A_\varepsilon(\mathbf z))
\le K_N\varepsilon^D,
\qquad 0<\varepsilon<1.
\label{eq:cap-D}
\end{equation}
We may therefore take $C_D\triangleq K_N$. The constant $C_D$ may depend on $D$ and on the fixed sequence $\{a_j\}_{j=1}^{\infty}$,
but it is independent of $\mathbf z$ or $\varepsilon$. No uniform bound on $K_N$ as
$N\to\infty$ is required: for each prescribed $D$, we choose one finite
$N$.


We next prove that, for every $\kappa\ge0$,
\begin{equation}
\mathcal M_\kappa
\triangleq\sup_{\mathbf z\in\overline B}
\int_B(1-Q(\mathbf z,\mathbf y))^{-\kappa}\dd\mu(\mathbf y)<\infty.
\label{eq:unweighted-moments}
\end{equation}

When $\kappa=0$, the integral equals $m$. Suppose that $\kappa>0$.
For every $t\in(0,1]$, elementary integration gives
\begin{equation*}
\begin{aligned}
t^{-\kappa}
=1+\kappa\int_t^1 u^{-\kappa-1}\dd u
=1+\kappa\int_0^1
u^{-\kappa-1}\ind_{\{t<u\}}\dd u.
\end{aligned}
\label{eq:layer-cake-scalar}
\end{equation*}
Apply this with $t=1-Q(\mathbf z,\mathbf y)$. Since the integrand on the right is
nonnegative, Tonelli's theorem gives
\begin{equation*}
\begin{aligned}
\int_B(1-Q(\mathbf z,\mathbf y))^{-\kappa}\dd\mu(\mathbf y)
&=m+\kappa\int_B\int_0^1
u^{-\kappa-1}\ind_{\{1-Q(\mathbf z,\mathbf y)<u\}}\dd u\dd\mu(\mathbf y)\\
&\quad=m+\kappa\int_0^1
u^{-\kappa-1}\mu(A_u(\mathbf z))\dd u.
\end{aligned}
\label{eq:layer-cake-integrated}
\end{equation*}
This equality is valid initially with values in $[0,\infty]$; thus we
have not assumed the integrability we are trying to prove.

Choose $D>\kappa$. Using \eqref{eq:cap-D}, we obtain
\begin{align*}
\int_B(1-Q(\mathbf z,\mathbf y))^{-\kappa}\dd\mu(\mathbf y)
\le m+\kappa C_D\int_0^1u^{D-\kappa-1}\dd u
=m+\frac{\kappa C_D}{D-\kappa}.
\end{align*}
The right-hand side is independent of $\mathbf z$. Hence
\begin{equation*}
\mathcal M_\kappa
\le m+\frac{\kappa C_D}{D-\kappa}<\infty,
\qquad D>\kappa.
\label{eq:moment-D-bound}
\end{equation*}


We first recall the radial integrability estimate used here. If
$t>-1$, then by Corollary \ref{cor:radial},
\begin{equation}
I_t=\int_B\delta(\mathbf y)^t\dd\mu(\mathbf y)<\infty.
\label{eq:radial-integrability}
\end{equation}

Now fix $s>-1$ and $L_0\ge0$. Choose $p_0>1$ such that
\begin{equation}
sp_0>-1,
\label{eq:holder-choice}
\end{equation}
and let $p_0'=p_0/(p_0-1)$.

Such a choice always exists:
\begin{itemize}
\item If $s\ge0$, we may take $p_0=2$.
\item If $-1<s<0$, condition~\eqref{eq:holder-choice} is equivalent to
$
1<p_0<-\frac1s.
$
This interval is nonempty because $s>-1$. For example,
$
p_0=\frac12\left(1-\frac1s\right)
$
satisfies
$
sp_0=\frac{s-1}{2}>-1.
$
\end{itemize}

H\"older's inequality now gives, for every $\mathbf z\in\overline B$,
\begin{equation*}
\begin{aligned}
\int_B\delta(\mathbf y)^s(1-Q(\mathbf z,\mathbf y))^{-L_0}\dd\mu(\mathbf y)
&\le
\left(\int_B\delta(\mathbf y)^{sp_0}\dd\mu(\mathbf y)\right)^{1/p_0}
\left(\int_B(1-Q(\mathbf z,\mathbf y))^{-L_0p_0'}\dd\mu(\mathbf y)\right)^{1/p_0'}\\
&\le I_{sp_0}^{1/p_0}\,
\mathcal M_{L_0p_0'}^{1/p_0'}.
\end{aligned}
\label{eq:holder-estimate}
\end{equation*}
The first factor is finite by \eqref{eq:radial-integrability}, because
$sp_0>-1$. The second factor is finite by
\eqref{eq:unweighted-moments}, and its bound is uniform in $\mathbf z$.
Taking the supremum proves
\[
\sup_{\mathbf z\in\overline B}
\int_B\delta(\mathbf y)^s(1-Q(\mathbf z,\mathbf y))^{-L_0}\dd\mu(\mathbf y)<\infty.
\]
 The proof is complete.
\end{proof}
Fixed $u,\, v>-\frac12$, and for all integers \(n,m\ge0\),  the integral operator $R_{n,m}$ on
$L^2(B,\mu)$ is defined by
\begin{equation}
(R_{n,m}f)(\mathbf z)
\triangleq\int_B
K_{n,m}(\mathbf z,\mathbf y)f(\mathbf y)\dd\mu(\mathbf y),\qquad \mathbf z\in B,\quad f\in L^2(B,\mu),
\label{eq:operator-definition}
\end{equation}
where
$
K_{n,m}(\mathbf z,\mathbf y)\triangleq \delta(\mathbf z)^u\delta(\mathbf y)^v
\ip{\mathbf z}{\mathbf y}^{n}\overline{\ip{\mathbf z}{\mathbf y}}^{\,m}\text{ for every } \mathbf y\in B.
$ We call $K_{n,m}$ the kernel of $R_{n,m}$.
\begin{lemma}
$R_{n,m}$ is a bounded linear operator from $L^2(B,\mu)$ into itself.
\end{lemma}
\begin{proof}

Set $d\triangleq n+m$. Recall \eqref{eq:Jdefinition},
\[
J_d(s)=\int_B\delta(\mathbf y)^s\norm{\mathbf y}^{2d}\dd\mu(\mathbf y),\qquad s>-1.
\]
Since $\norm{\mathbf y}<1$ for $\mathbf y\in B$,
\[
J_d(s)\le\int_B\delta(\mathbf y)^s\dd\mu(\mathbf y)<\infty,
\qquad s>-1,
\]
where the last assertion follows from Corollary~\ref{cor:radial}.
Thus
$J_d(2u)<\infty,$ and $J_d(2v)<\infty.
$

For each $\mathbf z,\mathbf y\in B$, the Cauchy--Schwarz inequality in $\ell^2$ gives
$
|\ip{\mathbf z}{\mathbf y}|\le\norm{\mathbf z}\norm{\mathbf y}.
$
Therefore,
\begin{align*}
|(R_{n,m}f)(\mathbf z)|\leq \int_B|K_{n,m}(\mathbf z,\mathbf y)f(\mathbf y)|\dd\mu(\mathbf y)
&=\delta(\mathbf z)^u\int_B\delta(\mathbf y)^v
|\ip{\mathbf z}{\mathbf y}|^{n+m}|f(\mathbf y)|\dd\mu(\mathbf y)\\
&\le\delta(\mathbf z)^u\norm{\mathbf z}^{d}
\int_B\delta(\mathbf y)^v\norm{\mathbf y}^{d}|f(\mathbf y)|\dd\mu(\mathbf y).
\end{align*}
Applying the Cauchy--Schwarz inequality in $L^2(B,\mu)$, we obtain
\begin{align*}
\int_B\delta(\mathbf y)^v\norm{\mathbf y}^{d}|f(\mathbf y)|\dd\mu(\mathbf y)
\le\left(\int_B\delta(\mathbf y)^{2v}\norm{\mathbf y}^{2d}
\dd\mu(\mathbf y)\right)^{1/2}\norm{f}_{L^2(B,\mu)}
=J_d(2v)^{1/2}\norm{f}_{L^2(B,\mu)}.
\end{align*}
Hence for every $\mathbf z\in B$, we have
\begin{equation}
|(R_{n,m}f)(\mathbf z)|
\le\delta(\mathbf z)^u\norm{\mathbf z}^{d}J_d(2v)^{1/2}\norm{f}_{L^2(B,\mu)}<\infty.
\label{eq:absolute-convergence}
\end{equation}
When $d=0$, the expression $\norm{\mathbf z}^{d}$ means~$1$, including at
$\mathbf z=\mathbf 0$.

Thus the integral converges absolutely for every $\mathbf z\in B$. Changing
$f$ on a $\mu$-null set leaves the integral unchanged, so the
definition respects the equivalence classes in $L^2(B,\mu)$.

Integrating \eqref{eq:absolute-convergence} with respect to $\mathbf z\in B$ gives
\begin{align*}
\norm{R_{n,m}f}_{L^2(B,\mu)}^2
\le J_d(2v)\norm{f}_{L^2(B,\mu)}^2
\int_B\delta(\mathbf z)^{2u}\norm{\mathbf z}^{2d}\dd\mu(\mathbf z)
=J_d(2u)J_d(2v)\norm{f}_{L^2(B,\mu)}^2<\infty.
\end{align*}
Consequently,
\begin{equation*}
\norm{R_{n,m}}_{\text{op}}
\le\sqrt{J_{n+m}(2u)J_{n+m}(2v)}<\infty,
\label{eq:operator-bound}
\end{equation*}
where $\norm{\cdot}_{\text{op}}$ denotes the operator norm from $L^2(B,\mu)$ into itself. This verifies that \eqref{eq:operator-definition} defines a bounded linear operator from $L^2(B,\mu)$ into itself.
\end{proof}
Furthermore, Lemma~\ref{lem:pieces} subsequently proves the stronger
assertion that $R_{n,m}$ is trace class.

 For vectors \(f,g\in L^2(B,\mu)\), the notation
\(f\otimes g\) means the rank-one operator
\[
h\mapsto\left(\int_B h \overline{g}\,\mathrm d \mu\right) f,\qquad h\in L^2(B,\mu).
\]
We shall use the identities
\begin{equation}
\norm{f\otimes g}_{\mathcal S_1}=\norm{f}_{L^2(B,\mu)}\norm{g}_{L^2(B,\mu)},
\qquad
\operatorname{Tr}(f\otimes g)=\ip{f}{g}_{L^2(B,\mu)}.
\label{eq:rank-one-identities}
\end{equation}

For the norm formula, the only possible nonzero singular value of
$f\otimes g$ is $\norm{f}_{L^2(B,\mu)}\norm{g}_{L^2(B,\mu)}$. Indeed,
\[
(f\otimes g)^*(f\otimes g)=\norm{f}_{L^2(B,\mu)}^2(g\otimes g),
\]
whose only possible nonzero eigenvalue is
$\norm{f}_{L^2(B,\mu)}^2\norm{g}_{L^2(B,\mu)}^2$.

For the trace formula, if $\{e_k\}_{k=1}^{\infty}$ is an orthonormal basis of $L^2(B,\mu)$, then
\begin{align*}
\operatorname{Tr}(f\otimes g)
=\sum_k\ip{(f\otimes g)e_k}{e_k}_{L^2(B,\mu)}
=\sum_k\ip{e_k}{g}_{L^2(B,\mu)}\ip{f}{e_k}_{L^2(B,\mu)}
=\ip{f}{g}_{L^2(B,\mu)}.
\end{align*}
We also need the following scalar inequality.
\begin{lemma}
For every $L>0$, there is a constant $A_L<\infty$
such that
\begin{equation}
x^{2d}\le A_L(1+d)^{-L}(1-x)^{-L},
\qquad d\in\N_0,\quad0\le x<1.
\label{eq:scalar-estimate}
\end{equation}
\end{lemma}
\begin{proof}
One may take
\[
A_L\triangleq\max\left\{1,\left(\frac Le\right)^L\right\}.
\]
If $d=0$, the left-hand side equals~$1$, while
$
A_L(1-x)^{-L}\ge1,
$
so the assertion follows.

Suppose that $d\ge1$. For $0<x<1$,
$
\log x\le-(1-x),
$
and hence
$
x^{2d}\le e^{-2d(1-x)}.
$
Therefore
\[
x^{2d}(1-x)^L\le e^{-2d(1-x)}(1-x)^L.
\]
Set $r=2d(1-x)$. Then
$
e^{-2d(1-x)}(1-x)^L=(2d)^{-L}r^Le^{-r}.
$
The function $r\mapsto r^Le^{-r}$ attains its maximum at $r=L$, so
$
r^Le^{-r}\le\left(\frac Le\right)^L.
$
Thus
\[
x^{2d}(1-x)^L
\le(2d)^{-L}\left(\frac Le\right)^L
\le(1+d)^{-L}\left(\frac Le\right)^L,
\]
where the last inequality uses $2d\ge1+d$. At $x=0$, the assertion
is immediate. This proves \eqref{eq:scalar-estimate}.
\end{proof}
For any $K\in L^2(B\times B,\mu\times\mu)$, its integral operator $T_K$ is defined by
$$
T_Kf(\mathbf z)\triangleq \int_Bf(\mathbf y)K(\mathbf z,\mathbf y)\,\mathrm d \mu(\mathbf y),\qquad f\in L^2(B,\mu).
$$
Indeed, Cauchy--Schwarz inequality gives
\begin{align*}
\norm{T_Kf}_{L^2(B,\mu)}^2
=\int_B\left|\int_BK(\mathbf z,\mathbf y)f(\mathbf y)\dd\mu(\mathbf y)\right|^2\dd\mu(\mathbf z)
\le\norm{f}_{L^2(B,\mu)}^2\int_B\int_B|K(\mathbf z,\mathbf y)|^2\dd\mu(\mathbf y)\dd\mu(\mathbf z).
\end{align*}
Thus $T_K$ is a bounded linear operator from $L^2(B,\mu)$ into itself with
\begin{equation}
\norm{T_K}_{\mathrm{op}}\le\norm{K}_{L^2(B\times B,\mu\times\mu)},
\label{eq:L2-operator-bound23}
\end{equation}
where $\norm{T_K}_{\mathrm{op}}$ denotes the operator norm of $T_K$ from $L^2(B,\mu)$ into itself.
\begin{lemma}\label{20261113lem2}
For any kernel $K\in L^2(B\times B,\mu\times\mu)$, its integral operator $T_K$ satisfies
\begin{equation}
T_K\in\mathcal S_2,\qquad \|T_K\|_{\mathcal S_2}=\norm{K}_{L^2(B\times B,\mu\times\mu)},\qquad\norm{T_K}_{\mathrm{op}}\le\norm{K}_{L^2(B\times B,\mu\times\mu)}.
\label{eq:L2-operator-bound}
\end{equation}
\end{lemma}
\begin{proof}
Firstly, the inequality is from \eqref{eq:L2-operator-bound23}. Let $\{e_k\}_{k=1}^{\infty}$ be an orthonormal basis of $L^2(B,\mu)$. For almost
every~$\mathbf z\in B$,
\[
(T_{K}e_k)(\mathbf z)
=\left\langle e_k,\overline{K (\mathbf z,\cdot)}\right\rangle_{L^2(B,\mu)}.
\]
Hence
\[
\sum_k|(T_{K}e_k)(\mathbf z)|^2
=\int_B|K (\mathbf z,\mathbf y)|^2\dd\mu(\mathbf y).
\]
Integrating in $\mathbf z\in B$ and applying Tonelli's theorem gives
\eqref{eq:L2-operator-bound}.

\end{proof}
\begin{lemma}\label{lem:pieces}
Every \(R_{n,m}\) is trace class and
\begin{equation}\label{eq:piecetracebound}
 \|R_{n,m}\|_{\mathcal S_1}
 \le\sqrt{J_{n+m}(2u)J_{n+m}(2v)},\qquad
 \operatorname{Tr}R_{n,m}=J_{n+m}(t).
\end{equation}
For every \(D>0\),
\begin{equation}\label{eq:piecerapid}
 \|R_{n,m}\|_{\mathcal S_2}
 \le C_D(1+n+m)^{-D}.
\end{equation}
\end{lemma}

\begin{proof}
Put \(d=n+m\) and \(t=u+v\).
For every $N\in\mathbb N,$ $\mathbf z=(z_i)_{i\in\mathbb N},\mathbf y=(y_i)_{i\in\mathbb N}\in B$, let
\[
P_N\mathbf z\triangleq \sum_{j=1}^Nz_j\mathbf e_j,\qquad
q_N(\mathbf z,\mathbf y)\triangleq \ip{P_N\mathbf z}{P_N\mathbf y}=\sum_{j=1}^Nz_j\overline{y_j},
\]
and
\[
K_{n,m}^N(\mathbf z,\mathbf y)\triangleq\delta(\mathbf z)^u\delta(\mathbf y)^v
q_N(\mathbf z,\mathbf y)^n\overline{q_N(\mathbf z,\mathbf y)}^{\,m}.
\]
Let $R_{n,m}^N$ be the integral operator with kernel $K_{n,m}^N$. Put
\[
\begin{split}
E_N\triangleq\left\{(\eta,\zeta)\in \mathbb{N}_0^{(\mathbb{N})}\times \mathbb{N}_0^{(\mathbb{N})}:
 \ |\eta|=n,\ |\zeta|=m, \,
 \operatorname{supp}\eta\cup\operatorname{supp}\zeta
\subset\{1,\ldots,N\}\right\}.
\end{split}
\]
This is a finite set. For $(\eta,\zeta)\in E_N$ and $\mathbf z\in B$, define
\[
c_{\eta,\zeta}\triangleq\frac{n!m!}{\eta!\zeta!},
\qquad F_{\eta,\zeta}(\mathbf z)\triangleq \delta(\mathbf z)^u \mathbf z^\eta\overline{\mathbf z}^{\,\zeta},
\qquad
G_{\eta,\zeta}(\mathbf z)\triangleq \delta(\mathbf z)^v \mathbf z^\eta\overline{\mathbf z}^{\,\zeta}.
\]
Both functions belong to $L^2(B,\mu)$, since
\[
\norm{F_{\eta,\zeta}}_{L^2(B,\mu)}^2
=M_{\eta+\zeta}(2u)\le I_{2u}<\infty,
\]
and similarly
\begin{equation}
\norm{G_{\eta,\zeta}}_{L^2(B,\mu)}^2
=M_{\eta+\zeta}(2v)\le I_{2v}<\infty.
\label{eq:vector-norms}
\end{equation}

The finite multinomial theorem gives
\[
q_N(\mathbf z,\mathbf y)^n\overline{q_N(\mathbf z,\mathbf y)}^{\,m}
=\sum_{(\eta,\zeta)\in E_N}
c_{\eta,\zeta}\mathbf z^\eta\overline{\mathbf z}^{\,\zeta}
\overline {\mathbf y}^{\,\eta}\mathbf y^\zeta.
\]
Moreover,
\[
F_{\eta,\zeta}(\mathbf z)\overline{G_{\eta,\zeta}(\mathbf y)}
=\delta(\mathbf z)^u\delta(\mathbf y)^v
\mathbf z^\eta\overline {\mathbf z}^{\,\zeta}\overline {\mathbf y}^{\,\eta}\mathbf y^\zeta.
\]
Thus
\begin{equation}
R_{n,m}^N\triangleq\sum_{(\eta,\zeta)\in E_N}
c_{\eta,\zeta}F_{\eta,\zeta}\otimes G_{\eta,\zeta}.
\label{eq:finite-rank-expansion}
\end{equation}
In particular, every $R_{n,m}^N$ has finite rank.


The key numerical identity is
\begin{equation}
\sum_{\substack{|\eta|=n\\|\zeta|=m}}
c_{\eta,\zeta}M_{\eta+\zeta}(s)=J_d(s),\qquad s>-1.
\label{eq:moment-sum}
\end{equation}
Let us prove this carefully.

For finite $N$, the multinomial theorem gives
\begin{equation}
\begin{aligned}
\sum_{(\eta,\zeta)\in E_N}
c_{\eta,\zeta}M_{\eta+\zeta}(s)
&=\int_B\delta(\mathbf z)^s
\sum_{(\eta,\zeta)\in E_N}
\frac{n!m!}{\eta!\zeta!}|\mathbf z^\eta|^2|\mathbf z^\zeta|^2\dd\mu(\mathbf z)\\
&\quad=\int_B\delta(\mathbf z)^s
\left(\sum_{j=1}^N|z_j|^2\right)^n
\left(\sum_{j=1}^N|z_j|^2\right)^m\dd\mu(\mathbf z)\\
&\quad=\int_B\delta(\mathbf z)^s\norm{P_N\mathbf z}^{2d}\dd\mu(\mathbf z).
\end{aligned}
\label{eq:finite-moment-sum}
\end{equation}
All summands are nonnegative. Also,
\[
\norm{P_N\mathbf z}^2\uparrow\norm{\mathbf z}^2,\qquad \text{as}\,N\to\infty.
\]
Therefore monotone convergence, together with
\[
\bigcup_{N=1}^{\infty}E_N
=\left\{(\eta,\zeta)\in \mathbb{N}_0^{(\mathbb{N})}\times \mathbb{N}_0^{(\mathbb{N})}:|\eta|=n,\ |\zeta|=m\right\},
\]
proves \eqref{eq:moment-sum}. When $d=0$, every zeroth power is understood as~$1$, and the same
argument applies.

Now use \eqref{eq:rank-one-identities} and \eqref{eq:vector-norms}:
\begin{align*}
\sum_{\substack{|\eta|=n\\|\zeta|=m}}
\norm{c_{\eta,\zeta}F_{\eta,\zeta}\otimes G_{\eta,\zeta}}_{\mathcal S_1}
=\sum_{\substack{|\eta|=n\\|\zeta|=m}}
c_{\eta,\zeta}
\sqrt{M_{\eta+\zeta}(2u)M_{\eta+\zeta}(2v)}.
\end{align*}
Apply Cauchy--Schwarz inequality to the numerical sequences
\[
\left(\sqrt{c_{\eta,\zeta}M_{\eta+\zeta}(2u)}\right)_{\eta,\zeta},
\qquad
\left(\sqrt{c_{\eta,\zeta}M_{\eta+\zeta}(2v)}\right)_{\eta,\zeta}.
\]
Using \eqref{eq:moment-sum}, we obtain
\begin{equation}
\begin{aligned}
&\sum_{\substack{|\eta|=n\\|\zeta|=m}}
c_{\eta,\zeta}\sqrt{M_{\eta+\zeta}(2u)M_{\eta+\zeta}(2v)}\\
&\qquad\le
\left(\sum_{\substack{|\eta|=n\\|\zeta|=m}}
c_{\eta,\zeta}M_{\eta+\zeta}(2u)\right)^{1/2}
\left(\sum_{\substack{|\eta|=n\\|\zeta|=m}}
c_{\eta,\zeta}M_{\eta+\zeta}(2v)\right)^{1/2}
=\sqrt{J_d(2u)J_d(2v)}<\infty.
\end{aligned}
\label{eq:trace-norm-sum}
\end{equation}

Since $\mathcal S_1$ is complete, \eqref{eq:trace-norm-sum} proves
that the operator series
\[
A=\sum_{\substack{|\eta|=n\\|\zeta|=m}}
c_{\eta,\zeta}F_{\eta,\zeta}\otimes G_{\eta,\zeta}
\]
converges absolutely in trace norm. Moreover,
\begin{equation}
\norm{A}_{\mathcal S_1}\le\sqrt{J_d(2u)J_d(2v)}.
\label{eq:limit-trace-bound}
\end{equation}
The finite sums $R_{n,m}^N$ converge to $A$ in trace norm, because
\begin{equation}
\norm{A-R_{n,m}^N}_{\mathcal S_1}
\le\sum_{\substack{|\eta|=n,\ |\zeta|=m\\(\eta,\zeta)\notin E_N}}
c_{\eta,\zeta}\norm{F_{\eta,\zeta}}_{L^2(B,\mu)}\norm{G_{\eta,\zeta}}_{L^2(B,\mu)}
\longrightarrow0.
\label{eq:trace-norm-convergence}
\end{equation}
This establishes actual convergence of the series. 


We must still verify that the operator $A$ constructed above is the
original integral operator $R_{n,m}$.

For every $\mathbf z,\mathbf y\in B$,
$
q_N(\mathbf z,\mathbf y)\longrightarrow\ip{\mathbf z}{\mathbf y}.
$
Hence
$
K_{n,m}^N(\mathbf z,\mathbf y)\longrightarrow K_{n,m}(\mathbf z,\mathbf y).
$
Furthermore,
\[
|q_N(\mathbf z,\mathbf y)|\le\norm{P_N\mathbf z}\norm{P_N\mathbf y}
\le\norm{\mathbf z}\norm{\mathbf y}<1,
\]
and the same bound holds for $|\ip{\mathbf z}{\mathbf y}|$. Therefore
\[
|K_{n,m}^N(\mathbf z,\mathbf y)|\le\delta(\mathbf z)^u\delta(\mathbf y)^v,\qquad
|K_{n,m}(\mathbf z,\mathbf y)|\le\delta(\mathbf z)^u\delta(\mathbf y)^v.
\]
It follows that
\begin{equation}
|K_{n,m}^N(\mathbf z,\mathbf y)-K_{n,m}(\mathbf z,\mathbf y)|^2
\le4\delta(\mathbf z)^{2u}\delta(\mathbf y)^{2v}.
\label{eq:kernel-majorant}
\end{equation}
The majorant is integrable on $B\times B$, since $2u>-1,2v>-1,$ and by Corollary \ref{cor:radial},
\[
\int_B\int_B\delta(\mathbf z)^{2u}\delta(\mathbf y)^{2v}
\dd\mu(\mathbf y)\dd\mu(\mathbf z)=I_{2u}I_{2v}<\infty.
\]
Dominated convergence consequently yields
\begin{equation}
\norm{K_{n,m}^N-K_{n,m}}_{L^2(B\times B,\mu\times\mu)}\longrightarrow0.
\label{eq:L2-kernel-limit}
\end{equation}

Applying \eqref{eq:L2-operator-bound} to $K=K_{n,m}^N-K_{n,m}$, we see
from \eqref{eq:L2-kernel-limit} that
\[
\|R_{n,m}^N- R_{n,m}\|_{\text{op}}\to 0.
\]
On the other hand, \eqref{eq:trace-norm-convergence} and
$
\norm{A-R_{n,m}^N}_{\text{op}}\le\norm{A-R_{n,m}^N}_{\mathcal S_1}
$
show that
\[
\|R_{n,m}^N-A\|_{\text{op}}\to 0.
\]
Uniqueness of the limit gives $A=R_{n,m}$.

Thus the countable expansion
\begin{equation}
R_{n,m}=\sum_{\substack{|\eta|=n\\|\zeta|=m}}
\frac{n!m!}{\eta!\zeta!}F_{\eta,\zeta}\otimes G_{\eta,\zeta}
\label{eq:infinite-rank-one-expansion}
\end{equation}
is valid in trace norm. From \eqref{eq:limit-trace-bound},
$
\norm{R_{n,m}}_{\mathcal S_1}
\le\sqrt{J_{n+m}(2u)J_{n+m}(2v)}.
$


For each pair $(\eta,\zeta)$, the rank-one trace formula gives
\begin{equation*}
\begin{aligned}
\operatorname{Tr}(F_{\eta,\zeta}\otimes G_{\eta,\zeta})
&=\ip{F_{\eta,\zeta}}{G_{\eta,\zeta}}_{L^2(B,\mu)}\\
&=\int_B\delta(\mathbf z)^{u+v}
\bigl(\mathbf z^\eta\overline {\mathbf z}^{\,\zeta}\bigr)
\overline{\bigl(\mathbf z^\eta\overline {\mathbf z}^{\,\zeta}\bigr)}\dd\mu(\mathbf z)\\
&=\int_B\delta(\mathbf z)^t|\mathbf z^{\eta+\zeta}|^2\dd\mu(\mathbf z)\\
&=M_{\eta+\zeta}(t),
\end{aligned}
\label{eq:rank-one-trace-moment}
\end{equation*}
where $t=u+v$. Using the finite expansion~\eqref{eq:finite-rank-expansion}, followed
by \eqref{eq:finite-moment-sum} with $s=t$, we obtain
\begin{equation*}
\begin{aligned}
\operatorname{Tr}R_{n,m}^N
=\sum_{(\eta,\zeta)\in E_N}c_{\eta,\zeta}M_{\eta+\zeta}(t)
=\int_B\delta(\mathbf z)^t\norm{P_N\mathbf z}^{2d}\dd\mu(\mathbf z).
\end{aligned}
\label{eq:finite-trace}
\end{equation*}
Since $t>-1$, monotone convergence gives
\begin{equation}
\operatorname{Tr}R_{n,m}^N\longrightarrow
\int_B\delta(\mathbf z)^t\norm{\mathbf z}^{2d}\dd\mu(\mathbf z)=J_d(t).
\label{eq:trace-moment-limit}
\end{equation}
At the same time, the trace is continuous in trace norm:
\[
\bigl|\operatorname{Tr}R_{n,m}^N-\operatorname{Tr}R_{n,m}\bigr|
\le\norm{R_{n,m}^N-R_{n,m}}_{\mathcal S_1}\longrightarrow0.
\]
Combining this with \eqref{eq:trace-moment-limit} proves
$
 \operatorname{Tr}R_{n,m}=J_d(t)=J_{n+m}(t).
$

This argument also justifies the diagonal formula in this particular
case:
\[
\operatorname{Tr}R_{n,m}=\int_BK_{n,m}(\mathbf z,\mathbf z)\dd\mu(\mathbf z),
\]
because
$
K_{n,m}(\mathbf z,\mathbf z)=\delta(\mathbf z)^t\norm{\mathbf z}^{2(n+m)}.
$

Recall for $\mathbf z=(z_i)_{i\in\mathbb N},\mathbf y=(y_i)_{i\in\mathbb N}\in B$,
$
Q(\mathbf z,\mathbf y)=\sum_{j=1}^{\infty}|z_j||y_j|,
$
and
$
|\ip{\mathbf z}{\mathbf y}|\le Q(\mathbf z,\mathbf y)<1.
$
Consequently, \eqref{eq:L2-operator-bound} implies
\begin{equation}
\norm{R_{n,m}}_{\mathcal S_2}^2
\le\int_B\int_B\delta(\mathbf z)^{2u}\delta(\mathbf y)^{2v}
Q(\mathbf z,\mathbf y)^{2d}\dd\mu(\mathbf y)\dd\mu(\mathbf z).
 \label{eq:HS-Q-bound}
\end{equation}

Fix $L>0$. By Lemma \ref{lem:moduluscap}, since $2v>-1$,
\begin{equation*}
\mathcal W_{2v,L}
\triangleq\sup_{\mathbf z\in\overline B}
\int_B\delta(\mathbf y)^{2v}(1-Q(\mathbf z,\mathbf y))^{-L}\dd\mu(\mathbf y)<\infty.
 \label{eq:weighted-modulus-moment}
\end{equation*}
Applying \eqref{eq:scalar-estimate} with $x=Q(\mathbf z,\mathbf y)$, we have
\[
Q(\mathbf z,\mathbf y)^{2d}
\le A_L(1+d)^{-L}(1-Q(\mathbf z,\mathbf y))^{-L}.
\]
Substituting this into \eqref{eq:HS-Q-bound} gives
\begin{equation*}
\begin{aligned}
\norm{R_{n,m}}_{\mathcal S_2}^2
&\le A_L(1+d)^{-L}\int_B\delta(\mathbf z)^{2u}\\
&\qquad\times\left[
\int_B\delta(\mathbf y)^{2v}(1-Q(\mathbf z,\mathbf y))^{-L}\dd\mu(\mathbf y)
\right]\dd\mu(\mathbf z)\\
&\le A_L(1+d)^{-L}\mathcal W_{2v,L}
\int_B\delta(\mathbf z)^{2u}\dd\mu(\mathbf z)\\
&=A_L\mathcal W_{2v,L}I_{2u}(1+d)^{-L}.
\end{aligned}
 \label{eq:HS-squared-decay}
\end{equation*}
Taking square roots, we obtain
\begin{equation*}
\norm{R_{n,m}}_{\mathcal S_2}
\le\bigl(A_L\mathcal W_{2v,L}I_{2u}\bigr)^{1/2}
(1+d)^{-L/2}.
\label{eq:HS-decay-L}
\end{equation*}

Finally, given any prescribed $D>0$, choose $L=2D$. Then
\[
\norm{R_{n,m}}_{\mathcal S_2}
\le C_D(1+n+m)^{-D},
\]
where $C_D\triangleq\bigl(A_{2D}\mathcal W_{2v,2D}I_{2u}\bigr)^{1/2}<\infty.$
The constant depends on $D$, the fixed parameters $u,v$, and the
Gaussian measure, but is independent of $n,m$. This proves
\eqref{eq:piecerapid}, and completes
the proof of Lemma~\ref{lem:pieces}.

\end{proof}

For $n\in\mathbb N_0$ and $d\in\mathbb R$, write
\begin{equation}\label{20261112for2}
A_n(d)\triangleq (d)_n/n!,\qquad \Lambda(d)\triangleq\sum_{n=0}^\infty|A_n(d)|,
\end{equation}
with \((d)_0\triangleq 1\) and
\((d)_n\triangleq d(d+1)\cdots(d+n-1)\). We also need the following estimates in this paper.
\begin{lemma}\label{20261112lem1}
If \(d\in\mathbb R\setminus\{0,-1,-2,\ldots\}\), then
\begin{equation}\label{eq:negative-asymptotic235}
|A_n(d)|
\sim\frac{n^{d-1}}{|\Gamma(d)|},
\qquad\text{as}\,\,\,n\to\infty.
\end{equation}
If \(d=-N\) for an integer \(N\ge1\), then
\begin{equation}\label{eq:terminating-coefficients35}
A_n(-N)=
\begin{cases}
(-1)^n\displaystyle\binom{N}{n},&0\le n\le N,\\[0.5em]
0,&n\ge N+1.
\end{cases}
\end{equation}
Moreover,
$\Lambda(d)<\infty
 \text{ for every }d<0.
$
\end{lemma}
\begin{proof}
For \(d>0\), every factor in \((d)_n= d(d+1)\cdots(d+n-1)\) is positive, and hence
\[
A_n(d)=\frac{ d(d+1)\cdots(d+n-1)}{n!}>0,\qquad n\ge0 .
\]
The Pochhammer--gamma identity (\cite[last display, in the proof of Theorem 1.2.7, p. 13]{AAR}) and the gamma-ratio asymptotic (\cite[equation (1.4.3), p. 20]{AAR}) give
\begin{equation*}\label{eq:positive-coeff-asymptotic2}
A_n(d)
=\frac{\Gamma(n+d)}{\Gamma(d)\Gamma(n+1)}
\sim\frac{n^{d-1}}{\Gamma(d)},
\qquad\text{as}\,\,\,n\to\infty.
\end{equation*}

Suppose \(d=-N\) for an integer \(N\ge 0\). Then
\begin{equation*}\label{eq:terminating-coefficients}
A_n(-N)=
\begin{cases}
(-1)^n\displaystyle\binom{N}{n},&0\le n\le N,\\[0.5em]
0,&n\ge N+1.
\end{cases}
\end{equation*}
Hence
\begin{equation*}\label{eq:terminating-l1}
\Lambda(-N)=\sum_{n=0}^N\binom{N}{n}=2^N<\infty.
\end{equation*}

Suppose \(d<0\) and \(d\notin\{-1,-2,\ldots\}\). Then \(\Gamma(d)\)
is finite and nonzero. The Pochhammer--gamma identity(\cite[last display, in the proof of Theorem 1.2.7, p. 13]{AAR}) gives
\[
A_n(d)=\frac{\Gamma(n+d)}{\Gamma(d)\Gamma(n+1)}.
\]
For all sufficiently large \(n\), we have \(n+d>0\), and therefore
\(\Gamma(n+d)>0\). Applying the gamma-ratio asymptotic(\cite[equation (1.4.3), p. 20]{AAR}) with
parameters \(d\) and \(1\), we get
\begin{equation*}\label{eq:negative-asymptotic2}
|A_n(d)|
=\frac{\Gamma(n+d)}{|\Gamma(d)|\Gamma(n+1)}
\sim\frac{n^{d-1}}{|\Gamma(d)|},
\qquad\text{as}\,\,\,n\to\infty.
\end{equation*}
In particular, there are constants \(C_d<\infty\) and \(n_d\in\mathbb N\) such that
\[
|A_n(d)|\le C_d n^{d-1},\qquad n\ge n_d .
\]
Since \(d-1<-1\),
$
\sum_{n=n_d}^\infty n^{d-1}<\infty.
$
Hence $\Lambda(d)<\infty.
$
\end{proof}
\begin{corollary}\label{20261113cor1}
For each fixed real $d$, there exists a positive constant $C_d$ and $N_d\in\mathbb N_0$ such that
\begin{equation*}\label{20261031for526}
 |A_n(d)|\leq C_d(1+n)^{N_d},\qquad n\in\mathbb N.
\end{equation*}
\end{corollary}
\begin{lemma}For every $k\in\mathbb N_0,$ and $ c>0$,
\begin{equation}\label{eq:convolution}
\sum_{n+m=k}A_n(c/2)A_m(c/2)=A_k(c).
\end{equation}
\end{lemma}
\begin{proof}
The binomial expansion gives, for \(|z|<1\),
\begin{equation}\label{eq:binomial-series}
(1-z)^{-c/2}
=\sum_{n=0}^\infty A_n(c/2)z^n.
\end{equation}
The series is absolutely convergent in
the open unit disk. Its Cauchy product with itself is therefore
legitimate and gives
\begin{align}
(1-z)^{-c}
&=\left(\sum_{n=0}^\infty A_n(c/2)z^n\right)
  \left(\sum_{m=0}^\infty A_m(c/2)z^m\right)\notag\\
&=\sum_{k=0}^\infty
\left(\sum_{n+m=k}A_n(c/2)A_m(c/2)\right)z^k.
\label{eq:cauchy-product}
\end{align}
On the other hand,
\[
(1-z)^{-c}=\sum_{k=0}^\infty A_k(c)z^k.
\]
Uniqueness of the coefficients of a power series now implies \eqref{eq:convolution}.
\end{proof}
The following operator inequality is well-known, for the reader's convenience we give a short proof here.
\begin{lemma}
Let $B_0$ be a trace-class operator on a Hilbert space with singular values $
s_1(B_0)\ge s_2(B_0)\ge\cdots\ge0.
$
 Then for every $r\in(1,\infty)$,
\begin{equation}\label{eq:interpolation}
 \|B_0\|_{\mathcal S_r}
\le\|B_0\|_{\mathrm{op}}^{\,1-1/r}\|B_0\|_{\mathcal S_1}^{\,1/r},
\end{equation}
and
\begin{equation}\label{20261030for2}
 \|B_0\|_{\mathrm{op}}=\sup_{i\in\mathbb N} s_i(B_0)\leq \|B_0\|_{\mathcal S_r},
\end{equation}
where $\|B_0\|_{\mathrm{op}}$ denotes the operator norm of $B_0$.
\end{lemma}
\begin{proof}
Since
\begin{equation*}
s_j(B_0)\le\|B_0\|_{\text{op}}=\sup_{i\in\mathbb N} s_i(B_0)\leq \|B_0\|_{\mathcal S_r}
=\left(\sum_i s_i(B_0)^r\right)^{\frac{1}{r}},\qquad j\in\mathbb N,\,r\in(1,\infty),
\end{equation*}
 we have
\begin{equation*}
\begin{aligned}
\|B_0\|_{\mathcal S_r}^r
=\sum_j s_j(B_0)^r
=\sum_j s_j(B_0)^{r-1}s_j(B_0)
\le\|B_0\|_{\text{op}}^{r-1}\sum_j s_j(B_0)
=\|B_0\|_{\text{op}}^{r-1}\|B_0\|_{\mathcal S_1}.
\end{aligned}
\end{equation*}
Taking the $r$-th root gives \eqref{eq:interpolation}.
\end{proof}
The measure $P$ is invariant under independent coordinate rotations,
but need not be invariant under general unitary transformations.
This distinction is essential to the following spectral argument.
\begin{lemma}\label{20261113lem1}
For every $u>-\frac{1}{2}$, \(\{\delta^u \mathbf z^\eta\}_{\eta\in \mathbb{N}_0^{(\mathbb{N})}}\) is an orthogonal family in $L^2(B,\mu)$.
\end{lemma}
\begin{proof}
Fix \(\eta,\xi\in \mathbb{N}_0^{(\mathbb{N})}\) with \(\eta\ne\xi\). Choose an index \(j\)
such that \(\eta_j\ne\xi_j\). For \(\theta\in\mathbb R\), let
\[
D_{j,\theta}\mathbf z
\triangleq(z_1,\ldots,z_{j-1},e^{i\theta}z_j,z_{j+1},\ldots),\qquad \mathbf z=(z_i)_{i\in\mathbb N}\in B.
\]
This map preserves \(B\) and \(\delta\). It also preserves \(P\),
and hence \(\mu\), because each coordinate Gaussian is invariant
under rotations. Restricting \(P\) to the rotation-invariant set
\(B\) preserves this invariance.

Consider
\[
I_{\eta,\xi}^{(u)}
\triangleq\int_B\delta(\mathbf z)^{2u}\mathbf z^\eta\overline{\mathbf z^\xi}\dd\mu(\mathbf z).
\]
This integral is absolutely convergent, since Cauchy--Schwarz inequality gives
\[
\int_B\delta(\mathbf z)^{2u}|\mathbf z^\eta||\mathbf z^\xi|\dd\mu(\mathbf z)
\le\sqrt{M_\eta(2u)M_\xi(2u)}<\infty.
\]
Changing variables by \(D_{j,\theta}\), we obtain
\begin{align}
I_{\eta,\xi}^{(u)}
=\int_B\delta(D_{j,\theta}\mathbf z)^{2u}
(D_{j,\theta}\mathbf z)^\eta
\overline{(D_{j,\theta}\mathbf z)^\xi}\dd\mu(\mathbf z)
=e^{i(\eta_j-\xi_j)\theta}I_{\eta,\xi}^{(u)}.
\label{eq:rotation}
\end{align}
Choose \(\theta\) with
\(e^{i(\eta_j-\xi_j)\theta}\ne1\). Then
\eqref{eq:rotation} forces \(I_{\eta,\xi}^{(u)}=0\).
Thus
\(\{\delta^u \mathbf z^\eta\}_{\eta\in \mathbb{N}_0^{(\mathbb{N})}}\) is an orthogonal family in $L^2(B,\mu)$.
\end{proof}

By Lemma \ref{20261113lem1}, the monomials \(\{\mathbf z^\eta\}_{\eta\in\mathbb{N}_0^{(\mathbb{N})}}\) form an orthogonal family in \(L^2(B,\mu)\). We define \(L_a^2(B,\mu)\) to be their closed linear span in \(L^2(B,\mu)\); that is,
\begin{equation}\label{20261114def1}
L_a^2(B,\mu)
\triangleq
\overline{\operatorname{span}\left\{
\mathbf z^\eta:\eta\in\mathbb{N}_0^{(\mathbb{N})}
\right\}}^{\,L^2(B,\mu)}.
\end{equation}
This space can be viewed as an infinite-dimensional analogue of the classical Bergman space on the unit ball of \(\mathbb C^n\); see \cite[3.1.1, p.~36]{Rud}.

\section{Proof of Theorem~\ref{thm:standard}}\label{20261105sbusection3}




We prove bounds for \(S_{a,b,c}\), which imply the same bounds for
\(T_{a,b,c}\) because
\[
 |\TT f(\mathbf z)|\le \SSS|f|(\mathbf z),\qquad \mathbf z\in B.
\]
For \(1<p\le\infty\), let \(p'\) be the conjugate exponent, with
\(p'=1\) when \(p=\infty\). In formulas involving \(\alpha/p\), put
\(\alpha/\infty\triangleq 0\).
For later use define
\begin{equation*}\label{eq:input}
 \mathcal I_p(b,\alpha)\triangleq \left\{(\alpha,b,p)\in \mathbb R\times\mathbb R\times [1,\infty]:
 \quad
 \begin{cases}
 b\ge\alpha,&p=1,\\
 p(b+1)>\alpha+1,&1<p<\infty,\\
 b>-1,&p=\infty.
 \end{cases}
\right\},\end{equation*}
\begin{equation*}\label{eq:output}
 \mathcal O_q(a,\beta)\triangleq \left\{(\beta,a,q)\in \mathbb R\times\mathbb R\times [1,\infty]:
 \quad
 \begin{cases}
 aq+\beta>-1,&q<\infty,\\
 a\ge0,&q=\infty.
 \end{cases}
\right\},
\end{equation*}
and
\[
 R_q(a,\beta)\triangleq
 \begin{cases}
 I_{aq+\beta}^{1/q},&q<\infty,\\
 1,&q=\infty\ \hbox{and }a\ge0,
 \end{cases}
\]
which is finite under \(\mathcal O_q(a,\beta)\).

\begin{proposition}\label{prop:suffp}
Suppose \(1<p\le\infty\), \(1\le q\le\infty\), and
\(\mathcal I_p(b,\alpha)\), \(\mathcal O_q(a,\beta)\) hold.
Then, for every \(c\in\R\), both operators are bounded, with
\begin{equation}\label{eq:pnorm}
 \norm{\TT}_{\Lt p\alpha\to\Lt q\beta},
 \ \norm{\SSS}_{\Lt p\alpha\to\Lt q\beta}
 \le
 R_q(a,\beta)\,
 W_{(b-\alpha/p)p',\,cp'}^{1/p'}.
\end{equation}
\end{proposition}

\begin{proof}
For \(1<p<\infty\), write the integrand using
\(f(\mathbf y)\delta(\mathbf y)^{\alpha/p}\).
H\"older's inequality with respect to \(\mu\) gives
\begin{align*}
 \SSS|f|(\mathbf z)
 \le \delta(\mathbf z)^a\norm f_{p,\alpha}
   \left(\int_B\delta(\mathbf y)^{(b-\alpha/p)p'}
                  d(\mathbf z,\mathbf y)^{-cp'}\dd\mu(\mathbf y)\right)^{1/p'}
 \le\delta(\mathbf z)^a\norm f_{p,\alpha}
       W_{(b-\alpha/p)p',\,cp'}^{1/p'}.
\end{align*}
Note that
\[
 (b-\alpha/p)p'>-1
 \quad\Longleftrightarrow\quad
 p(b+1)>\alpha+1.
\]
For \(p=\infty\), the same estimate is simply
\[
 \SSS|f|(\mathbf z)\le\delta(\mathbf z)^a\norm f_\infty W_{b,c},
\]
and its integrability condition is \(b>-1\).
Theorem~\ref{thm:FR} gives finiteness of the relevant moment.
Taking the target norm proves \eqref{eq:pnorm}.
The argument also proves absolute convergence at each \(\mathbf z\in B\).
\end{proof}

\begin{proposition}\label{prop:suff1}
Suppose \(1\le q<\infty\), \(b\ge\alpha\), and \(aq+\beta>-1\).
For every \(c\in\R\), both operators are bounded, and
\begin{equation*}\label{eq:1norm}
 \norm{\TT}_{\Lt1\alpha\to\Lt q\beta},
 \ \norm{\SSS}_{\Lt1\alpha\to\Lt q\beta}
 \le W_{aq+\beta,\,cq}^{1/q}.
\end{equation*}
\end{proposition}

\begin{proof}
Relative to \(\nu_\alpha\), the positive kernel is
\[
 \delta(\mathbf z)^a\delta(\mathbf y)^{b-\alpha}d(\mathbf z,\mathbf y)^{-c}.
\]
Minkowski's integral inequality gives
\begin{align*}
 \norm{\SSS|f|}_{q,\beta}
 \le\int_B |f(\mathbf y)|\delta(\mathbf y)^{b-\alpha}
       \left(\int_B\delta(\mathbf z)^{aq+\beta}
                        d(\mathbf z,\mathbf y)^{-cq}\dd\mu(\mathbf z)\right)^{1/q}
                   \dd\nu_\alpha(\mathbf y)
 \le W_{aq+\beta,\,cq}^{1/q}\norm f_{1,\alpha}.
\end{align*}
Here \(\delta(\mathbf y)^{b-\alpha}\le1\).
\end{proof}



Fix \(c\in\R\). Continuity of \(w\mapsto(1-w)^{-c}\) at \(0\)
gives a number \(0<\rho<1\) such that
\begin{equation}\label{eq:realpart}
 \Ree(1-w)^{-c}\ge\tfrac12,
 \qquad |w|<\rho.
\end{equation}

\begin{proposition} \label{prop:necessaryout}
If either \(T_{a,b,c}\) or \(S_{a,b,c}\) satisfies a boundedness
inequality from \(\Lt p\alpha\) to \(\Lt q\beta\), then
\(\mathcal O_q(a,\beta)\) holds.
\end{proposition}

\begin{proof}
Take \(0<r<\rho\) and
\[
 f_0(\mathbf y)\triangleq\delta(\mathbf y)^{-b}\ind_{B_r}(\mathbf y),\qquad \mathbf y\in B.
\]
This function belongs to \(\mathcal D\), and therefore to every $\Lt p\alpha$, since all weights are bounded above and below
on \(B_r\).
For every \(\mathbf z\in B\) and \(\mathbf y\in B_r\), \(|\langle \mathbf z,\mathbf y\rangle|<\rho\).
Equation~\eqref{eq:realpart} consequently gives
\[
 |\TT f_0(\mathbf z)|
 \ge\Ree\TT f_0(\mathbf z)
 \ge\tfrac12\mu(B_r)\delta(\mathbf z)^a.
\]
For \(S\), its kernel modulus is at least the real part in
\eqref{eq:realpart}, so the same lower bound holds.
Thus \(\delta^a\) must belong to \(\Lt q\beta\).
For finite \(q\), this is exactly \(aq+\beta>-1\);
for \(q=\infty\), it is exactly \(a\ge0\), by
Corollary~\ref{cor:radial}.
\end{proof}

\begin{proposition}\label{prop:necessaryin}
If either operator satisfies a boundedness inequality from
\(\Lt p\alpha\) to \(\Lt q\beta\), then
\(\mathcal I_p(b,\alpha)\) holds.
\end{proposition}

\begin{proof}
Let \(E=B_r\), where \(0<r<\rho\).
For every nonnegative \(f\in\mathcal D\), \eqref{eq:realpart} gives
\[
 |\TT f(\mathbf z)|\ge
 \tfrac12\delta(\mathbf z)^a\int_B\delta(\mathbf y)^b f(\mathbf y)\dd\mu(\mathbf y),
 \qquad \mathbf z\in E.
\]
The same estimate holds for \(S\).
The norm of \(\delta^a\ind_E\) in \(\Lt q\beta\) is finite and
strictly positive, including when \(q=\infty\).
The asserted operator inequality therefore forces
\begin{equation}\label{eq:functional}
 \int_Bf(\mathbf y)\delta(\mathbf y)^{b-\alpha}\dd\nu_\alpha(\mathbf y)
 \le C\norm f_{p,\alpha},
 \qquad 0\le f\in\mathcal D,
\end{equation}
for some finite \(C\).

For \(p=1\), \eqref{eq:functional} requires
\(\delta^{b-\alpha}\in L^\infty(\nu_\alpha)\); by equivalence of
measures and Corollary~\ref{cor:radial}, this means \(b\ge\alpha\).
For \(1<p<\infty\), the usual integral-functional criterion gives
\[
 \int_B\bigl(\delta(\mathbf y)^{b-\alpha}\bigr)^{p'}
                    \dd\nu_\alpha(\mathbf y)
 = I_{(b-\alpha)p'+\alpha}
 = I_{(b-\alpha/p)p'}<\infty.
\]
By \eqref{eq:radialthreshold}, this is
\(p(b+1)>\alpha+1\).
For \(p=\infty\), using increasing inner balls in
\eqref{eq:functional} gives
\[
 \int_B\delta(\mathbf y)^b\dd\mu(\mathbf y)=I_b<\infty,
\]
so \(b>-1\).

\end{proof}



\begin{theorem} \label{thm:endpoint}
For all real \(a,b,c,\alpha,\beta\), the following are equivalent:
\begin{itemize}
\item[$\mathrm{1)}$]
$\TT:\Lt1\alpha\longrightarrow\Lt\infty\beta
 \text{ is bounded};
$
\item[$\mathrm{2)}$]
$\SSS:\Lt1\alpha\longrightarrow\Lt\infty\beta
 \text{ is bounded};
$
\item[$\mathrm{3)}$]
$a\ge0,\quad b-\alpha\ge0,\quad\text{and}\quad c\le a+b-\alpha.
$
\end{itemize}
Moreover the two operator norms, interpreted as extended values,
are equal to
\begin{equation}\label{eq:endpointnorm}
 \sup_{\mathbf z,\mathbf y\in B}
 \frac{\delta(\mathbf z)^a\delta(\mathbf y)^{b-\alpha}}{|1-\langle \mathbf z,\mathbf y\rangle|^c}.
\end{equation}
Under conclusion $\mathrm{3)}$, this supremum is at most \(2^{|c|}\).
\end{theorem}

\begin{proof}
Put \(v=b-\alpha\) and consider the continuous kernel relative to
\(\nu_\alpha\),
\[
 H(\mathbf z,\mathbf y)\triangleq\delta(\mathbf z)^a\delta(\mathbf y)^v(1-\langle \mathbf z,\mathbf y\rangle)^{-c},
 \qquad \mathbf z,\mathbf y\in B.
\]
A finite supremum of \(|H|\) immediately bounds the
\(L^1\to L^\infty\) norm from above.

For the reverse inequality, fix \(\mathbf z_0,\mathbf y_0\in B\),
let \(h_0\triangleq|H(\mathbf z_0,\mathbf y_0)|\), and choose \(\theta\in\R\) so that
\(e^{-i\theta}H(\mathbf z_0,\mathbf y_0)=h_0\).
Given \(0<\eta<h_0\), choose relatively norm-open neighborhoods
\(U,V\) of these points, contained in some inner ball of $B$, with
\[
 \Ree\bigl(e^{-i\theta}H(\mathbf z,\mathbf y)\bigr)>h_0-\eta,
 \qquad \mathbf z\in U,\ \mathbf y\in V.
\]
The weighted measures of $U$ and $V$ are positive and finite by
Proposition~\ref{prop:support2}.
The function
\[
 f(\mathbf y)\triangleq e^{-i\theta}\frac{\ind_V(\mathbf y)}{\nu_\alpha(V)},\qquad \mathbf y\in B,
\]
belongs to \(\mathcal D\) and has \(\norm f_{1,\alpha}=1\).
For every \(\mathbf z\in U\), \(|\TT f(\mathbf z)|>h_0-\eta\).
Since \(\nu_\beta(U)>0\), $\|\TT f\|_{\Lt\infty\beta}>h_0-\eta$. Taking the supremum over \(\mathbf z_0,\mathbf y_0\) and
letting \(\eta\downarrow0\) proves \eqref{eq:endpointnorm}.
The same proof applies to the positive kernel \(|H|\).

It remains to determine when \eqref{eq:endpointnorm} is finite.
At \(\mathbf y=\mathbf 0\) the modulus is \(\delta(\mathbf z)^a\), forcing \(a\ge0\).
At \(\mathbf z=\mathbf 0\) it is \(\delta(\mathbf y)^v\), forcing \(v\ge0\).
For \(\mathbf z=\mathbf y=r \mathbf e_1\), \(0<r<1\), it equals
\[
 (1-r^2)^{a+v-c}.
\]
Finiteness as \(r\uparrow1\) forces \(c\le a+v\).

Conversely, assume \(a,v\ge0\) and \(c\le a+v\).
If \(c\le0\), then \(\delta(\mathbf z)^a,\delta(\mathbf y)^v\le1\) and
\(d(\mathbf z,\mathbf y)\le2\), so \(|H(\mathbf z,\mathbf y)|\le2^{-c}\).
If \(c>0\), put \(u=\delta(\mathbf z)\), \(w=\delta(\mathbf y)\), and
\(\tau=\max(u,w)\). We have
\begin{align*}
 d(\mathbf z,\mathbf y)&\ge1-\norm {\mathbf z}\norm {\mathbf y}
 \ge\frac{(1-\norm {\mathbf z}^2)+(1-\norm{\mathbf  y}^2)}{2}
 =\frac{u+w}{2}\ge\frac{\tau}{2}.
\end{align*}
Since \(u^aw^v\le\tau^{a+v}\) and \(0<\tau\le1\),
$
 |H(\mathbf z,\mathbf y)|\le2^c\tau^{a+v-c}\le2^c.
$
This proves sufficiency and the stated bound.
\end{proof}

\begin{proof}[Proof of Theorem~\ref{thm:standard}]
Propositions~\ref{prop:suffp} and~\ref{prop:suff1} prove sufficiency
for every pair except \(p=1,q=\infty\).
Propositions~\ref{prop:necessaryout} and~\ref{prop:necessaryin}
give the necessary conditions directly for each operator.
Theorem~\ref{thm:endpoint} settles the remaining endpoint, including
its equality case. These statements exhaust Table~\ref{tab:standard}.
\end{proof}

\section{Proof of Theorem~\ref{thm:schatten2}}\label{20261105sbusection4}
The first boundedness assertion in
Theorem~\ref{thm:schatten2} is the \(p=q=2\) case of
Theorem~\ref{thm:standard}. We assume \eqref{eq:hilbertadmissible}
throughout the remaining proof of Theorem~\ref{thm:schatten2}, so \(t>-1\).

The operators
\begin{equation*}
A_T\triangleq U_\beta T_{a,b,c}U_\alpha^{-1}\quad\text{and}\quad
 A_S\triangleq U_\beta S_{a,b,c}U_\alpha^{-1}
\end{equation*}
 have kernels
\begin{equation}\label{eq:unitarykernels}
 K_T(\mathbf z,\mathbf y)\triangleq\frac{\delta(\mathbf z)^u\delta(\mathbf y)^v}{(1-\langle \mathbf z,\mathbf y\rangle)^c},
 \qquad
 K_S(\mathbf z,\mathbf y)\triangleq\frac{\delta(\mathbf z)^u\delta(\mathbf y)^v}{|1-\langle \mathbf z,\mathbf y\rangle|^c},
 \qquad \mathbf z,\mathbf y\in B,
\end{equation}
relative to \(\mu\). By Corollary \ref{cor:radial} and Theorem \ref{thm:FR}, $K_T,K_S\in L^2(B\times B,\mu\times \mu)$.  By Lemma \ref{20261113lem2},
\begin{equation}\label{eq:HSformula23}
 \|T_{a,b,c}\|_{\mathcal S_2}^2
 =\|S_{a,b,c}\|_{\mathcal S_2}^2= \|A_T\|_{\mathcal S_2}^2
 =\|A_S\|_{\mathcal S_2}^2
 =\int_B\!\int_B
 \frac{\delta(\mathbf z)^{2u}\delta(\mathbf y)^{2v}}{|1-\ip {\mathbf z}{\mathbf y}|^{2c}}
 \dd\mu(\mathbf y)\dd\mu(\mathbf z)<\infty.
\end{equation}

The binomial formulas
give
\begin{align}
 A_T&=\sum_{n\ge0}A_n(c)R_{n,0},\label{eq:Tseries}\\
 A_S&=\sum_{n,m\ge0}A_n(c/2)A_m(c/2)R_{n,m},
 \label{eq:Sseries}
\end{align}
where $A_n(c)$ is defined at \eqref{20261112for2}.
The identities in \eqref{eq:Tseries}--\eqref{eq:Sseries} hold in
the Hilbert--Schmidt norm. Indeed, Corollary~\ref{20261113cor1}
gives polynomial growth of the coefficients $A_n(c)$ and
$A_n(c/2)$. Hence, for some $M\geq0$ and $C>0$,
\[
\sum_{n\geq0}|A_n(c)|\,|\langle\mathbf z,\mathbf y\rangle|^n
\leq C(1-Q(\mathbf z,\mathbf y))^{-M},
\]
and the analogous product estimate holds for the double series.
Lemma~\ref{lem:moduluscap} implies
\[
\int_B\!\int_B
\delta(\mathbf z)^{2u}\delta(\mathbf y)^{2v}
(1-Q(\mathbf z,\mathbf y))^{-2M}\,
\mathrm d\mu(\mathbf y)\mathrm d\mu(\mathbf z)<\infty.
\]
Thus the kernel expansions in \eqref{eq:Tseries} and \eqref{eq:Sseries} converge in $L^2(B\times B,\mu\times\mu)$
by dominated convergence. By Lemma \ref{20261113lem2},
\begin{equation}\label{20261115for1}
\sum_{n=0}^NA_n(c)R_{n,0}\rightarrow A_T,\qquad
\sum_{n,m=0}^N A_n(c/2)A_m(c/2)R_{n,m}\rightarrow A_S,\qquad \text{in}\,\,\mathcal S_2,\quad \text{as}\,N\to\infty.
\end{equation}





\begin{proposition} \label{thm:svd}
Under \eqref{eq:hilbertadmissible}, the nonzero singular values of
\(T_{a,b,c}:\mathcal H_\alpha\to\mathcal H_\beta\), with
multiplicity, are the nonzero members of
\begin{equation}\label{eq:singularvalues}
 \sigma_\eta=
 \frac{|(c)_{|\eta|}|}{\eta!}
 \sqrt{M_\eta(2u)M_\eta(2v)},\qquad\eta\in \mathbb{N}_0^{(\mathbb{N})}.
\end{equation}
In particular, for any \(r>0\),
\[
 \|T_{a,b,c}\|_{\mathcal S_r}^r
 =\sum_{\eta\in \mathbb{N}_0^{(\mathbb{N})}}
 \left(\frac{|(c)_{|\eta|}|}{\eta!}\right)^r
 \bigl(M_\eta(2u)M_\eta(2v)\bigr)^{r/2},
\]
where both sides may be infinite.
\end{proposition}

\begin{proof}
By Lemma \ref{20261113lem1}, the families
\(\{\delta^u \mathbf z^\eta\}_{\eta\in \mathbb{N}_0^{(\mathbb{N})}}\) and
\(\{\delta^v \mathbf z^\eta\}_{\eta\in \mathbb{N}_0^{(\mathbb{N})}}\) are orthogonal families in $L^2(B,\mu)$. Equations
\eqref{eq:infinite-rank-one-expansion} and \eqref{eq:Tseries} give
\[
 A_T
 =\sum_{\eta\in \mathbb{N}_0^{(\mathbb{N})}}\frac{(c)_{|\eta|}}{\eta!}\sqrt{M_\eta(2u)M_\eta(2v)}
 \left(\frac{\delta^u \mathbf z^\eta}{\sqrt{M_\eta(2u)}}\right)\otimes\left(\frac{\delta^v \mathbf z^\eta}{\sqrt{M_\eta(2v)}}\right).
\]
This is a singular-value decomposition which proves the assertion.
\end{proof}

For \(A_S\), the mixed monomials in \eqref{eq:Sseries} need not
be orthogonal within a fixed character. Its coefficients therefore
cannot simply be read as eigenvalues or singular values. We use
operator-ideal estimates and character blocks instead.
\begin{proposition}\label{20261030prop1}
Under \eqref{eq:hilbertadmissible},
\begin{equation}\label{eq:all-membership}
T_{a,b,c},S_{a,b,c}
\in\bigcap_{1<r<\infty}\mathcal S_r(H_\alpha,H_\beta).
\end{equation}
\end{proposition}
\begin{proof}

Firstly, \eqref{eq:piecerapid} and \eqref{20261030for2} imply that
\begin{equation}\label{eq:rapid-op}
\|R_{n,m}\|_{\text{op}}
\le C_D(1+n+m)^{-D},
\qquad\text{for every }D>0,\,n,m\in\mathbb N_0.
\end{equation}
By the definition of $J_{n+m}(2u)$ at \eqref{eq:Jdefinition} and \eqref{eq:piecetracebound}, we have
\begin{equation}\label{20261030for3}
\sup_{m,n\in\mathbb N_0}\|R_{n,m}\|_{\mathcal S_1}\leq
\sup_{m,n\in\mathbb N_0}\sqrt{J_{n+m}(2u)J_{n+m}(2v)}\leq C_0\triangleq \sqrt{I_{2u}I_{2v}}<\infty.
\end{equation}
Fix $r>1$. Applying~\eqref{eq:interpolation} to $R_{n,m}$ and
using~\eqref{eq:rapid-op}--\eqref{20261030for3}, we obtain
\begin{equation}\label{20261030for5}
\begin{aligned}
\|R_{n,m}\|_{\mathcal S_r}
\le\|R_{n,m}\|_{\text{op}}^{\,1-1/r}
\|R_{n,m}\|_{\mathcal S_1}^{\,1/r}
\le C_D^{\,1-1/r}C_0^{\,1/r}
(1+n+m)^{-D(1-1/r)}.
\end{aligned}
\end{equation}
Now let $L$ be any positive number. Choose
$
D\triangleq\frac{Lr}{r-1}.
$
Then
$
D\left(1-\frac1r\right)=L.
$
Equation~\eqref{20261030for5} therefore becomes
\begin{equation}\label{eq:rapid-Sr}
\|R_{n,m}\|_{\mathcal S_r}
\le C_{r,L}(1+n+m)^{-L},
\qquad n,m\ge0 ,
\end{equation}
where
$
C_{r,L}\triangleq C_0^{1/r}C_{Lr/(r-1)}^{\,1-1/r}.
$
The constants may depend on $r,L,u,v$, and the Gaussian measure.
They do not depend on $n,m$.

By Corollary \ref{20261113cor1}, there exist a positive constant $C_c$ and nonnegative integer $N_c$ such that
\begin{equation}\label{eq:polynomial-coefficients}
\begin{aligned}
|A_n(c)|&\le C_c(1+n)^{N_c},\qquad
|A_n(c/2)|&\le C_c(1+n)^{N_c},\qquad n\in\mathbb N.
\end{aligned}
\end{equation}

Fix $r>1$. From~\eqref{eq:rapid-Sr}
and~\eqref{eq:polynomial-coefficients},
\begin{equation*}
\begin{aligned}
\|A_n(c)R_{n,0}\|_{\mathcal S_r}
=|A_n(c)|\,\|R_{n,0}\|_{\mathcal S_r}
\le C_cC_{r,L}(1+n)^{N_c-L}.
\end{aligned}
\end{equation*}
Choose
$
L>N_c+1.
$
Then $N_c-L<-1$, so
$
\sum_{n=0}^\infty(1+n)^{N_c-L}<\infty.
$
It follows that
\begin{equation}\label{eq:T-absolute}
\sum_{n=0}^\infty
\|A_n(c)R_{n,0}\|_{\mathcal S_r}<\infty.
\end{equation}
For $r\ge1$, $\mathcal S_r$, equipped with its Schatten norm, is a Banach
space. Thus absolute convergence
in~\eqref{eq:T-absolute} gives an $\mathcal S_r$-limit of the operator series.

Let
\[
P_N^T\triangleq \sum_{n=0}^N A_n(c)R_{n,0},\qquad N\in\mathbb N,
\]
and denote its $\mathcal S_r$-limit temporarily by $B_T$. Thus
\begin{equation}\label{20261030for6}
\|P_N^T-B_T\|_{\mathcal S_r}\longrightarrow0.
\end{equation}
We must verify that this limit is the original operator $A_T$.

By \eqref{20261030for2} and \eqref{20261115for1},
\[
\begin{aligned}
\|B_T-A_T\|_{ \text{op}}
\le\|B_T-P_N^T\|_{\text{op}}
+\|P_N^T-A_T\|_{\text{op}}
\le\|B_T-P_N^T\|_{\mathcal S_r}
+\|P_N^T-A_T\|_{ \mathcal S_2}
\longrightarrow0.
\end{aligned}
\]
Hence $B_T=A_T$. We have proved
$
A_T\in \mathcal S_r.
$

For the double series, \eqref{eq:rapid-Sr}
and~\eqref{eq:polynomial-coefficients} give
\begin{equation*}
\begin{aligned}
\|A_n(c/2)A_m(c/2)R_{n,m}\|_{\mathcal S_r}
&=|A_n(c/2)|\,|A_m(c/2)|\,\|R_{n,m}\|_{\mathcal S_r}\\
&\le C_c^2C_{r,L}
(1+n)^{N_c}(1+m)^{N_c}(1+n+m)^{-L}\\
&\le C_c^2C_{r,L}(1+n+m)^{2N_c-L}.
\end{aligned}
\end{equation*}
We now sum over all $n,m\ge0$. For each $k\ge0$, there are exactly
$k+1$ ordered pairs satisfying
\[
n+m=k:
\qquad(0,k),(1,k-1),\ldots,(k,0).
\]
Therefore,
\begin{equation*}
\begin{aligned}
\sum_{n,m\ge0}
\|A_n(c/2)A_m(c/2)R_{n,m}\|_{\mathcal S_r}
&\le C\sum_{k=0}^\infty
\sum_{n+m=k}(1+k)^{2N_c-L}\\
&=C\sum_{k=0}^\infty
(k+1)(1+k)^{2N_c-L}\\
&=C\sum_{k=0}^\infty
(1+k)^{2N_c+1-L}.
\end{aligned}
\end{equation*}
Choose
$
L>2N_c+2.
$
Then the final series converges, and hence
\begin{equation*}
\sum_{n,m\ge0}
\|A_n(c/2)A_m(c/2)R_{n,m}\|_{\mathcal S_r}<\infty.
\end{equation*}
By similar arguments for $A_T$, we have
$
A_S\in \mathcal S_r.
$
Unitary invariance then gives \eqref{eq:all-membership}.
\end{proof}
Proposition \ref{20261030prop1} proves conclusion 1 of Theorem~\ref{thm:schatten2}.



\begin{proposition} \label{20261031prop1}
Under \eqref{eq:hilbertadmissible},
assume
\begin{equation}\label{eq:assumptions}
c<t+1.
\end{equation}
Then
\begin{equation}\label{eq:conclusion}
 T_{a,b,c},S_{a,b,c}\in \mathcal S_1 (\mathcal{H}_{\alpha},\mathcal{H}_{\beta}).
\end{equation}
Moreover, the right-side series of  \eqref{eq:Tseries} and \eqref{eq:Sseries}
converge absolutely in trace norm.
\end{proposition}
\begin{proof}
For \(c>0\), every factor in \((c)_n= c(c+1)\cdots(c+n-1)\) is positive, and hence
\[
A_n(c)=\frac{ c(c+1)\cdots(c+n-1)}{n!}>0,\qquad n\ge0 .
\]
Then \eqref{eq:negative-asymptotic235} gives
\begin{equation}\label{eq:positive-coeff-asymptotic}
A_n(c)
=\frac{\Gamma(n+c)}{\Gamma(c)\Gamma(n+1)}
\sim\frac{n^{c-1}}{\Gamma(c)},
\qquad\text{as}\,\,\,n\to\infty.
\end{equation}
By  \eqref{eq:Jasymptotic}, for every fixed \(s>-1\), we have
\begin{equation*}\label{eq:J-asymptotic}
J_n(s)\sim \widetilde{h}(1)\Gamma(s+1)n^{-s-1},
\qquad\text{as}\,\,\,n\to\infty.
\end{equation*}
Applying this with \(s=2u\) and \(s=2v\) yields
\begin{align}
W_n
&\triangleq\sqrt{J_n(2u)J_n(2v)}\notag\\
&\sim
\widetilde{h}(1)\sqrt{\Gamma(2u+1)\Gamma(2v+1)}
n^{-u-v-1}\notag\\
&=
\widetilde{h}(1)\sqrt{\Gamma(2u+1)\Gamma(2v+1)}
n^{-t-1}.
\label{eq:W-asymptotic}
\end{align}
Here $t+1>0$, because \(2u+1>0\)
and \(2v+1>0\).

Define
\begin{equation*}\label{eq:bn-definition}
b_n\triangleq A_n(c)W_n,\qquad
C_{u,v,c}
\triangleq\frac{\widetilde{h}(1)\sqrt{\Gamma(2u+1)\Gamma(2v+1)}}{\Gamma(c)}>0.
\end{equation*}
Multiplying \eqref{eq:positive-coeff-asymptotic}
and \eqref{eq:W-asymptotic}, we obtain the more precise form:
\begin{equation}\label{eq:bn-asymptotic}
b_n\sim C_{u,v,c}\,n^{c-t-2}
\qquad\text{as}\,\,\,n\to\infty.
\end{equation}
In particular,
\begin{equation*}\label{eq:bn-comparability}
b_n\asymp(1+n)^{c-t-2},\qquad n\ge0 .
\end{equation*}
Here \(\asymp\) means that there exist constants \(0<C_1\le C_2<\infty\),
independent of \(n\), for which
\[
C_1(1+n)^{c-t-2}\le b_n\le C_2(1+n)^{c-t-2}.
\]
Indeed, \eqref{eq:bn-asymptotic} gives such bounds for all sufficiently
large \(n\). Each of the remaining finitely many terms, including
\(b_0=C_0\), is finite and strictly positive, so the constants can be
adjusted to cover all \(n\ge0\).

Let
$
\tau\triangleq t+1-c>0.
$
Then \(c-t-2=-1-\tau\), so
\[
\sum_{n=0}^\infty(1+n)^{c-t-2}
=\sum_{n=0}^\infty(1+n)^{-1-\tau}<\infty.
\]
It follows that
\begin{equation*}\label{eq:Bc-definition}
\mathcal B_c\triangleq \sum_{n=0}^\infty b_n
=\sum_{n=0}^\infty A_n(c)W_n<\infty.
\end{equation*}

Using \eqref{eq:piecetracebound} and positivity of \(A_n(c)\), we conclude
that the right-side series of  \eqref{eq:Tseries}
converge absolutely in trace norm:
\begin{equation}\label{eq:T-positive-summability}
\sum_{n=0}^\infty
\|A_n(c)R_{n,0}\|_{\mathcal S_1}
\le\sum_{n=0}^\infty A_n(c)W_n
=\mathcal B_c<\infty.
\end{equation}

Because \(c/2>0\), all coefficients \(A_n(c/2)\) are positive. By
\eqref{eq:piecetracebound},
\[
\|A_n(c/2)A_m(c/2)R_{n,m}\|_{\mathcal S_1}
\le A_n(c/2)A_m(c/2)W_{n+m}.
\]
By \eqref{eq:convolution}, we obtain
\begin{align}
\sum_{n,m\ge0}
\|A_n(c/2)A_m(c/2)R_{n,m}\|_{\mathcal S_1}\notag
&\le
\sum_{n,m\ge0}A_n(c/2)A_m(c/2)W_{n+m}\notag\\
&=
\sum_{k=0}^\infty W_k
\sum_{n+m=k}A_n(c/2)A_m(c/2)\notag\\
&=
\sum_{k=0}^\infty A_k(c)W_k
=\mathcal B_c<\infty.
\label{eq:S-positive-summability}
\end{align}

For \(c<0\), then both \(c\) and \(c/2\) are negative. Applying
the final conclusion of Lemma \ref{20261112lem1},
we have
\[
\Lambda(c)<\infty,\qquad\Lambda(c/2)<\infty.
\]
Using the uniform bound in \eqref{20261030for3}, we obtain
\begin{equation}\label{eq:T-negative-summability}
\sum_{n=0}^\infty\|A_n(c)R_{n,0}\|_{\mathcal S_1}
\le C_0\sum_{n=0}^\infty|A_n(c)|
=C_0\Lambda(c)<\infty.
\end{equation}
Similarly, Tonelli's theorem gives
\begin{align}
\sum_{n,m\ge0}
\|A_n(c/2)A_m(c/2)R_{n,m}\|_{\mathcal S_1}
&\le
C_0\sum_{n,m\ge0}|A_n(c/2)|\,|A_m(c/2)|\notag\\
&=
C_0\left(\sum_{n=0}^\infty|A_n(c/2)|\right)^2
=C_0\Lambda(c/2)^2<\infty.
\label{eq:S-negative-summability}
\end{align}
Both operator series are therefore absolutely summable in trace norm
for every \(c<0\).

For \(c=0\),
\[
A_0(0)=1,\qquad A_n(0)=0,\quad n\ge1.
\]
In this case, we have
\begin{equation}\label{eq:rank-one}
A_T=A_S=R_{0,0}=\delta^u\otimes\delta^v.
\end{equation}
Indeed,
\[
(R_{0,0}h)(\mathbf z)
=\delta(\mathbf z)^u\int_B\delta(\mathbf y)^v h(\mathbf y)\dd\mu(\mathbf y)
=\langle h,\delta^v\rangle_{L^2(B,\mu)}\delta(\mathbf z)^u,\qquad \mathbf z\in B,\,\,h\in L^2(B,\mu).
\]
The vectors \(\delta^u,\delta^v\) belong to \(L^2(B,\mu)\) by \eqref{eq:radialthreshold} and
\eqref{eq:assumptions}, and both are nonzero because
\(\delta>0\) on \(B\) and \(\mu(B)>0\).
Thus the operator has rank one, with its unique nonzero singular
value equal to
\begin{equation}\label{eq:rank-one-norm}
\|\delta^u\|_{L^2(B,\mu)}\|\delta^v\|_{L^2(B,\mu)}
=\sqrt{I_{2u}I_{2v}}=C_0.
\end{equation}
In particular, both operators are trace class in this case.

Notice also that \(u>-\tfrac12\) and \(v>-\tfrac12\) imply
\begin{equation}\label{eq:t-positive}
t=u+v>-1,\qquad t+1>0.
\end{equation}
Consequently every \(c\le0\) automatically satisfies \(c<t+1\).

For \(c>0\) with \(c<t+1\), absolute trace-norm summability was
proved in \eqref{eq:T-positive-summability}
and \eqref{eq:S-positive-summability}. For \(c<0\), it was proved in
\eqref{eq:T-negative-summability}
and \eqref{eq:S-negative-summability}. The case \(c=0\) was treated
directly in \eqref{eq:rank-one}.

Similar to the argument around \eqref{20261030for6}, define
\begin{equation}\label{eq:partial-sums}
P_N^T\triangleq \sum_{n=0}^N A_n(c)R_{n,0},
\qquad
P_N^S\triangleq\sum_{\substack{n,m\ge0\\n+m\le N}}
A_n(c/2)A_m(c/2)R_{n,m},\qquad N\in\mathbb N.
\end{equation}
Since \(\mathcal S_1 \) is a Banach space, there
exist \(B_T,B_S\in \mathcal S_1 \) such that
\begin{equation}\label{eq:S1-limits}
\|P_N^T-B_T\|_{\mathcal S_1}\longrightarrow0,\qquad
\|P_N^S-B_S\|_{\mathcal S_1}\longrightarrow0.
\end{equation}

For any trace-class operator \(C\), the elementary inequality for its
singular values gives
\[
\|C\|_{\mathcal S_2}
\le \|C\|_{\mathcal S_1}.
\]
Thus \eqref{eq:S1-limits} also gives convergence in \(\mathcal S_2\). By the proof of Proposition \ref{20261030prop1}, we have
$$
\|P_N^T-A_T\|_{\mathcal S_2}\longrightarrow0,\qquad
\|P_N^S-A_S\|_{\mathcal S_2}\longrightarrow0,
$$
which imply that
\[
\|B_T-A_T\|_{\mathcal S_2}
\le\|B_T-P_N^T\|_{\mathcal S_1}
+\|P_N^T-A_T\|_{\mathcal S_2}\longrightarrow0,
\]
and
\[
\|B_S-A_S\|_{\mathcal S_2}
\le\|B_S-P_N^S\|_{\mathcal S_1}
+\|P_N^S-A_S\|_{\mathcal S_2}\longrightarrow0.
\]
Hence
$
B_T=A_T,\, B_S=A_S\in \mathcal S_1.
$

Finally, the unitary equivalences \eqref{20261031for1} imply \eqref{eq:conclusion}.
\end{proof}
For every $\varepsilon\in(0,1)$, let \(\mathscr M_\varepsilon f\triangleq \ind_{\{\mathbf z\in B:\delta(\mathbf z)\ge\varepsilon\}}\cdot f\), for every \(f\in L^2(B,\mu)\).
 Define
\[
 J_n^\varepsilon(s)\triangleq
 \int_{\{\mathbf z\in B:\delta(\mathbf z)\geq\varepsilon\}}\delta(\mathbf z)^s\|\mathbf z\|^{2n}\,\mathrm d\mu(\mathbf z),\qquad \varepsilon\in(0,\infty),\,s\in\mathbb R,\,n\in\mathbb N_0.
\]
Since $\delta(\mathbf z)\geq\varepsilon$ implies $\|\mathbf z\|^2\leq1-\varepsilon$,
\[
 J_n^\varepsilon(s)\leq
 C_{\varepsilon,s}(1-\varepsilon)^n,\qquad
 C_{\varepsilon,s}\triangleq \int_{\{\mathbf z\in B:\delta(\mathbf z)\geq\varepsilon\}}\delta^s\,\mathrm d\mu<\infty.                                                    \]
This is finite for every real $s$, since $\delta^s$ is bounded on
$\{\mathbf z\in B:\delta(\mathbf z)\geq\varepsilon\}$.

By \eqref{eq:rank-one-identities}, apply the
triangle inequality to \eqref{eq:infinite-rank-one-expansion}, and then Cauchy--Schwarz inequality to the sum:
\begin{align*}
 \|\mathscr M_\varepsilon R_{n,m}\mathscr M_\varepsilon\|_{\mathcal S_1}
 &\leq
 \left(\sum_{|\eta|=n,|\zeta|=m}
 \frac{n!m!}{\eta!\zeta!}
 \|\mathscr M_\varepsilon(\delta^u \mathbf z^{\eta+\zeta})\|_{L^2(B,\mu)}^2\right)^{1/2}\\
 &\quad\times
 \left(\sum_{|\eta|=n,|\zeta|=m}
 \frac{n!m!}{\eta!\zeta!}
 \|\mathscr M_\varepsilon(\delta^v \mathbf z^{\eta+\zeta})\|_{L^2(B,\mu)}^2\right)^{1/2}.
\end{align*}
The multinomial
identity
\[
 \sum_{|\eta|=n}\frac{n!}{\eta!}|\mathbf z^\eta|^2=\|\mathbf z\|^{2n}
\]
(first for finite sums, then by monotone convergence) and Tonelli's
theorem show that the two factors squared are
$J_{n+m}^\varepsilon(2u)$ and $J_{n+m}^\varepsilon(2v)$.
Consequently
\begin{equation}\label{20261031for6}
 \|\mathscr M_\varepsilon R_{n,m}\mathscr M_\varepsilon\|_{\mathcal S_1}
 \leq C_\varepsilon(1-\varepsilon)^{n+m}.
\end{equation}
In particular, the compressed pieces are trace class. The geometric
power is $n+m$, because each moment contributes that power before
taking the square root of the product.
\begin{proposition}
\label{20261031prop2}
Under \eqref{eq:hilbertadmissible}, if either
\(T_{a,b,c}\) or \(S_{a,b,c}\) belongs to
\(\mathcal S_1(\mathcal H_\alpha,\mathcal H_\beta)\),
then \(c<t+1\).
\end{proposition}
\begin{proof}
Recall \eqref{eq:Tseries} and \eqref{eq:Sseries},
\[
 A_T=\sum_{n\geq0}A_n(c)R_{n,0},\qquad
 A_S=\sum_{n,m\geq0}A_n(c/2)A_m(c/2)R_{n,m}.
\]
By Corollary \ref{20261113cor1} and \eqref{20261031for6},
\[
 \sum_n |A_n(c)|\cdot \|\mathscr M_\varepsilon R_{n,0}\mathscr M_\varepsilon\|_{\mathcal S_1}
 \leq C_\varepsilon\sum_n|A_n(c)|(1-\varepsilon)^n\leq C_{\varepsilon,c}\sum_n(1+n)^{N_c}(1-\varepsilon)^n<\infty,
\]
and
\begin{equation}\label{20161031for8}
 \sum_{n,m}|A_n(c/2)A_m(c/2)|\cdot
 \|\mathscr M_\varepsilon R_{n,m}\mathscr M_\varepsilon\|_{\mathcal S_1}
 \leq C_\varepsilon
 \left(\sum_n|A_n(c/2)|(1-\varepsilon)^n\right)^2<\infty.
\end{equation}
Thus both compressed expansions converge absolutely in trace norm
for every real $c$.

These sums represent the compressed operators. Indeed, for all  $\mathbf z,\mathbf y\in\{\mathbf z\in B:\delta(\mathbf z)\geq\varepsilon\}$,
\[
 |\langle \mathbf z,\mathbf y\rangle|\leq\|\mathbf z\|\|\mathbf y\|\leq1-\varepsilon<1.
\]
Hence the binomial series for $(1-\langle \mathbf z,\mathbf y\rangle)^{-c}$ and
the product of the two binomial series for
$|1-\langle \mathbf z,\mathbf y\rangle|^{-c}$ converge absolutely and uniformly.
Their sums are precisely the kernels above. The trace-norm limits
are also operator-norm limits, so they equal
$\mathscr M_\varepsilon A_T \mathscr M_\varepsilon$ and $\mathscr M_\varepsilon A_S \mathscr M_\varepsilon$.

Taking traces in \eqref{eq:infinite-rank-one-expansion} is justified by \eqref{20261031for6}. Since
$\operatorname{Tr}(f\otimes g)=\langle f,g\rangle_{L^2(B,\mu)}$, for all $ f,g\in L^2(B,\mu)$, we obtain
\begin{align*}
 \operatorname{Tr}(\mathscr M_\varepsilon R_{n,m}\mathscr M_\varepsilon)
 &=\int_{\{\mathbf z\in B:\delta(\mathbf z)\geq\varepsilon\}}\delta(\mathbf z)^{u+v}
 \sum_{\eta,\zeta}
 \frac{n!m!}{\eta!\zeta!}|\mathbf z^\eta\overline{\mathbf z^\zeta}|^2\,\mathrm d\mu(\mathbf z)\\
 &=\int_{\{\mathbf z\in B:\delta(\mathbf z)\geq\varepsilon\}}
 \delta(\mathbf z)^t\|\mathbf z\|^{2(n+m)}\,\mathrm d\mu(\mathbf z).
\end{align*}
The second equality is the multinomial identity applied separately
to $\eta$ and $\zeta$.

The trace is continuous on $\mathcal S_1$, so the trace-norm convergence
in \eqref{20161031for8} permits termwise tracing. Thus
\[
 \operatorname{Tr}(\mathscr M_\varepsilon A_T \mathscr M_\varepsilon)
 =\int_{\{\mathbf z\in B:\delta(\mathbf z)\geq\varepsilon\}}\delta(\mathbf z)^t
       \sum_{n\geq0}A_n(c)\|\mathbf z\|^{2n}\,\mathrm d\mu(\mathbf z),
\]
and
\[
 \operatorname{Tr}(\mathscr M_\varepsilon A_S \mathscr M_\varepsilon)
 =\int_{\{\mathbf z\in B:\delta(\mathbf z)\geq\varepsilon\}}\delta(\mathbf z)^t
       \sum_{n,m\geq0}A_n(c/2)A_m(c/2)
                   \|\mathbf z\|^{2(n+m)}\,\mathrm d\mu(\mathbf z).
\]
The interchanges are also justified by uniform absolute convergence
on the truncated set $\{\mathbf z\in B:\delta(\mathbf z)\geq\varepsilon\}$. For $0\leq x<1$, the binomial identity gives
\[
 \sum_{n\geq0}A_n(c)x^n=(1-x)^{-c},
\]
and
\[
 \sum_{n,m\geq0}A_n(c/2)A_m(c/2)x^{n+m}
 =\left(\sum_{n\geq0}A_n(c/2)x^n\right)^2
 =(1-x)^{-c}.
\]
Taking $x=\|\mathbf z\|^2$ yields the key identity
\begin{equation}\label{20261031for10}
 \operatorname{Tr}(\mathscr M_\varepsilon A_T \mathscr M_\varepsilon)
 =\operatorname{Tr}(\mathscr M_\varepsilon A_S \mathscr M_\varepsilon)
 =\int_{\{\delta\geq\varepsilon\}}\delta^{t-c}\,\mathrm d\mu.
\end{equation}

Suppose $A$ is either operator on $L^2(B,\mu)$ and $A\in\mathcal S_1$. The ideal property
and $\|\mathscr M_\varepsilon\|_{\rm op}\leq1$ imply
\[
 \|\mathscr M_\varepsilon A \mathscr M_\varepsilon\|_{\mathcal S_1}
 \leq\|A\|_{\mathcal S_1}.
\]
Since $|\operatorname{Tr} (\mathscr M_\varepsilon A \mathscr M_\varepsilon)|\leq\|\mathscr M_\varepsilon A \mathscr M_\varepsilon\|_{\mathcal S_1}$, \eqref{20261031for10} gives
\[
 0\leq\int_{\{\delta\geq\varepsilon\}}
       \delta(\mathbf z)^{t-c}\,\mathrm d\mu(\mathbf z)\leq\|A\|_{\mathcal S_1}.
\]
As $\varepsilon\downarrow0$, the sets $\{\delta\geq\varepsilon\}$
increase to $B$. Monotone convergence therefore gives
\[
 \int_B\delta(\mathbf z)^{t-c}\,\mathrm d\mu(\mathbf z)\leq \|A\|_{\mathcal S_1}<\infty.
\]
By \eqref{eq:radialthreshold}, this forces $t-c>-1$, i.e. $c<t+1$.
\end{proof}

Assume $c<t+1$. Proposition \ref{20261031prop1} has already established that
$ T_{a,b,c}$ and $S_{a,b,c}$ are trace class in this case. Also
$\mathscr M_\varepsilon\to I$ strongly as $\varepsilon\downarrow0$.
For every trace-class operator $A$ on $L^2(B,\mu)$,
\begin{equation}\label{20261031for11}
 \|\mathscr M_\varepsilon A \mathscr M_\varepsilon-A\|_{\mathcal S_1}\longrightarrow0.
 \end{equation} To see this, first take $A$ finite rank: strong convergence of
$\mathscr M_\varepsilon$ implies convergence in trace norm on a finite
dimensional range, both for left and right multiplication. For
general trace-class $A$, approximate it in trace norm by finite-rank
operators and use $\|\mathscr M_\varepsilon\|_{\rm op}\leq1$ and the ideal
property. This proves \eqref{20261031for11}.

Continuity of the trace in trace norm now yields
\[
 \operatorname{Tr} A=\lim_{\varepsilon\downarrow0}
       \operatorname{Tr}(\mathscr M_\varepsilon A \mathscr M_\varepsilon).
\]
Using \eqref{20261031for10}  and monotone convergence,
\begin{equation}\label{20261031for12}
 \operatorname{Tr} A_T=\operatorname{Tr} A_S
 =\lim_{\varepsilon\downarrow0}
 \int_{\{\mathbf z\in B:\delta(\mathbf z)\geq\varepsilon\}}\delta(\mathbf z)^{a+b+\frac{\beta-\alpha}{2}-c}\,\mathrm d\mu(\mathbf z)
 =I_{a+b+\frac{\beta-\alpha}{2}-c}.
\end{equation}
Together with the necessity proved in Proposition \ref{20261031prop2}, the sufficiency in
Proposition \ref{20261031prop1} and \eqref{20261031for12}, this completes the proof of Conclusion 2 of Theorem \ref{thm:schatten2}.

We now prove the final assertion of Theorem~\ref{thm:schatten2} using techniques from Fourier analysis.

For every $n\in\mathbb N$, we denote all $C^{\infty}$ complex-valued functions on $\mathbb C^n$ with compact support by $C_c^{\infty}(\mathbb{C}^n)$. Since any $f\in C_c^{\infty}(\mathbb{C}^n)$ can be regarded as a cylinder function on $\ell^2$ that depends only on the first $n$ variables, we write
\[
\mathscr {C}_c^{\infty}\triangleq \bigcup_{k=1}^{\infty}C_c^{\infty}(\mathbb{C}^k),
\]
which can be viewed as the set of all smooth cylindrical functions on $\ell^2$. Let $\mathcal H\triangleq L^2(B,\mu),$ and
\begin{equation}
(V_\theta f)(\mathbf z)\triangleq f(e^{i\theta}\mathbf z),\quad \theta\in\mathbb R,\,\,f\in \mathscr C_c^{\infty},\,\,\mathbf z=(z_j)_{j\in\mathbb N}\in B,
\end{equation}
where
$
e^{i\theta}\mathbf z=(e^{i\theta}z_j)_{j\geq1}.
$
Invariance of $\mu$ gives
\[
\begin{aligned}
\|V_\theta f\|_{\mathcal H}^2
=\int_{B}|f(e^{i\theta}\mathbf z)|^2\,d\mu(\mathbf z)
=\int_{B}|f(\mathbf z)|^2\,\mathrm d\mu(\mathbf z)
=\|f\|_{\mathcal H}^2,\qquad f\in \mathscr C_c^{\infty}.
\end{aligned}
\]
By \cite[Proposition 2.4, p. 528]{YZ}, $\mathscr C_c^{\infty}$ is dense in $\mathcal H$, $V_{\theta}$ extends to a unitary operator on $\mathcal H$.

\begin{lemma}\label{20261101lem1}
\[
\lim_{\theta\to0}\|V_\theta f-f\|_{\mathcal H}=0,
\qquad f\in\mathcal H .
\]
\end{lemma}
\begin{proof}
For $f\in \mathscr C_c^{\infty}$,
$f(e^{i\theta}\mathbf z)\to f(\mathbf z)$ pointwise as $\theta\to0$, and
\[
|f(e^{i\theta}\mathbf z)-f(\mathbf z)|^2\leq4\left(\sup_{\mathbf z\in B} |f(\mathbf z)|\right)^2.
\]
The dominated convergence theorem therefore gives
\[
\|V_\theta f-f\|_{\mathcal H}\longrightarrow0.
\]
Since $\mathscr C_c^{\infty}$ is dense in $\mathcal H$, for a
general $f\in\mathcal H$, we can choose $g\in \mathscr C_c^{\infty}$ with
$\|f-g\|_{\mathcal H}<\varepsilon$. Since every $V_\theta$ is a
unitary operator,
\[
\begin{aligned}
\|V_\theta f-f\|_{\mathcal H}
\leq \|V_\theta(f-g)\|_{\mathcal H}
   +\|V_\theta g-g\|_{\mathcal H}+\|g-f\|_{\mathcal H}
\leq2\varepsilon+\|V_\theta g-g\|_{\mathcal H}.
\end{aligned}
\]
Letting $\theta\to0$ and then $\varepsilon\downarrow0$ proves strong
continuity.
\end{proof}
In particular, $\theta\mapsto V_\theta f$ is continuous
as an $\mathcal H$-valued function. Combining this with \cite[Theorem 1, p. 133]{Yo80} implies that the following integral exists in the sense of Bochner as in \cite[pp. 132--133]{Yo80},
\[
P_Nf\triangleq \frac1{2\pi}\int_0^{2\pi}
             e^{-iN\theta}V_\theta f\,d\theta,\qquad f\in \mathcal H,\,N\in\mathbb Z.
\]
Furthermore,
\begin{equation}\label{20261101for1}
V_\theta^*=V_{-\theta},\qquad V_{\theta+\phi}=V_\theta V_\phi,\qquad \theta,\phi\in\mathbb R,\qquad
V_0=I,\qquad V_{2\pi}=I,
\end{equation}
where $I$ denotes the identity operator on $\mathcal H$.


For $f,g\in\mathcal H$, \eqref{20261101for1} and changing variables give
\[
\begin{aligned}
\langle P_Nf, g\rangle_{\mathcal H}
&=\frac1{2\pi}\int_0^{2\pi}
                  \langle e^{-iN\theta}V_{ \theta}f,g\rangle_{\mathcal H}\,\mathrm d\theta\\
&=\frac1{2\pi}\int_0^{2\pi}
                  \langle f,e^{iN\theta}V_{-\theta}g\rangle_{\mathcal H}\,\mathrm d\theta\\
&=\frac1{2\pi}\int_0^{2\pi}
                  \langle f,e^{-iN\theta}V_{\theta}g\rangle_{\mathcal H}\,\mathrm d\theta
=\langle f, P_Ng\rangle_{\mathcal H}.
\end{aligned}
\]

For $M,N\in\mathbb Z$, Fubini's theorem and \eqref{20261101for1} yield
\[
\begin{aligned}
\langle P_NP_Mf, g \rangle_{\mathcal H}
&=\frac1{(2\pi)^2}\int_0^{2\pi}\int_0^{2\pi}
 e^{-iN\theta-iM\phi}\langle V_{\theta+\phi}f, g\rangle_{\mathcal H}\,\mathrm d\phi\,\mathrm d\theta.
\end{aligned}
\]
Set $\psi=\theta+\phi$ in the inner integral, with all variables
understood modulo $2\pi$. Then
\[
\begin{aligned}
\langle P_NP_Mf, g \rangle_{\mathcal H}
&=\frac1{2\pi}\int_0^{2\pi}e^{-iN\psi}\langle  V_\psi f, g \rangle_{\mathcal H}\,\mathrm d\psi
 \left(\frac1{2\pi}\int_0^{2\pi}
                    e^{i(N-M)\phi}\,\mathrm d\phi\right)\\
&=\begin{cases}
\langle P_Nf, g\rangle_{\mathcal H},&N=M,\\
0,&N\ne M.
\end{cases}
\end{aligned}
\]
Consequently,
\begin{equation}
P_N^*=P_N,\qquad P_N^2=P_N,\qquad
P_NP_M=0\quad(N\ne M).
\label{20261101for2}
\end{equation}
Thus the $\{P_N:N\in\mathbb Z\}$ are mutually orthogonal projections.
Define
\[
E_N\triangleq \{f\in\mathcal H:V_\theta f=e^{iN\theta}f
                   \text{ for every }\theta\in\mathbb R\},\qquad N\in\mathbb Z.
\]
\begin{proposition}\label{20261109prop1}It holds that
$
 \mathcal H=\bigoplus_{N\in\mathbb Z}E_N
$
and
$$
 \sum_{N\in\mathbb Z}P_N=I
       ,
 $$
in the strong operator topology.

\end{proposition}
\begin{proof}


For any $\phi\in\mathbb R,N\in\mathbb Z,$ and $f\in\mathcal H$, \eqref{20261101for1} gives
\[
\begin{aligned}
V_\phi P_Nf
=\frac1{2\pi}\int_0^{2\pi}
       e^{-iN\theta}V_{\theta+\phi}f\,\mathrm d\theta
=e^{iN\phi}\frac1{2\pi}\int_0^{2\pi}
       e^{-iN(\theta+\phi)}V_{\theta+\phi}f\,\dd\theta
=e^{iN\phi}P_Nf.
\end{aligned}
\]
Hence $\operatorname{Ran}(P_N)\subset E_N$. Conversely, if
$f\in E_N$, then
\[
P_Nf=\frac1{2\pi}\int_0^{2\pi}
e^{-iN\theta}e^{iN\theta}f\,\dd\theta=f\in \operatorname{Ran}(P_N).
\]
Therefore,
\begin{equation}
 \operatorname{Ran}(P_N)=E_N.                \label{20261101for3}
\end{equation}


For $m\in\mathbb N_0$, define the Fej\'er kernel
\[
F_m(\theta)\triangleq\frac1{m+1}
       \left|\sum_{j=0}^{m}e^{ij\theta}\right|^2,\qquad \theta\in\mathbb R.
\]
Its Fourier expansion is
\begin{equation}
F_m(\theta)=
\sum_{|N|\leq m}
\left(1-\frac{|N|}{m+1}\right)e^{iN\theta}.
\label{20261101for4}
\end{equation}
The following properties (one may see \cite[Chapter 2, $\S$ 5.2, pp. 52--54]{SS}) hold:
\begin{equation}
F_m(\theta)\geq0,\qquad
\frac1{2\pi}\int_{-\pi}^{\pi}F_m(\theta)\,d\theta=1,
\label{20261101for5}
\end{equation}
and, for each $0<\delta<\pi$,
\begin{equation}
\lim_{m\to\infty}
\frac1{2\pi}\int_{\delta\leq|\theta|\leq\pi}
                         F_m(\theta)\,d\theta=0.
\label{20261101for6}
\end{equation}
Indeed, for $\theta\not=0$ modulo $2\pi$,
$
F_m(\theta)=\frac1{m+1}
 \frac{\sin^2((m+1)\theta/2)}{\sin^2(\theta/2)}.
$
Thus, whenever $\delta\leq|\theta|\leq\pi$,
\[
0\leq F_m(\theta)\leq
\frac1{(m+1)\sin^2(\delta/2)},
\]
which proves \eqref{20261101for6}. The normalization in \eqref{20261101for5} follows from \eqref{20261101for4}.

For $f\in\mathcal H$, define its Fej\'er mean by
\[
\sigma_m(f)\triangleq\frac1{2\pi}\int_{-\pi}^{\pi}
                         F_m(\theta)V_\theta f\,\dd\theta,\qquad m\in\mathbb N.
\]
Using the finite Fourier expansion \eqref{20261101for4} and the definition of $P_N$,
we obtain
\begin{equation}
\sigma_m(f)=\sum_{|N|\leq m}
 \left(1-\frac{|N|}{m+1}\right)P_Nf.
\label{20261101for8}
\end{equation}

We claim that $\sigma_m(f)\to f$ in $\mathcal H$. Since the Fej\'er
kernel has normalized integral one,
\[
\sigma_m(f)-f
=\frac1{2\pi}\int_{-\pi}^{\pi}
 F_m(\theta)(V_\theta f-f)\,\dd\theta.
\]
The norm inequality for Bochner integrals (\cite[Corollary 1, p. 133]{Yo80}) gives
\begin{equation}
\|\sigma_m(f)-f\|_{\mathcal H}
\leq\frac1{2\pi}\int_{-\pi}^{\pi}
 F_m(\theta)\|V_\theta f-f\|_{\mathcal H}\,\dd\theta.
\label{20261101for9}
\end{equation}

Let $\varepsilon>0$. Lemma \ref{20261101lem1} provides
$0<\delta<\pi$ such that
\[
\|V_\theta f-f\|_{\mathcal H}<\varepsilon,
\qquad |\theta|<\delta .
\]
On the complementary set we use
\[
\|V_\theta f-f\|_{\mathcal H}\leq2\|f\|_{\mathcal H}.
\]
Splitting the integral in \eqref{20261101for9}, and using \eqref{20261101for5}, yields
\[
\begin{aligned}
\|\sigma_m(f)-f\|_{\mathcal H}
&\leq \varepsilon
+\frac{2\|f\|_{\mathcal H}}{2\pi}
  \int_{\delta\leq|\theta|\leq\pi}F_m(\theta)\,d\theta.
\end{aligned}
\]
By \eqref{20261101for6}, the second term tends to zero as $m\to\infty$.
Since $\varepsilon$ is arbitrary,
\[
\lim_{m\to\infty}\|\sigma_m(f)-f\|_{\mathcal H}=0.
\]

By \eqref{20261101for8}, each $\sigma_m(f)$ belongs to the finite linear span of
the spaces $E_N$. Since $\sigma_m(f)\to f$, it follows that
$
\overline{\operatorname{span}\{E_N:N\in\mathbb Z\}}
=\mathcal H.
$
Together with the mutual orthogonality established in \eqref{20261101for2} and
the range identity \eqref{20261101for3}, this proves
\[
 \mathcal H=\bigoplus_{N\in\mathbb Z}E_N,
\]
which implies that the
symmetric partial sums converge in norm:
\[
\lim_{m\to\infty}
\left\|\sum_{|N|\leq m}P_Nf-f\right\|_{\mathcal H}=0,
\qquad f\in\mathcal H.
\]
This is exactly the assertion
\[
\sum_{N\in\mathbb Z}P_N=I
\quad\text{strongly}.
\]
\end{proof}
\begin{remark}
ChatGPT pointed out that the above decomposition can be obtained by the Peter--Weyl Theorem(one may see \cite[Theorem 5.12, p. 143]{Fol}) in the representation theory of compact groups. For completeness, we still provide a short proof here which is based on basic techniques in Fourier analysis.
\end{remark}
Since both kernels in \eqref{eq:unitarykernels} are invariant
under $V_{\theta}$, $A_XV_\theta=V_\theta A_X$ for all $X=T,S$ and $\theta\in\mathbb R$.
Consequently each $E_N$ reduces $A_X$ and
\begin{equation}
 A_X=\bigoplus_{N\in\mathbb Z}P_NA_XP_N,\qquad X=T,S.   \label{202611011.2}
\end{equation}

For $n,m\geq0$, since $R_{n,m}$ have kernel
\[
 \delta(\mathbf z)^u\delta(\mathbf y)^v
 \langle\mathbf  z,\mathbf y\rangle^n\overline{\langle \mathbf z,\mathbf y\rangle}^{\,m}.
\]
Its rank-one expansion is
\[
 R_{n,m}=\sum_{\substack{|\eta|=n\\|\zeta|=m}}
 \frac{n!m!}{\eta!\zeta!}
 (\delta^u \mathbf z^\eta\overline{\mathbf z^\zeta})
 \otimes(\delta^v \mathbf z^\eta\overline{\mathbf z^\zeta}).
\]
Hence $R_{n,m}$ annihilates $E_N$ unless
$N=n-m$, and its range is contained in $E_{n-m}$. In particular,
$R_{n,m}=P_{n-m}R_{n,m}P_{n-m}$. Since
$A_T=\sum_{n\geq0}A_n(c)R_{n,0}$,
\begin{equation}
 P_NA_TP_N=A_N(c)R_{N,0},\quad N\geq0,\qquad
 P_NA_TP_N=0,\quad N<0.                                \label{20261101for1.4}
\end{equation}

\begin{lemma}\label{20261113lem3}
 For both $X=T,S$ and every $D>0$,
\begin{equation*}\label{eq:blockrapid*}
\|P_NA_XP_N\|_{\mathrm{op}}
       \leq C_D(1+|N|)^{-D},\qquad N\in\mathbb Z,
\end{equation*}
where $C_D$ is a positive constant depending only on \(D\) and the fixed parameters \(u,v,c\), as well as the underlying Gaussian measure, but independent of \(N\).
\end{lemma}
\begin{proof}
Recall \eqref{eq:rapid-Sr}, for every $L>0$,
\begin{equation}
 \|R_{n,m}\|_{\text{op}}\leq\|R_{n,m}\|_{\mathcal S_2}
 \leq C_L(1+n+m)^{-L},\qquad n,m\in\mathbb N_0.                                 \label{202611011.4a}
\end{equation}
By Corollary \ref{20261113cor1}, the coefficients $A_N(c)$ have polynomial growth. Thus \eqref{20261101for1.4}
immediately gives
\[
 \|P_NA_TP_N\|_{\text{op}}
 \leq C_D(1+N)^{-D},\qquad N\geq0,
\]
for every $D>0$.

For $A_S$, its expansion is
\[
 A_S=\sum_{n,m\geq0}A_n(c/2)A_m(c/2)R_{n,m}.
\]
Since $R_{n,m}=P_{n-m}R_{n,m}P_{n-m}$, the $N$th block is
\[
 P_NA_SP_N
 =\sum_{\substack{n,m\geq0\\n-m=N}}
 A_n(c/2)A_m(c/2)R_{n,m}.
\]
By Corollary \ref{20261113cor1}, so for some
$M\geq0$,
\[
 |A_n(c/2)A_m(c/2)|\leq C(1+n+m)^M,\qquad n,m\in\mathbb N_0.
\]
Consequently, by \eqref{202611011.2} and \eqref{202611011.4a},
\[
 \|P_NA_SP_N\|_{\text{op}}
 \leq C_L\sum_{n-m=N}(1+n+m)^{M-L}.
\]
If $N\geq0$, write $n=m+N$; if $N<0$, write $m=n+|N|$.
In either case the right-hand side of above inequality is bounded by
\[
 C_L\sum_{q=0}^\infty(1+|N|+2q)^{M-L}.
\]
For $L>M+D+1$, comparison with an integral gives
\[
 \sum_{q=0}^\infty(1+|N|+2q)^{M-L}
 \leq C_{D,L}(1+|N|)^{-D}.
\]
Thus, for both $X=T,S$ and every $D>0$,
$
\|P_NA_XP_N\|_{\text{op}}
       \leq C_D(1+|N|)^{-D}.\
$
\end{proof}

We use the following elementary singular-value inequality. If $B\in\mathcal S_r$
and $0<r<1$, then
\begin{equation}
 \|B\|_{\mathcal S_1}\leq
 \|B\|_{\text{op}}^{\,1-r}\|B\|_{\mathcal S_r}^{\,r}.                 \label{202611011.5}
\end{equation}
Indeed, if $s_j(B)$ are the singular values, then
$s_j(B)\leq\|B\|_{\text{op}
}$ and
\[
 \sum_js_j(B)
 =\sum_js_j(B)^r s_j(B)^{1-r}
 \leq\|B\|_{\text{op}}^{1-r}\sum_js_j(B)^r.
\]
The same argument applies to finite-rank compressions, and then by
monotone convergence to compact operators.

\begin{lemma} \label{20261101lem2}
Suppose $A$ is a bounded operator from
$L^2(B,\mu)$ into itself such that $ A =\bigoplus_{N\in\mathbb Z}P_NA P_N$, and for every $D>0$,
\begin{equation}
 \|P_NAP_N\|_{\mathrm{op}}\leq C_D(1+|N|)^{-D}.                \label{202611011.6}
\end{equation}
If $A\in\mathcal S_r$ for some $0<r<1$, then $A\in\mathcal S_1$. Moreover,
\[
 N\longmapsto\text{Tr}(P_NAP_N)
\]
decreases faster than every negative power of $|N|$, and
\begin{equation}
 \theta\longmapsto\text{Tr}(V_\theta A)
 =\sum_{N\in\mathbb Z}e^{iN\theta}\text{Tr}(P_NAP_N)
                                                               \label{202611011.7}
\end{equation}
is a $C^\infty$ periodic function.
\end{lemma}

\begin{proof}
Since $ A =\bigoplus_{N\in\mathbb Z}P_NA P_N$,
\[
 \|P_NAP_N\|_{\mathcal S_r}\leq\|A\|_{\mathcal S_r},\qquad N\in\mathbb Z.
\]
Applying  \eqref{202611011.5} to $B_N\triangleq P_NAP_N$ and using \eqref{202611011.6}, we get
\begin{equation}
 \|B_N\|_{\mathcal S_1}
 \leq C_D^{1-r}(1+|N|)^{-D(1-r)}
             \|A\|_{\mathcal S_r}^{\,r},\qquad N\in\mathbb Z.                       \label{202611011.8}
\end{equation}
Choose $D$ so large that $D(1-r)>1$. Summing \eqref{202611011.8} over $N$
shows $\sum_N\|B_N\|_{\mathcal S_1}<\infty$. The operator $A$
is therefore trace class and
\[
 \|A\|_{\mathcal S_1}=\sum_N\|B_N\|_{\mathcal S_1}.
\]
For each $N$,
\[
 |\text{Tr} B_N|\leq\|B_N\|_{\mathcal S_1}
 \leq C_D^{1-r}(1+|N|)^{-D(1-r)}\|A\|_{\mathcal S_r}^{\,r}.
\]
Since $D$ in \eqref{202611011.6} is arbitrary, the block traces have rapid decay.

On $E_N$, $V_\theta$ is multiplication by $e^{iN\theta}$.
The trace-class block decomposition thus gives  \eqref{202611011.7}, with absolute
and uniform convergence. For every integer $k\geq0$, choose $D$
large enough that
$\sum_N|N|^k|\text{Tr} B_N|<\infty$. The differentiated Fourier series
then converges uniformly. Hence the function in \eqref{202611011.7} is $C^\infty$.
\end{proof}

\begin{lemma}
If $A_S$ is trace class, then
\begin{equation}
 G(\theta)\triangleq\mathrm{Tr}(V_\theta A_S)
 =\int_0^1(1-x)^t \widetilde{h}(x)|1-xe^{i\theta}|^{-c}\,\dd x,\label{202611012.5}
\end{equation}
where $\widetilde{h}$ is the density of $X$ which appeared in Theorem \ref{lem:moments}.
\end{lemma}

\begin{proof}
The trace-class homogeneous expansion in the proof of Proposition
\ref{20261031prop2} permits
termwise tracing after first using the truncations
$\mathscr M_\varepsilon$ and then letting $\varepsilon\downarrow0$ in trace
norm. The operator $V_\theta A_S$ has kernel obtained by replacing
the output variable $\mathbf z$ in the kernel of $A_S$ by $e^{i\theta}\mathbf z$.
On the diagonal $\mathbf y=\mathbf z$ its kernel is
\[
 \delta(\mathbf z)^{u+v}
 \left|1-\langle e^{i\theta}\mathbf z,\mathbf z\rangle\right|^{-c}
 =\delta(\mathbf z)^t\left|1-e^{i\theta}\|\mathbf z\|^2\right|^{-c}.
\]
The trace is the integral of this diagonal value; this identity is
justified by the
finite sums followed by trace-norm convergence in the proof of Proposition
\ref{20261031prop2}. Precisely, let
$$
A_S^{N}\triangleq \sum_{n,m\geq 0, n+m\leq N}A_n(c/2)A_m(c/2)R_{n,m},\qquad N\in\mathbb N.
$$
Since \(A_S\) is trace class,
Proposition~\ref{20261031prop2} gives \(c<t+1\). By the proof of Proposition~\ref{20261031prop1},
$$
 \|A_S^{N}-A_S\|_{\mathcal S_1}\to 0,
$$
and noting that $|\text{Tr}(V_{\theta}(A_S^{N}-A_S))|\leq \|A_S^{N}-A_S\|_{\mathcal S_1}$. Hence,
$$
\lim_{N\to \infty}\text{Tr}(V_{\theta}A_S^{N})=\text{Tr}(V_{\theta}A_S).
$$
For each finite truncation, the trace-norm convergent rank-one expansions of the operators \(R_{n,m}\), established in Lemma \ref{lem:pieces}, give
\[
\operatorname{Tr}(V_\theta A_S^N)
=
\int_0^1(1-x)^t\widetilde h(x)
\sum_{\substack{n,m\geq0\\n+m\leq N}}
A_n(c/2)A_m(c/2)e^{i(n-m)\theta}x^{n+m}\,\mathrm dx.
\]
If \(c>0\), the coefficients \(A_n(c/2)\) are nonnegative, and the absolute value of the summand is dominated by
\[
(1-x)^t\widetilde h(x)
\sum_{n,m\geq0}A_n(c/2)A_m(c/2)x^{n+m}
=
(1-x)^{t-c}\widetilde h(x),
\]
which is integrable because \(c<t+1\).
If \(c<0\), Lemma~\ref{20261112lem1} gives
\[
\sum_{n,m\geq0}
|A_n(c/2)A_m(c/2)|
=
\Lambda(c/2)^2<\infty.
\]
Consequently, the integrand is dominated by
\[
\Lambda(c/2)^2(1-x)^t\widetilde h(x),
\]
which is integrable because \(t>-1\) and \(\widetilde h\) is bounded. The case \(c=0\) is immediate. Dominated convergence therefore yields the asserted integral formula.
Finally, the
law of $X$ has density $\widetilde{h}$, so changing variables
$x=\|\mathbf z\|^2$ gives \eqref{202611012.5}.
\end{proof}
We need the following lemma in the proof of Proposition \ref{prop:nonterminating}.
\begin{lemma}
\label{lem:even-smooth}
Let \(G\in C^\infty((-a,a))\), where \(a>0\), and suppose that \(G\) is even:
$
G(-x)=G(x),\, |x|<a.
$
Define
$
H(s)\triangleq G(\sqrt{s}),\, 0\leq s<a^2.
$
Then \(H\) belongs to \(C^\infty([0,a^2))\), where smoothness at
\(s=0\) is understood in the one-sided sense. More precisely, for every
integer \(k\geq 0\), there is an even function
\(G_k\in C^\infty((-a,a))\) such that
$
H^{(k)}(s)=G_k(\sqrt{s}),\, 0<s<a^2,
$
and consequently
$
H^{(k)}(0)=G_k(0).
$
\end{lemma}

\begin{proof}
We construct the functions \(G_k\) recursively. Set \(G_0\triangleq G\).
Suppose that \(G_k\) is smooth and even. Since \(G_k\) is even,
\(G_k'(0)=0\), and \(G_k'\) is odd. For \(x\neq 0\), define
\[
G_{k+1}(x)\triangleq \frac{G_k'(x)}{2x}.
\]
At \(x=0\), define
$
G_{k+1}(0)\triangleq \frac12G_k''(0).
$
We first verify that \(G_{k+1}\) is smooth and even. By the fundamental
theorem of calculus and \(G_k'(0)=0\),
\[
G_k'(x)
=x\int_0^1G_k''(tx)\,\mathrm dt.
\]
Thus, for every \(x\in(-a,a)\),
\[
G_{k+1}(x)
=\frac12\int_0^1G_k''(tx)\,\mathrm dt.
\]
The right-hand side defines a smooth function of \(x\), since
differentiation under the integral is valid on compact subintervals.
Moreover, \(G_k''\) is even, so \(G_{k+1}\) is even. Its value at zero
is exactly \(\frac12G_k''(0)\), as prescribed.

For \(s>0\), the chain rule gives
\[
\frac{\mathrm d}{\mathrm ds}G_k(\sqrt{s})
=\frac{G_k'(\sqrt{s})}{2\sqrt{s}}
=G_{k+1}(\sqrt{s}).
\]
Starting with \(H(s)=G_0(\sqrt{s})\), induction therefore yields
\[
H^{(k)}(s)=G_k(\sqrt{s}),\qquad s>0.
\]
Since \(G_k\) is continuous at zero,
\[
\lim_{s\downarrow0}H^{(k)}(s)=G_k(0).
\]
These continuous extensions are indeed the successive right
derivatives at zero. To see this inductively, if \(H^{(k)}\) on
\((0,a^2)\) has a continuous extension to zero, then the fundamental
theorem of calculus gives
\[
H^{(k-1)}(s)-H^{(k-1)}(0)
=\int_0^s H^{(k)}(t)\,\mathrm dt,
\]
where \(H^{(k-1)}(0)\) is its continuous limiting value. Dividing by
\(s\) and letting \(s\downarrow0\) shows that the right derivative at
zero equals \(H^{(k)}(0)\). This proves the assertion for every \(k\).
\end{proof}
\begin{proposition}
\label{prop:nonterminating}
For \(0<r<1\),
\[
 c\notin\{0,-1,-2,\ldots\}\ \Longrightarrow\ A_T\notin\mathcal S_r,
 \qquad
 c\notin\{0,-2,-4,\ldots\}\ \Longrightarrow\ A_S\notin\mathcal S_r.
\]
\end{proposition}

\begin{proof}
For \(T\),
\(P_NA_TP_N=A_N(c)R_{N,0},\,N\ge0\). By \eqref{eq:Jasymptotic}, \eqref{eq:piecetracebound},
\eqref{eq:negative-asymptotic235}, for $c\notin\{0,-1,-2,\ldots\}$,
\[
| \operatorname{Tr}(P_NA_TP_N)|
 =|A_N(c)J_N(t)|
 \sim \frac{\widetilde{h}(1)\Gamma(t+1)}{|\Gamma(c)|}N^{c-t-2},\qquad \text{as}\,N\to\infty.
\]
The coefficient is nonzero. These traces do not have rapid decay,
combining Lemmas~\ref{20261113lem3} and \ref{20261101lem2} imply $A_T\notin\mathcal S_r$.

For \(S\), suppose to the contrary that \(A_S\in\mathcal S_r\). Then Lemmas~\ref{20261113lem3} and \ref{20261101lem2} imply $A_S\in\mathcal S_1$ and the conclusion 2 of Theorem \ref{thm:schatten2} that has already been proved gives \(c<t+1\). By Lemma \ref{20261101lem2}, the function $G$ defined in \eqref{202611012.5} is an even $C^{\infty}$ function. Write \(d=c/2\) and
\begin{equation}\label{eq:Hphase}
 w= 2(1-\cos\theta),\qquad
 F(w)=\int_0^1(1-x)^t \widetilde{h}(x)
       \bigl((1-x)^2+xw\bigr)^{-d}\dd x,\qquad w>0.
\end{equation}
Since
\[
 |1-xe^{i\theta}|^2=(1-x)^2+2x(1-\cos\theta),
\]
equation \eqref{202611012.5} becomes
\begin{equation}
 G(\theta)=F(2(1-\cos\theta)).                          \label{202611012.7}
\end{equation}
For each fixed $w>0$, the quantity $(1-x)^2+xw$ is bounded away
from zero on $[0,1]$. Differentiation under the integral is therefore
legitimate, and for every integer $k\geq0$,
\begin{equation}
 F^{(k)}(w)=(-1)^k(d)_k
 \int_0^1(1-x)^t \widetilde{h}(x)x^k
       \big((1-x)^2+xw\big)^{-d-k}\,\dd x.                 \label{202611012.8}
\end{equation}
For clarity, on a compact interval $w\in[w_0,w_1]\subset(0,\infty)$
the integrand and each of its $w$-derivatives are bounded by a
constant times $(1-x)^t \widetilde{h}(x)$; this is integrable because $t>-1$
and $h$ is bounded. Thus dominated convergence proves \eqref{202611012.8}.

The hypothesis $c\notin\{0,-2,-4,\ldots\}$ is equivalent to
$d\notin\{0,-1,-2,\ldots\}$. Hence $(d)_k\neq0$ for every $k$.
Choose $k$ so large that
\begin{equation}
 d+k>0,\qquad 2(d+k)\geq t+1.                           \label{202611012.9}
\end{equation}
This is possible because $k$ can be arbitrarily large. As
$w\downarrow0$, the positive integrand in \eqref{202611012.8} increases pointwise
to
\[
 (1-x)^{t-2(d+k)}\widetilde{h}(x)x^k.
\]
By Lemma \ref{lem:radialdensity}, $\widetilde{h}$ is bounded below by a positive constant on a
left-neighborhood of $1$. Also $x^k$ is bounded below there by a
positive constant. The limiting integral therefore dominates a
positive multiple of
\[
 \int_{x_0}^1(1-x)^{t-2(d+k)}\,\dd x,\qquad x_0\in(0,1),
\]
which diverges by \eqref{202611012.9}, since the exponent is at most $-1$.
Monotone convergence gives
\begin{equation}
 |F^{(k)}(w)|\longrightarrow\infty,
 \qquad w\downarrow0.                                 \label{202611012.10}
\end{equation}

We now explain why \eqref{202611012.10} contradicts smoothness of $G$.
The function $G$ is even, because the right-hand side of \eqref{202611012.5}
depends on $\theta$ only through $\cos\theta$. If a $C^\infty$
function $G$ is even, then by Lemma \ref{lem:even-smooth},
\[
 H(s)\triangleq G(\sqrt{s}),\qquad s\geq0,
\]
is $C^\infty$ up to $s=0$. Since
\[
 2(1-\cos\sqrt{s})=s+O(s^2)
\]
and its derivative at $s=0$ is $1$, the map
$s\mapsto2(1-\cos\sqrt{s})$ is a smooth local diffeomorphism of
$[0,\infty)$ near zero. Indeed, it admits a smooth extension to an open interval containing zero through its power series. Equation \eqref{202611012.7} therefore implies that
$F\in C^{\infty}([0,\varepsilon))$ for sufficiently small $\varepsilon>0$. In particular $F^{(k)}$ must remain
bounded near zero, contradicting \eqref{202611012.10}. This proves Proposition
\ref{prop:nonterminating}.
\end{proof}


\begin{proposition}\label{prop:polynomial}
Fix \(0<r<1\) and an integer \(k\ge1\). Then
\[
 T_{a,b,-k}\in\mathcal S_r
 \quad\Longleftrightarrow\quad\sum_j\vartheta_j^r<\infty,
 \qquad
 S_{a,b,-2k}\in\mathcal S_r
 \quad\Longleftrightarrow\quad\sum_j\vartheta_j^r<\infty.
\]
Both operators have infinite rank.
\end{proposition}
\begin{proof}
Let $c=-k$. Then
\[
 (1-\langle \mathbf z,\mathbf y\rangle)^k
 =\sum_{n=0}^k A_n(-k)\langle \mathbf z,\mathbf y\rangle^n,\qquad
 A_n(-k)=(-1)^n\binom{k}{n}.
\]
By Proposition \ref{thm:svd}, the nonzero singular values of $A_T$ are
\[
 \sigma_\eta=
 \frac{|(-k)_{|\eta|}|}{\eta!}
 \bigl(M_\eta(2u)M_\eta(2v)\bigr)^{1/2},
 \qquad |\eta|\leq k.
\]
Using \eqref{eq:momentupper},
\[
 \sigma_\eta\leq C_k\vartheta^\eta,\qquad |\eta|\leq k.
\]
Thus
\[
 \|A_T\|_{\mathcal S_r}^r
 \leq C_k\sum_{n=0}^k\sum_{|\eta|=n}
                      (\vartheta^\eta)^r
 \leq C_k\sum_{n=0}^k
              \left(\sum_j\vartheta_j^r\right)^n<\infty,
\]
whenever $\sum_j\vartheta_j^r<\infty$. Conversely, for $\eta=\mathbf e_j$, \eqref{eq:firstmoment} implies
\begin{equation}
 \sigma_{\mathbf e_j}
 =k\bigl(M_{\mathbf e_j}(2u)M_{\mathbf e_j}(2v)\bigr)^{1/2}
 \sim k(I_{2u}I_{2v})^{1/2}\vartheta_j,\qquad \text{as}\,\,j\to\infty.                \label{202611013.3}
\end{equation}
If $A_T\in\mathcal S_r$, the sum of $\sigma_{\mathbf e_j}^r$ is finite, forcing
$\sum_j\vartheta_j^r<\infty$. Therefore
\[
 A_T\in\mathcal S_r\quad\Longleftrightarrow\quad
 \sum_j\vartheta_j^r<\infty,\qquad c=-k,\ k\geq1.
\]
Since every $\vartheta_j>0$, \eqref{202611013.3} also proves that $A_T$ has
infinite rank.

Let $c=-2k$. Then
\[
 |1-\langle \mathbf z,\mathbf y\rangle|^{2k}
 =(1-\langle \mathbf z,\mathbf y\rangle)^k
  (1-\overline{\langle \mathbf z,\mathbf y\rangle})^k.
\]
The expansion is finite in $n,m\leq k$. The rank-one terms of $A_S$, indexed
by $|\eta|=n$, $|\zeta|=m$, have vectors
\[
 f_{\eta,\zeta}=\delta^u \mathbf z^\eta\overline{\mathbf z^\zeta},
 \qquad
 g_{\eta,\zeta}=\delta^v \mathbf z^\eta\overline{\mathbf z^\zeta},
\]
and scalar coefficient
$A_n(-k)A_m(-k)n!m!/(\eta!\zeta!)$.
By \eqref{eq:momentupper}, their rank-one Schatten quasi-norm is
\begin{equation}
 \|f_{\eta,\zeta}\otimes g_{\eta,\zeta}\|_{\mathcal S_r}
 =\bigl(M_{\eta+\zeta}(2u)M_{\eta+\zeta}(2v)\bigr)^{1/2}
 \leq C(\eta+\zeta)!\vartheta^{\eta+\zeta}.             \label{202611013.5}
\end{equation}
Furthermore,
\begin{equation}
 n!m!\frac{(\eta+\zeta)!}{\eta!\zeta!}
 =n!m!\prod_j\binom{\eta_j+\zeta_j}{\eta_j}
 \leq n!m!\binom{n+m}{n}=(n+m)!\leq(2k)!.
                                                               \label{202611013.6}
\end{equation}
The product inequality follows by comparing the coefficient of
$x^n$ in $$\prod_j(1+x)^{\eta_j+\zeta_j}=(1+x)^{n+m}$$ and observing that the product of binomial coefficients is one particular nonnegative summand of that coefficient.

For $0<r<1$, assume \(\sum_{j=1}^{\infty}\vartheta_j^r<\infty\).
To justify the use of the Schatten quasi-triangle inequality for the
countable monomial expansions, introduce the truncations
\[
 A_{S,L}\triangleq
 \sum_{n,m=0}^k A_n(-k)A_m(-k)R_{n,m}^{L},\qquad L\in\mathbb N,
\]
where \(R_{n,m}^{L}\) is the finite-rank operator defined in the proof
of Lemma~\ref{lem:pieces}. For each \(L\), its rank-one expansion is finite.
Applying the Schatten \(r\)-triangle inequality and
\eqref{202611013.5}--\eqref{202611013.6}, we obtain
\[
 \|A_{S,L}\|_{\mathcal S_r}^r
 \leq C_k\sum_{n,m=0}^k
       \left(\sum_{j=1}^{L}\vartheta_j^r\right)^{n+m}
 \leq C_k\sum_{n,m=0}^k
       \left(\sum_{j=1}^{\infty}\vartheta_j^r\right)^{n+m}<\infty.
\]
The constant is independent of \(L\).
The proof of Lemma~\ref{lem:pieces} gives \(R_{n,m}^{L}\to R_{n,m}\) in trace norm;
therefore \(A_{S,L}\to A_S\) in operator norm. For every fixed \(j\),
the singular-value inequality
\[
 |s_j(A_{S,L})-s_j(A_S)|\leq\|A_{S,L}-A_S\|_{\mathrm{op}}
\]
implies \(s_j(A_{S,L})\to s_j(A_S)\). Fatou's lemma for series now yields
\[
 \sum_{j=1}^{\infty}s_j(A_S)^r
 \leq\liminf_{L\to\infty}\|A_{S,L}\|_{\mathcal S_r}^r<\infty.
\]
Consequently,
\begin{equation}\label{202611013.7}
 \sum_j\vartheta_j^r<\infty
 \quad\Longrightarrow\quad A_S\in\mathcal S_r,
 \qquad c=-2k.
\end{equation}

For each $j\in\mathbb N$, consider the unitary action on $L^2(B,\mu)$:
\[
 (V_\theta^{(j)}f)(\mathbf z)
\triangleq f(z_1,\ldots,e^{i\theta}z_j,\ldots),\qquad \theta\in\mathbb R,\,\mathbf z=(z_i)_{i\in\mathbb N}\in B,\,f\in L^2(B,\mu).
\]
Let $Q_j$ be the orthogonal projection on $L^2(B,\mu)$ onto the following subspace
$$
\left\{f\in L^2(B,\mu):V_\theta^{(j)}f=e^{i\theta}f,V_\theta^{(k)}f= f,\,\text{ for all }k\in\mathbb N\setminus\{j\},\,\theta\in\mathbb R\right\}.
$$
The projections $\{Q_j\}_{j=1}^{\infty}$ are mutually
orthogonal and every $Q_j$ reduces $A_S$.

For $c=-2k$, the finite expansion is
\[
 A_S=\sum_{n,m=0}^kA_n(-k)A_m(-k)R_{n,m}=\sum_{n,m=0}^kA_n(-k)A_m(-k)\sum_{\substack{|\eta|=n\\|\zeta|=m}}
 \frac{n!m!}{\eta!\zeta!}
 (\delta^u \mathbf z^\eta\overline{\mathbf z^\zeta})
 \otimes(\delta^v \mathbf z^\eta\overline{\mathbf z^\zeta}),
\]
where $A_n(-k)=(-1)^n\binom{k}{n}$. A term in the $Q_j$-block must have $\eta-\zeta=\mathbf e_j$.
Consequently $n=|\eta|=|\zeta|+1=m+1$. Thus $n+m=2m+1$
is odd, and every coefficient in this block has the same negative
sign:
\[
 A_n(-k)A_m(-k)
 =-\binom{k}{n}\binom{k}{m}<0.
\]
Use the unit vectors
\[
 f_j\triangleq\frac{\delta^v z_j}{M_{\mathbf e_j}(2v)^{1/2}},\qquad
 g_j\triangleq\frac{\delta^u z_j}{M_{\mathbf e_j}(2u)^{1/2}}.
\]
Direct computation shows that
$$
\ip{Q_jA_SQ_jf_j}{g_j}
 =-\frac1{D_j}\sum_{m=0}^{k-1}\ \sum_{|\zeta|=m}
      b_{k,m,\zeta}
      M_{\zeta+e_j}(2u)M_{\zeta+e_j}(2v),
$$
where
$$
b_{k,m,\zeta}
\triangleq \binom{k}{m+1}\binom{k}{m}
     \frac{(m+1)!m!}{(\zeta+\mathbf e_j)!\zeta!}>0,\qquad \text{and}\qquad D_j\triangleq \bigl(M_{\mathbf e_j}(2u)M_{\mathbf e_j}(2v)\bigr)^{1/2}>0.
$$
Since for $m=0$ and $\zeta=0$,
$
 b_{k,0,0}=\binom{k}{1}\binom{k}{0}\frac{1!0!}{\mathbf e_j!0!}=k,
$
\begin{equation*}
 \begin{split}
 -\ip{Q_jA_SQ_jf_j}{g_j}
 &=kD_j+\mathcal R_j,\\
 \mathcal R_j
 &\triangleq\frac1{D_j}\sum_{m=1}^{k-1}\ \sum_{|\zeta|=m}
       b_{k,m,\zeta}
       M_{\zeta+\mathbf e_j}(2u)M_{\zeta+\mathbf e_j}(2v)\geq0.
 \end{split}
\end{equation*}
Hence
$$
 |\langle Q_jA_SQ_jf_j,g_j \rangle_{L^2(B,\mu)}| \geq kD_j= k\bigl(M_{\mathbf e_j}(2u)M_{\mathbf e_j}(2v)\bigr)^{1/2},
$$
so
\begin{equation}\label{202611013.10}
 \begin{split}
 \|Q_jA_SQ_j\|_{\text{op}}&\geq  \|Q_jA_SQ_jf_j\|_{L^2(B,\mu)} \geq   |\langle Q_jA_SQ_jf_j,g_j \rangle_{L^2(B,\mu)}|\\
& \geq k\bigl(M_{\mathbf e_j}(2u)M_{\mathbf e_j}(2v)\bigr)^{1/2}
 \sim k(I_{2u}I_{2v})^{1/2}\vartheta_j.
\end{split}
\end{equation}
If $A_S\in\mathcal S_r$, orthogonality and reduction imply
\[
 \sum_j\|Q_jA_SQ_j\|_{\mathcal S_r}^r\leq\|A_S\|_{\mathcal S_r}^r.
\]
Since $\|Q_jA_SQ_j\|_{\mathcal S_r}\geq\|Q_jA_SQ_j\|_{\text{op}}$, \eqref{202611013.10} forces
$\sum_j\vartheta_j^r<\infty$. Combined with \eqref{202611013.7}, this proves
\[
 A_S\in\mathcal S_r\quad\Longleftrightarrow\quad
 \sum_j\vartheta_j^r<\infty,\qquad c=-2k,\ k\geq1.
\]
The lower bound \eqref{202611013.10} also proves that $A_S$ has infinite rank.

\end{proof}

If $c=0$, both kernels equal $\delta(\mathbf z)^u\delta(\mathbf y)^v$. Therefore
\[
 A_T=A_S=\delta^u\otimes\delta^v,
\]
is a rank-one operator whose only nonzero singular value is
\[
 \|\delta^u\|_{L^2(B,\mu)}\|\delta^v\|_{L^2(B,\mu)}
 =(I_{2u}I_{2v})^{1/2}.
\]
It belongs to every Schatten class $\mathcal S_r$, $r>0$. Combining this with Proposition \ref{prop:nonterminating} and Proposition \ref{prop:polynomial} proves the final conclusion of Theorem \ref{thm:schatten2}.

\section*{Acknowledgements}
This work is partially supported by National Key R$\&$D Program of China under grant 2024YFA1013400, and by the Science Development Project of Sichuan University under grant 2020SCUNL201.
\section*{AI Declaration}
The authors gratefully acknowledge ChatGPT (GPT-6 Astra Pro) 
for improving the language and
grammatical presentation of the manuscript and suggesting an approach to the proof of the final
assertion of Theorem~\ref{thm:schatten2}.
The authors independently developed the mathematical ideas and proofs for all other parts of the paper.

\end{document}